\documentclass[11pt , a4paper, fleqn]{article}
\usepackage [utf8] {inputenc} 
\usepackage[T1]{fontenc}
\usepackage {amsmath,mathtools}
\usepackage {hyperref} 
\usepackage {graphicx}
\usepackage {fancyhdr}
\usepackage {float} 
\usepackage {xcolor}
\usepackage{lastpage}
\usepackage[top=2cm, bottom=2cm, left=3cm , right=3cm]{geometry}
\usepackage{vmargin}
\usepackage{listings}
\usepackage{lscape}
\usepackage{amsthm, amssymb}
\usepackage{enumerate}
\usepackage[12pt]{extsizes}
\usepackage{upgreek}
\usepackage{authblk}
 \usepackage[toc,page]{appendix}

\usepackage{graphicx}
\usepackage{subcaption}
\usepackage{algorithmic,algorithm}

\def\build#1_#2^#3{\mathrel{\mathop{\kern 0 pt#1}\limits_{#2}^{#3}}}

{ \begin{list}%
	{$\bullet$}%
	{\setlength{\labelwidth}{30pt}%
	 \setlength{\leftmargin}{35pt}%
	 \setlength{\itemsep}{\parsep}}}%
{ \end{list} }

\def\build#1_#2^#3{\mathrel{\mathop{\kern 0 pt#1}\limits_{#2}^{#3}}}

\definecolor{hellgelb}{rgb}{1,1,0.8}
\definecolor{colKeys}{rgb}{0,0,1}
\definecolor{colIdentifier}{rgb}{0,0,0}
\definecolor{colComments}{rgb}{0,0.5,0}
\definecolor{colString}{rgb}{0.62,0.12,0.94}

\definecolor{orange}{rgb}{1,0.5,0}

\theoremstyle{definition}
\newtheorem{definition}{Definition }[section]
\newtheorem{example}{Example }[section]
\newtheorem{remark}{Remark}[section]
\newtheorem{assumption}{Assumption}[section]
\newtheorem{notation}{Notation}[section]

\theoremstyle{plain}
\newtheorem{proposition}{Proposition }[section]
\newtheorem{theorem}{Theorem }[section]

\newcommand{\pfrac}[2]{\frac{\partial#1}{\partial #2}}
\newcommand{\dleq}[3]{#1 \leq #2 \leq #3}

\def\build#1_#2^#3{\mathrel{\mathop{\kern 0 pt#1}\limits_{#2}^{#3}}}

\newcommand{\mspan}{\mathrm{span}\,}
\newcommand{\mrk}{\mathrm{rk}\,}

\newcommand{\mcork}{\mathrm{cork}\,}
\newcommand{\mmod}{\ \mathrm{mod}\,}

\newcommand{\md}{\mathrm{d}}

\def \vp {\varphi}
\def \a {\alpha}
\def \b {\beta}

\def \S {\Sigma}

\def \R {\mathbb{R}}

\def \t {\tilde}

\newcommand{\red}[1]{\textcolor{black}{#1}}

\newcommand{\green}[1]{\textcolor{black}{#1}}
\newcommand{\cyan}[1]{\textcolor{black}{#1}}

\newcommand{\nf}[1]{{\normalfont #1}}

\setmarginsrb{2.5cm}{2.5cm}{2cm}{2cm}{0cm}{0cm}{1cm}{1cm}

\fancypagestyle{plain}{
}

\date{}

\title{Dynamic feedback linearization of two-input control systems via successive one-fold prolongations}

\author[]{Florentina Nicolau
 \thanks{QUARTZ EA7393 Laboratory, ENSEA, 6 Avenue du Ponceau, 95014 Cergy-Pontoise, France.
         {\tt\small florentina.nicolau@ensea.fr}
        }
}

\author[]{Witold Respondek
\thanks{Institute of Automatic Control, Lodz University of Technology, \L{}\'od\'z, Poland.  Emeritus Professor at Normandie Universit\'e, INSA Rouen,
France.
        {\tt\small witold.respondek@insa-rouen.fr}}%
}

 \author[]{Shunjie Li
 \thanks{School of Mathematics and Statistics,
Nanjing University of Information Science and Technology, China.
     {\tt\small shunjie.li@nuist.edu.cn }}
 }

\affil[]{}
\begin{document}

\maketitle

%

\begin{abstract}
 In this paper, we propose a constructive algorithm to dynamically linearize two-input control systems via successive one-fold prolongations of a control that has to be suitably chosen at each step of the algorithm. Linearization via successive one-fold prolongations  requires special properties of the linearizability distributions $\mathcal{D}^0 \subset \mathcal{D}^1 \subset\mathcal{D}^2 \subset \cdots$. Contrary to the case of  {static} feedback \cyan{linearizability}, they need not be involutive but the first noninvolutive one {has} to contain 
 an involutive subdistribution of corank one. \red{The main idea of the proposed algorithm is  to replace, \cyan{at each step,} the first noninvolutive distribution  by its  involutive subdistribution of corank one},
 thus for the {prolonged} system we gain at least one new involutive distribution. Our algorithm is constructive, gives sufficient conditions for flatness, and can be seen as the dual of the dynamic feedback linearization algorithm of Battilotti and Califano~{\cite{battilotti2003further,califano2005dynamic}}. 
\end{abstract}



%
%
%



\section{Introduction}
{I}{n} this paper, we study flatness of nonlinear control systems of the form
$$
\Xi\,:~\ \dot x=~F(x,u),
$$
where $x$ is the state defined on a open subset $X$ of $\mathbb{R}^n$  and~$u$ is the control taking values in an open subset $U$ of $\mathbb{R}^m$ (more generally, an $n$-dimensional manifold $ X$ and an~$m$-dimensional manifold $U,$ resp.). The dynamics $F$ \red{are such that $\mrk \pfrac{F}{u} = m$, and are supposed smooth} (the word smooth will always mean $\mathcal{C}^\infty$-smooth).
The notion of flatness was introduced in control theory in the 1990's, by Fliess, L\'evine, Martin and Rouchon \cite{fliess61vine}, see also
{\cite{isidori1986sufficient,jakubczyk1993invariants,levine2009analysis, martin1992phd, pomet1995differential, schlacher2007construction}, 
and the references therein,}  
and has attracted a lot of attention because of its multiple applications in the problem of constructive controllability and motion planning {(for some recent applications, see \cite{ brahmi2021flatness, kolar2017time, li2021quadruped, louembet2009path, nguyen2020flat, wipke1999advisor}).} 
Flat systems form a class of control systems whose set of trajectories can be parametrized by~$m$ functions and their time-derivatives,~$m$ being the number of controls. More precisely, the system  $\Xi\,: \dot x=~F(x,u)$
is \textit{flat} if we can find~$m$ functions, $\varphi_i(x,u,\dots,u^{( l)})$ such that
\begin{equation} \label{descript}
x=\gamma(\varphi,\dots,\varphi^{(s-1)}) \mbox{ and }
u=\delta(\varphi,\dots,\varphi^{(s)}),
\end{equation}
for certain {integers $l$ and $s$}, where $\varphi=(\varphi_1,\dots,\varphi_m)$ is called \textit{a flat output}. Therefore the time-evolution of all state and control variables can be determined from that of flat outputs without integration and all trajectories of the system can be completely parameterized.
Flat systems can be seen as a generalization of {static feedback linearizable systems}. Namely, they are  linearizable via dynamic, invertible and endogenous feedback, see {\cite{fliess1992systemes, fliess61vine, martin1992phd,pomet1995differential}}.
The simplest flat systems that are not static feedback linearizable are control-affine systems that become static feedback linearizable after an invertible one-fold prolongation of a suitably chosen control  (which is the simplest dynamic feedback). A complete geometric characterization of that class of systems has been proposed by the authors in \cite{nicolau2016two} for the two-input case, and  \cite{nicolau2017flatness} for the multi-input case {(see also \cite{gstottner2020linearization, gstottner2021necessary, nicolau2016flatness} for related results on linearization via a two-fold prolongation).}
Inspired by those {results}, the goal of this paper is to propose a constructive algorithm  to dynamically linearize two-input control systems via successive one-fold prolongations of a control that has to be suitably chosen at each step of the algorithm. {Preliminary results, on which this paper is based, were presented in~\cite{nicolau2021dynamic}.}

Linearization via successive one-fold prolongations  requires special properties of the linearizability distributions $\mathcal{D}^0 \subset \mathcal{D}^1 \subset\mathcal{D}^2 \subset \cdots$. Contrary to the case of  {static} feedback linearization, they need not be involutive but the first noninvolutive one {has} to contain a sufficiently large  {(actually of corank one)} involutive subdistribution.
 The main idea of the proposed algorithm is, at each new step,  to prolong the system in such a way that the first noninvolutive distribution~{$\mathcal{D}^k$} is replaced by its largest involutive subdistribution~{$\mathcal{H}^k$}, thus for the prolonged system we gain at least one new involutive distribution.
 {Working with the first noninvolutive distribution guarantees \red{that}
all objects that we use are  invariant with respect to feedback transformations.}
 {If the algorithm
 terminates with a system that is static feedback linearizable, then the starting system of the algorithm will be called ``linearizable  via successive one-fold prolongations with involutivity gain'' (shortly, LSOPI)}.
In the present paper, {using  necessary and sufficient conditions for~$\mathcal{H}^k$ to exist, we provide conditions for that step-by-step prolongation to dynamically linearize the system}.  The conditions are easy to check, their verification involves differentiation and algebraic operations only, without solving PDE's or bringing the system into a normal form.
If such an involutive subdistribution exists, we uniquely identify it in {almost} all cases and provide an explicit construction {(except for one particular subcase, where {the algorithm may not be conclusive})}. Such an involutive subdistribution has to be calculated at each step of the algorithm and plays a crucial role for linearization via successive one-fold prolongations with involutivity gain. Indeed, it will take the place of the first noninvolutive distribution and, moreover, it enables us to define \red a to-be-prolonged control.

Our {solution is related to important results of  Battilotti and Califano~ \cite{battilotti2003further, califano2005dynamic} (see also \cite{battilotti2004constructive, battilotti2008geometric} for the multi-input case),
where a constructive condition for dynamic feedback linearization of two-input nonlinear systems, also based on prolongations, is presented. The algorithm that we propose can be seen as the dual of the results of~\cite{battilotti2003further,califano2005dynamic} {whose authors} identify the inputs that {\it should not be extended
through prolongations} with the help of the {\it last}  noninvolutive linearizability distribution while here, we work with the {\it first} noninvolutive distribution  and we identify the inputs that {\it have to be prolonged}\red.
%
\red At each step of our algorithm, general static feedback transformations are allowed (involving the original states and the controls prolonged at the previous steps). \red In
 \cite{battilotti2003further},
 the class of allowed dynamic feedback laws is restricted to prolongations  of the original inputs plus a final static feedback on
the extended system, \red{but} {in~\cite{califano2005dynamic}, \red{the framework is extended to
general feedback transformations and}
two static state feedback actions (called  the
direction feedback and the reduction feedback) may also be applied at each step of the algorithm.}
In this paper, we work under constant rank assumptions, \red{which simplifies certain aspects of the analysis. In \cite{battilotti2004constructive}, the authors also investigate the  more general situation
where the linearizability distributions may lose constant rank.}

{Other  results  related to dynamic feedback linearizability via prolongations can be found, for instance, in \cite{fossas2000linearization, franch2005linearization, levine2025differential} (see also the references therein);} {in the aforementioned papers only original controls are prolonged, so the problem reduces to the classical static feedback linearization, namely, to checking involutivity of some distributions.}

\red{Some systems, although dynamically linearizable via successive one-fold prolongations (referred to as LSOP-systems), may  however, not satisfy the LSOPI conditions.
In other words, 
we may gain no involutive distribution or the involutivity of some distributions may even be destroyed by
 the prolongation.
 This distinction was not made in our previous paper~\cite{nicolau2021dynamic}. 
 In Section~\ref{sec:interplay} we examine the interplay between LSOPI, LSOP and flatness. Clearly, LSOPI-systems form a subclass of LSOP-systems
 and we show that this inclusion is strict.}
 Specifically, we demonstrate that certain systems can be dynamically linearized via successive one-fold prolongations
{but do not admit any
sequence of prolongations where each new prolongation would yield} at least one new involutive distribution (\red{hence} these systems are LSOP without being {LSOPI}).
  \red{Furthermore, for two-input control systems, we establish that flatness is equivalent to LSOP-linearizability and explore its connections with linearization via an $\ell$-fold prolongation of a suitably chosen control. In general, LSOP and linearization via  $\ell$-fold prolongation require knowledge of the system’s flat output to determine which control should be prolonged. When no flat output is known, this method necessitates either transforming the system into a specific normal form or solving certain PDEs, making it generally non-constructive. This limitation precisely motivates the introduction of the LSOPI notion in this paper, which proves to be constructive in almost all cases. As discussed above, the construction relies, at each step of the LSOPI-algorithm, on an involutive subdistribution~$\mathcal{H}^k$ of corank one in~$\mathcal{D}^k$, the first non-involutive distribution of the system at the given step of the algorithm.
%
$\mathcal{H}^k$ 
} 
plays a crucial role {in} the proposed algorithm. Its construction depends on how the involutivity of~$\mathcal{D}^k$ is lost and we
distinguish three mutually exclusive cases.
While our previous work~\cite{nicolau2021dynamic} considered only the first two cases, this paper {also studies} the third {more involved} case and provides  an explicit construction of~$\mathcal{H}^k$ in nearly all subcases (with the exception of one particular subcase, where {it may be challenging to check whether the algorithm 
{is conclusive)}.} 
Unlike~\cite{nicolau2021dynamic}, the current paper presents {detailed algorithms and} complete proofs of all results.} {Additionally, it introduces an equivalent geometric characterization of~$\mathcal{H}^k$ in the {second} case avoiding concepts like the Engel rank  used in~\cite{nicolau2021dynamic}. 

The paper is organized as follows. In Section~\ref{sec_flatness}, we recall the {notion} of flatness.
\red{In Section~\ref{sec_main results}, we formalize and present our main result:\vspace{-0.2cm}
\begin{list}{-}{}
 \item first, we define the notion of LSOPI in Section~\ref{ssec: lsopi def};\vspace{-0.2cm}
 \item second,  we propose an algorithm, based on Theorem~\ref{thm: algo}, for dynamic linearization of LSOPI-systems in Section~\ref{ssec: algo}; we also discuss its limitations and  explain
 \cyan{what is needed to render it fully constructive;}\vspace{-0.2cm}
\item third, we
give verifiable conditions to check the existence of the subdistribution~$\mathcal{H}^k$, on which the LSOPI-algorithm is based, and if~$\mathcal{H}^k$ exists, an explicit way to construct it in Section~\ref{ssec: construction H};\vspace{-0.2cm}
\item finally combining all above results, we formulate Theorem~\ref{thm: algo constructive} and constructive Algorithm~\ref{algo detailed constructive} in Section~\ref{ssec: algo constructive}.\vspace{-0.2cm}
\end{list}
%
%
In Section~\ref{sec:interplay} we discuss relations between the notions of LSOPI, LSOP, linearization via an $\ell$-fold prolongation and flatness.
In addition to the examples provided in Section~\ref{sec:interplay}, we illustrate the LSOPI-algorithm via the chained form in Section~\ref{sec_examples}. Finally,  we} provide \cyan{all proofs (except for those of Theorems~\ref{thm: algo} and~\ref{thm: equiv flat LSOP}  proved directly after their statements)} in Section~\ref{sec_proofs}. {Appendices A and~B
 further elaborate} on the properties and construction of~$\mathcal{H}^k$.}

\section{Flatness} \label{sec_flatness}

Fix an integer $l \geq -1$ and denote $U^l = U \times \mathbb{R}^{ml}$ and $\bar u^l = (u,\dot  u, \dots, u^{(l)})$. For $l=-1$,  the set $ U^{-1} $ is empty and $\bar u^{-1}$ in an empty sequence.

\begin{definition} \label{Def platitude}The system $ \Xi: \dot x= F(x,u)$ is \textit{flat} at $(x_0, \bar u_0^l) \in X \times U^l$, for  $l\geq -1$, if there exists 
\cyan{an open set $\mathcal O^{l}$ containing} $(x_0, \bar u_0^l)$
and~$m$ smooth functions $\varphi_i=\varphi_i(x, u,\dot  u, \dots, u^{(\cyan q)})$, \cyan{with $q\leq l$}, $1 \leq i\leq m$, defined in $\mathcal O^l$, having the following property: there exist an integer $s$ and smooth functions~$ \gamma_i$, $1 \leq i \leq n$, and~$\delta_j$, $ 1 \leq j \leq m$, such that
$$
x_i = \gamma_i (\varphi, \dot \varphi, \dots, \varphi^{(s-1)}) \mbox{ and } u_j=  \delta_j (\varphi, \dot \varphi, \dots, \varphi^{(s)})
$$
for any $C^{l+s}$-control $u(t)$ and corresponding trajectory $x(t)$
that satisfy $(x(t),$ $ u(t), \dots ,$ $u^{(l)}(t)) \in~\mathcal O^l$, where $\varphi=(\varphi_1, \dots, \varphi_m)$ and is  called a \textit{flat output}.
\end{definition}

\cyan{It can be proven using the extension algorithm, see, e.g., \cite{di1989rank,respondek1990right} that the differentials of the components of
 flat outputs and their time-derivatives are independent, implying that they  do not satisfy any differential equation of the form $\Phi(\vp, \dot \vp, \ldots, \vp^{(r)}) = 0$.}

\red{It is well known
that systems linearizable via invertible static feedback are flat.}
To see that, consider the control-affine system
\begin{equation} \label{Sigma}
\Sigma\ :\ \dot x = f(x)+\sum_{i=1}^mu_ig_i(x),
\end{equation}
where $f$ and $g_1, \ldots, g_m$ are smooth vector fields on~$X$. The system~$\Sigma$
is {said} linearizable by static {invertible} feedback if it is equivalent, via a diffeomorphism $z= \phi(x)$ and an invertible {static} feedback transformation $u=\alpha(x)+\beta(x)v$, to a linear controllable system of the form $\dot z=Az+Bv$.
The problem of static feedback linearization  was solved by 
Brockett~\cite{brockett1978feedback} (for a smaller class of transformations) 
and then by 
Jakubczyk and Respondek \cite{jakubczyk1980on} and, independently, by Hunt and Su~\cite{hunt1981linear}, who gave geometric necessary and sufficient conditions recalled in Theorem~\ref{thm F-lin systems} below. 
%
For the control-affine system~$\S$, we define inductively the sequence of distributions $\mathcal{D}^{j+1} = \mathcal{D}^j + [f,\mathcal{D}^j]$, where~$\mathcal{D}^0$ is given by  $\mathcal{D}^0= \mspan \{g_1,\ldots, g_m\}$ and {we} denote $[f,\mathcal{D}^j] =\{[f,\xi]: \red{\text{ for all }} \xi \in \mathcal{D}^j\}$.
In other words \red{(see \cite{jakubczyk1980on} where expression~\eqref{eq: sequence D j} is deduced from the above definitions)}
\begin{equation}\label{eq: sequence D j}
\mathcal{D}^j  = \mspan\{ad_f^q g_i, \dleq 1 i m, \dleq 0 q j\}, \mbox{ for } j\geq 0.
\end{equation}
The distributions $\mathcal{D}^j$ {are} called linearizability distributions {and the integer $j$ is called the order of the distribution $\mathcal{D}^j$.}

\begin{theorem}[{See} \cite{hunt1981linear,jakubczyk1980on}]
\label{thm F-lin systems}
 {\it $\Sigma$ is locally static feedback linearizable, around $x_0 \in~X$, if and only  if
 \begin{enumerate}[\normalfont ({FL}1)]
  \item For {all} $j \geq 0$, the distributions $\mathcal{D}^j$, around $x_0 \in~X$, are involutive and of constant rank.
  \item $\mathcal{D}^{n-1} = TX$.
 \end{enumerate}}
\end{theorem}
The geometry of static feedback linearizable systems is {thus} 
given by the following sequence of nested involutive distributions:%
$$
\mathcal{D}^0 \subset \mathcal{D}^1\subset \cdots \subset\mathcal{D}^{n-1} = TX. %
$$

In this paper, 
we propose a constructive algorithm  to dynamically linearize two-input control systems via successive one-fold prolongations.
We make the following assumption:
\begin{assumption}\label{assum: constant rk and acces}
  {\it 
  Unless stated otherwise, we assume that \emph{all ranks are constant} in a neighborhood of the nominal point around which we work. {All presented results  are {thus} valid on an open and dense subset of $X\times U \times \mathbb{R}^{m\ell}$} (with the integer~$\ell$   being related to 
   \red{the highest order of the time-derivative $u^{(\ell)}$ involved in}
  the  proposed algorithm), where all ranks are constant,
and hold locally, around {any} given point of that subset.
}
 \end{assumption}

Before giving our main results, let us introduce \cyan{some notations as well as} the notion of corank that will be frequently used in the rest of the paper.

\begin{notation} \label{def corank}
 Let $\mathcal{D}$ and $\mathcal{E}$ be two distributions of constant rank \red{on a smooth manifold $X$,   and $f$ a vector field on $X$.} Denote $[\mathcal{D},\mathcal{E}] = \{[\xi, \zeta]:  \red{\text{ for all }} \xi \in \mathcal{D}$ $\red{\text{and }} \zeta \in \mathcal{E}\}$ and  $[f,\mathcal{E}] = \{[f,\zeta]:  \red{\text{ for all }}\zeta \in \mathcal{E}\}$.  If $\mathcal{D}\subset \mathcal{E}$, \cyan{then $\mcork(\mathcal{D}\subset\mathcal{E}) = e- d$ is
 the corank of the inclusion $\mathcal{D}\subset \mathcal{E}$,
 where $\mrk \mathcal{D} = d$ and $\mrk \mathcal{E} = e$.}
\end{notation}

\begin{notation}
\red{For} unit vector fields in $x$-coordinates, we  write $e_i = (0, \ldots, 0, 1,$ $ 0, \ldots, 0)^\top =~\partial_{x_i}$, and thus for a vector field $f$ on $X$, we have $ f = \sum_{i=1}^n f_i(x)\partial_{x_i}$.
\end{notation}

Since
we deal with {two-input systems of the form~$\S$, i.e., with $m=2$,} that are not static feedback linearizable,
and for which, according to  Assumption~\ref{assum: constant rk and acces} above,   all ranks are supposed constant,
then either 
{all distributions~$\mathcal{D}^j$ are involutive but $\mrk \mathcal{D}^{n-1} <n$, that is, (FL2) is not satisfied, or one of the distributions~$\mathcal{D}^j$ is not involutive  (condition (FL2) may or may not be verified) implying that  (FL1) fails. In the former case, the system is not flat (because it is not controllable). Therefore, from now on, we assume that}
there exists an \red{order}~$k$ such that~$\mathcal{D}^k$ is not involutive, 
and we
suppose that   \red{\textit{$k$ is the smallest order for which~$\mathcal{D}^k$ is not involutive}. We call $k$ the {\it non-involutivity index of $\S$}.}
%
\red{We say that the  non-involutivity index of $\S$
does not exist
if all distributions  $\mathcal{D}^j$, $j\geq 0$, are involutive (i.e., $\S$ is static feedback linearizable).}

\section{Main results} 
\label{sec_main results}

%

Throughout, we will consider two-input systems
and, to simplify the exposition of the paper, we will treat only the control-affine case. That is, we work with
systems 
\begin{equation} \label{eq: Sigma}
 \Sigma: \dot x = f(x)+u_1g_1(x)+u_2g_2(x),
\end{equation} 
where~$x\in X$, $u=(u_1,u_2)^\top\in \R^2$,  and $f$,~$g_1$,~$g_2$ are $\mathcal{C}^{\infty}$-smooth vector fields. 

\cyan{Below, when we say that    {an object} is feedback invariant, we mean that it is invariant under invertible static feedback transformations of the form $f \mapsto f + g\a$ and $g \mapsto g \b$, \red{where  $g= (g_1, g_2)$.} 
The conditions of all results presented in the paper do not depend on coordinates and are feedback invariant.
The following 
proposition, whose proof is given in Section~\ref{sec_proofs}, gives a necessary condition for flatness of two-input control systems that are not static feedback linearizable.
\begin{proposition}\label{prop: nec cond flatness}
 Consider 
 $\S : \dot x = f(x)+u_1g_1(x)+u_2g_2(x)$ and suppose that its first noninvolutive distribution is~$\mathcal{D}^k$. Then
 \red{the following conditions hold:}
 \begin{enumerate}[\normalfont (1)]
  \item  The distributions $\mathcal{D}^j$, for $\dleq 0 j k$, are feedback invariant.
  \item
\red{If} $\S$ is locally flat at $(x_0, \bar u_0^l) \in X \times U^l$, for  a certain $l\geq -1$, then $\mcork (\mathcal{D}^{k-1}\subset \mathcal{D}^{k}) =2$ around~$x_0$.
 \end{enumerate}
%
\end{proposition}
}

\subsection{\red{LSOPI-definition}}\label{ssec: lsopi def}

\begin{definition}[Prolongation]\label{def: prolongation}
      Given a control system~$\S$ 
      of {the} {form~\eqref{eq: Sigma}}, 
      we construct its prolongation ${\widetilde \S^{(1,0)}}$ by, first applying an invertible feedback of the form $u= {\a(x) + } \beta(x) \tilde u$, where \red{$\t u = (\t u_1^p, \t u_2)$ and} $\b(x)$ is an invertible $2\times2$-matrix,  to get
       $$
      {\widetilde \S} : \dot x = {\t f(x)}+\red{\t u_1^p}\t  g_1(x)+\t u_2\red{\t g_2^p}(x),
      $$
      with {$\tilde f = f +g  \alpha $ and  $\tilde g = g\beta$, where $g = (g_1, g_2)$ and $\tilde g = (\tilde g_1,\red{\tilde g_2^p})$}, and then prolonging the first (feedback modified) control as $\red{\dot {\t u}_1^p} = v_1$, which justifies the notation ${\widetilde \S^{(1,0)}}$,  and keeping the second control $\t u_2 = v_2$ to obtain
         \begin{small}  $$
      {\widetilde \S^{(1,0)}} : \dot{x}_p =
\begin{pmatrix}
    \dot{x} \\
      \dot {\t u}_1^p
\end{pmatrix}
=
\begin{pmatrix}
    \tilde{f}(x) + \tilde{u}_1^p \tilde{g}_1(x) \\
       0
\end{pmatrix}
+ v_1
\begin{pmatrix}
    0 \\
       1
\end{pmatrix}
+ v_2
\begin{pmatrix}
    \tilde{g}_2^p(x) \\
    0
\end{pmatrix}
      $$
      \end{small}
     %
      where ${x_p = (x^\top, \red{\t u_1^p})^\top}$ and $(v_1, v_2)^{\top}\in \mathbb{R}^2$.
      \red{The control $\t u_1^p$ is called  a \emph{to-be-prolonged control}.}
      Denote by $\mathcal{D}^j_p$ 
      the linearizability distributions of ${\widetilde \S^{(1,0)}}$.
\end{definition}

In this paper, we propose an algorithm that constructs a sequence of systems~$\S_i$, for $\dleq 0 i \ell$, \red{with the following properties:}
%

\begin{definition}[{LSOPI}] \label{def: LSOPI} 
The control system $ \S: \dot x = f(x)+u_1g_1(x)+u_2g_2(x)$ is said to be  {\it dynamically linearizable  {via} successive one-fold prolongations with involutivity gain} (shortly, called {LSOPI})
 if there exist an integer $\ell$ and a sequence of systems ${\S= \S_0 }, \S_1, \ldots, \S_\ell$, where $\S_{i+1}$ is the prolongation   $\S_{i+1} =  {\widetilde \S_i^{(1,0)}} $, for $\dleq 0 i \ell-1$,  such that
  {$k_0 < k_1<\cdots < k_{\ell -1}$}, 
  with $k_i$
  \red{the non-involutivity index of~$\S_i$,}
  and $\S_\ell$ is locally static feedback linearizable  {(i.e.,  $k_{\ell}$ does not exist).}
\end{definition}

\begin{notation} \label{notation Sigma i}
{The  {subindex} 0 of  $\Sigma_0$
indicates the fact that $\Sigma_0$ is the original control system~{$\S$} that we would like, {if possible},  to dynamically linearize via successive one-fold  prolongations. 
For~$\S_i$, $\dleq 0 i \ell$,  we will write
\begin{equation}\label{eq: Sigma i}
 \dot x^i = f^i (x^i)+u_1^i g_1^i(x^i)+u_2^i g_2^i (x^i), \quad x^i\in X^i,
\end{equation}
with the upper-index~$i$ corresponding to the state, state-space, controls, and vector fields of~$\S_i$. 
To each system~$\S_i$ 
we will associate the  distributions $\mathcal{D}_i^j$, for $j\geq 0$, and {denote} by $\mathcal{D}_i^{k_i}$ the first noninvolutive distribution {of~$\S_i$.}
{We say that $k_i$
does not exist}
if all distributions  $\mathcal{D}_i^j$ are involutive, {that is, for   the sequence of Definition~\ref{def: LSOPI}, we have $i = \ell$}.}
 \end{notation}

\begin{remark}[LSOPI versus LSOP]\label{rk:LSOP}
 \red{Systems $ \S: \dot x = f(x)+u_1g_1(x)+u_2g_2(x)$  verifying Definition~\ref{def: LSOPI} \cyan{(that is,  $\S_{i+1} =  {\widetilde \S_i^{(1,0)}} $ and $\S_\ell$ is static feedback linearizable)} but for which we drop the requirement $k_0 < k_1<\cdots < k_{\ell -1}$, are called {\it dynamically linearizable  {via} successive one-fold prolongations} (shortly, called {LSOP}).}
 In other words, \red{for LSOP-systems,}
 it is possible to have  $k_i \geq k_{i+1}$ for some  integers~$i$ (i.e., for the system $\S_{i+1}$ constructed at step  $i$,  we gain no involutive distribution or the involutivity of some distributions is destroyed by
 the prolongation), see Examples~\ref{ex: not necessary conditions} \red{and~\ref{ex:involutivity destroyed} in Section~\ref{sec:interplay}.
  This also explains why the linearization property introduced in Definition~\ref{def: LSOPI} is referred to as ``with involutivity gain'': indeed, in the construction of the sequence ${\S = \S_0}, \S_1, \ldots, \S_\ell$, the (new) system $\S_{i+1}$ constructed at step  $i$
     admits at least one additional involutive distribution compared to its predecessor $\S_i$, as illustrated by  the following diagram: 
%
 \begin{center}
    \begin{footnotesize}
    $$
\hspace{-0.2cm}\begin{array}{r @{\,} c @{\,} l}
\S_i &
\xrightarrow{\tiny \begin{array}{l}
                    (\a^i,\beta^i)\\
                    u^{i+1}_1 = \dot {\t u}_1^{i,p}, u^{i+1}_2 = {\t u}^i_2
                    \end{array}
}
&
\S_{i+1} =  \widetilde \S_i^{(1,0)}\\[1ex]
\begin{array}{r}
\mathcal{D}_i^0 \subset  \cdots \subset \mathcal{D}_i^{k_i-1}\subset\red{\mathcal{D}_i^{k_i}}\\
\red{\tiny \text{1st non-inv}}
\end{array}
&
&
\begin{array}{l@{}c@{}r}
\mathcal{D}_{i+1}^0 \subset  \cdots \subset& {\mathcal{D}_{i+1}^{{k_i}}}&\subset
\cdots \subset
\mathcal{D}_{i+1}^{k_{i+1}-1}\subset\red{\mathcal{D}_{i+1}^{{k_{i+1}}}}\\
& {\tiny \text{inv.}}& \red{\tiny \text{1st non-inv}}
\end{array}
\end{array}
$$
\end{footnotesize}
\end{center}
 The classes of systems  possessing the above-defined linearization properties are related as follows (where we use calligraphic font to denote classes of systems):
  $$\mathcal{LSOPI} \subset \mathcal{LSOP}$$ 
 with the inclusion being strict, see Theorem~\ref{thm: equiv flat LSOP} in Section~\ref{sec:interplay} where we also discuss these  properties with respect to flatness.
}
\end{remark}

 The importance of the proposed notion of  linearizability via successive one-fold prolongations with involutivity gain (LSOPI)
 %
 %
 is two-fold. First, the assumption $k_i<k_{i+1}$ \red{(that is, the ``involutivity gain'')} implies that passing from~$\S_i$ to $\S_{i+1}$ decreases the dimension of the nonlinear subsystem.
 Namely for~$\S_i$ evolving on $X^i$, with $\dim X^i = n_i$, that dimension is $n_i -2k_i$, while for $\S_{i+1}$ evolving on $X^{i+1}$, with $\dim X^{i+1} = n_{i+1}= n_i +1$, that dimension is $n_{i+1} -2k_{i+1}$, which is smaller if $k_i< k_{i+1}$, compare with the proof \cyan{of} Theorem~\ref{thm: algo}.
 Second, the construction of a prolongation $\S_{i+1} =  {\widetilde \S_i^{(1,0)}} $, for which $k_i< k_{i+1}$, as required by Definition~\ref{def: LSOPI}, is closely
related (see Proposition~\ref{prop: equiv H feedback} below) to the existence of an involutive subdistribution of corank one in the first noninvolutive distribution ${\mathcal{D}_i^{k_i}}$ associated to $\S_{i}$. This last property can be checked for a large class of systems with two inputs via differentiation and algebraic operations only, which we present in detail in Section~\ref{ssec: construction H}.

To summarize, for 
LSOPI-systems we propose a constructive algorithm (due to the second of the above properties {and to an explicit construction of the prolongation~$\S_{i+1}$}) that dynamically linearizes the system in a finite number of steps (due to the first of the above properties) and we will provide such \red{a constructive algorithm} in Section~\ref{ssec: algo constructive}\red.

Finally, the {extension} of our results to the case of general nonlinear control systems 
$ \Xi : \dot x= F(x,u)$,~$x\in X$, $u\in \R^2$, is straightforward. Indeed, by prolonging both its inputs,  $ \Xi$ {is} (dynamically) transformed into the  control-affine system
$$
 \Xi^{(1,1)}:\left\{
\begin{array}{lcl}
 \dot x &= & F(x,u_1, u_2)\\
 \dot u_1  &= & v_1\\
  \dot u_2 &= & v_2,
\end{array}
\right.
$$
and it is the control-affine prolongation $ \Xi^{(1,1)}$ {to which} we apply the  constructive  algorithm proposed in Section~\ref{ssec: algo constructive}  {by taking} $ \S_0 = \Xi^{(1,1)}$.

\subsection{Algorithm for dynamic linearization {of LSOPI-systems}}\label{ssec: algo}

{From now on, we use the notation~$\S$, without a subindex, to refer to any control-affine system, not necessarily the starting system of the proposed linearization algorithm. For instance, Proposition~\ref{prop: equiv H feedback} below applies to any control-affine system~$\S$ and will be applied to  $\S = \S_i$ at each step of the linearization algorithm. Recall that we denote the {linearizability} distributions associated to~$\S$ by $\mathcal{D}^j$ (without a \red{subscript}), and those of its prolongation  $\widetilde{\S}^{(1, 0)}$ by  $\mathcal{D}_p^j$ (where the subscript $p$ refers to the prolongation).}

\begin{proposition}\label{prop: equiv H feedback}
Consider the system $\S : \dot x = f(x)+u_1g_1(x)+u_2g_2(x)$ and \red{let $k$ be its non-involutivity index, i.e.,}
its first noninvolutive distribution is~$\mathcal{D}^k$.
For~$\S$, the following \red{conditions} are equivalent:
\begin{enumerate}[\normalfont (i)]
 \item There exists an involutive distribution~$\mathcal{H}^k$ verifying   $\mathcal{D}^{k-1} \subset \mathcal{H}^k \subset \mathcal{D}^k$, with both inclusions of corank one 
 {(and where $\mathcal{D}^{k-1} = 0$ if $k=0$)};
 \item There exists an invertible static feedback transformation $u = {\a(x)} + \beta(x)\tilde u$ such that {the linearizability distributions  $\mathcal{D}^0_p\subset \cdots\subset \mathcal{D}_p^k$ of} the prolongation
$$
\widetilde{\S}^{(1, 0)} :
\left \{
\begin{array}{lcl}
\dot x &= &\t f(x)+ \red{\t u_1^p} \tilde g_1(x)+v_2 \red{\tilde  g_{2}^p}(x) \\
\red{\dot {\t u}_1^p}&= &v_1
\end{array}
\right.
$$
{are involutive and satisfy} $\mrk  \mathcal{D}_p^k = 2k+2$.
\end{enumerate}
Moreover, if  condition \nf{(i)}, equivalently condition \nf{(ii)}, is satisfied,
then there exists $\red{\widetilde {\mathcal{G}}_2^p} = \mspan\{\red{\tilde  g_{2}^p}\}$ of {corank} one in~$\mathcal{D}^0$, 
such that $ \mathcal{H}^k  =  \mathcal{D}^{k-1} + \mspan\{ad_f^k \red{\tilde  g_{2}^p}\}$
and, moreover, {if $k \geq 1$, then} all distributions $\mathcal{H}^j = \mathcal{D}^{j-1}+\mspan \{ad_f^j \red{\tilde  g_{2}^p}\}$, for $\red 0\leq j \leq k-1$, are also involutive.
\end{proposition}

\begin{remark}
 {We use the notation  $\mathcal{H}^j$, for all $\dleq 0 j k$, because all distributions satisfy $\mathcal{D}^{j-1} \subset \mathcal{H}^j \subset \mathcal{D}^j$ with both inclusions of corank one \red{(recall that we put $\mathcal{D}^{-1}= 0$)}. The difference is that $\mathcal{D}^j$, for $\dleq 0 j k-1$, are involutive but~$\mathcal{D}^k$ is not.}
\end{remark}

{Recall} that the distributions $\mathcal{D}^0, \ldots, \mathcal{D}^k$ are feedback invariant and 
{if~$\mathcal{H}^k$, satisfying (i), is unique, then~$\mathcal{H}^k$ depends neither on~$\a$  nor on the choice of~$\red{\t g_2^p}$ defining $ \red{\widetilde {\mathcal{G}}_2^p} $ (i.e., neither on transformations of the form $f \mapsto f + g\a$ nor on   $(g_1, \red{\t g_2^p})  \mapsto(g_1, \red{\t g_2^p}) \beta$, with $\beta = (\beta_i^j)_{\dleq 1 {i,j} 2}$ invertible and such that $\beta_1^2 = 0$).}
%
%
The involutive subdistribution~$\mathcal{H}^k$ of corank one in~$\mathcal{D}^k$ plays a crucial role for
{LSOPI}.  Indeed, as we will see {at the end of this subsection}, with the help of~$\mathcal H^k$, we can define {an}  invertible static feedback transformation $u = {\a(x) }+ \beta(x)\tilde u$ leading to a prolongation $\widetilde{\S}^{(1, 0)}$ whose {distributions $\mathcal{D}_p^{j} = \mspan\{\partial_{\red{\t u_1^p}}\} + \mathcal{H}^j$, for $\dleq 0 j k$, are involutive}
%
{(see the proof of Proposition~\ref{prop: equiv H feedback}, as well as Algorithm~\ref{algo detailed} below)}, so the involutive {distribution $\mathcal{D}_p^{k}$ (whose construction is based on the involutive subdistribution~$\mathcal{H}^k$)} takes the {role} of the noninvolutive one~$\mathcal{D}^k$ of~$\S$.
It is namely {the above} fact that the  algorithm {for {LSOPI}} is based {on}. {Given the importance of the involutive subdistribution~$\mathcal{H}^k$, we will call it a LSOPI-distribution, {which we formalize as follows.}

\begin{definition}[LSOPI-distribution]\label{def: LSOPI-distribution}
 Consider the control system $\S : \dot x = f(x)+u_1g_1(x)+u_2g_2(x)$ and
 \red{let $k$ be its non-involutivity index}.
 \cyan{A} distribution $\mathcal{H}^k \subset \mathcal{D}^k$
 is called a LSOPI-distribution if it is
 involutive and satisfies $\mathcal{D}^{k-1} \subset \mathcal{H}^k \subset \mathcal{D}^k$, with both inclusions of corank one.
\end{definition}
}

For~$\S$, a LSOPI-distribution either exists and is unique, or there exist many, or there does not exist any. We formalize it in the following \red{remark.} 

\begin{remark}\label{rk: cases LSOPI-distribution}
 Consider a 
  control system~$\S$ of the form~\eqref{eq: Sigma}. \red{It is immediate that} one of the following mutually exclusive conditions holds:
  \begin{enumerate}[\normalfont (i)]
  \item $k$ does not exist,
  that is,  all 
  $\mathcal{D}^j$  are involutive;
  \item  $k$   exists and there exists a unique LSOPI-distribution $\mathcal{H}^{k}$;
  \item  $k$   exists  but a LSOPI-distribution $\mathcal{H}^{k}$ does not exist;
    \item  $k$   exists and there \cyan{exist many} LSOPI-distribution~$\mathcal{H}^{k}$.
 \end{enumerate}
\end{remark}

Now we propose an algorithm identifying all LSOPI-systems provided {that} at each step either (i) or~(ii) or (iii) of \red{Remark~\ref{rk: cases LSOPI-distribution}} holds.

 \color{black}

\begin{theorem}\label{thm: algo}
 Consider a two-input control system~$\S$ of the form~\eqref{eq: Sigma},
 and assume that a sequence of systems $\S = \S_0, \S_1, \ldots,\S_{i-1},  \S_i$ has been constructed such that for each $\dleq 0 s i-1$, there exist $k_s$ and a unique LSOPI-distribution $\mathcal{H}_s^{k_s}$ {of $\S_s$} and that
 $\S_{s+1} = \widetilde \S_s^{(1,0)}$. 
 Assume that one of the three following cases holds for~$\S_i$:
 \begin{enumerate}[\normalfont (i)]
  \item either {$k_i$ does not exist},
  that is,  all linearizability distributions $\mathcal{D}_i^j$  are involutive;
  \item or $k_i$   exists and there exists a unique LSOPI-distribution~$\mathcal{H}_i^{k_i}$; in this case, set $\S_{i + 1} = \widetilde \S_i^{(1,0)}$;
  \item or  $k_i$   exists  but a LSOPI-distribution~$\mathcal{H}_i^{k_i}$ does not exist;
 \end{enumerate}

  Then the sequence $\S_0, \S_1, \ldots,  \S_i$ always terminates at, say,~$\S_{i^*}$ which satisfies either {\normalfont (i)} or {\normalfont (iii)}; if~$\S_{i^*}$ satisfies \nf{(i)} {and there exists $\rho$ such that $\mathcal{D}_{i^*}^{\rho} = TX^{i^*}$, where $X^{i^*}$ is the state-space of~$\S_{i^*}$}, then~$\S_{i^*}$ is static feedback linearizable, implying that~$\S$ is LSOPI, \cyan{and in that case  $i^* = \ell$};
  {if~$\S_{i^*}$ satisfies~\nf{(i)} and  $\rho$ does not exist, then~$\S$ is not flat (in particular, not  LSOPI);}
  if~$\S_{i^*}$ satisfies  {\normalfont (iii)}, then~$\S$ is not LSOPI.

\end{theorem}

{The two basic issues of, first, how to check {(i), (ii), and  (iii),} as well as, second, how to treat the case of a non-unique $\mathcal{H}_i^{k_i}$, are discussed in Section~\ref{ssec: construction H} and allow to strengthen Theorem~\ref{thm: algo} to Theorem~\ref{thm: algo constructive}.}
We present the proof of the above theorem here, as it serves to justify the LSOPI {Algorithm~\ref{alg iteration i}}.

\begin{proof}
Assume that~$\S_i$ satisfies (i). Then we set $i^* = i$ and it is immediate to conclude that~$\S_{i^*}$ is static feedback linearizable (implying that $\S = \S_0$ is LSOPI) if $\rho$ exists, \cyan{and in that case $i^* = \ell$,} and that is not flat (in particular, not  LSOPI) if $\rho$ does not exist because~$\S_{i^*}$ is not controllable.
If~$\S_i$ satisfies (iii), then we set $i^* = i$ and according to Proposition~\ref{prop: equiv H feedback}, the system~$\S_{i^*}$ is  not  LSOPI.
Finally, assume that~$\S_i$ satisfies~(ii). Then $\mathcal{D}_i^{k_i}$ is the first noninvolutive distribution, all $\mathcal{D}_i^j$, for $\dleq 0 j k_i-1$,  are involutive, $\mrk \mathcal{D}_i^{k_i-1} = 2k_i$, and~$\S_i$ can be transformed
 (via a change of coordinates and an invertible static feedback) into $\dot {\bar x}^i = \bar f(\bar x^i, z^1), \dot z^1 = z^2, \dots, \dot z^{k_i} = v$, where  $z^j=(z_1^j, z_2^j)^{ \top}$ and  $v = (v_1, v_2)^{ \top}$.
 The ${\bar x}^i $-subsystem is nonlinear (since $\mathcal{D}_i^{k_i}$ is noninvolutive) of  dimension $\dim \bar x^i = n_i - 2k_i$, where  the state-space $X^i$ of~$\S_i$ is of dimension~$n_i$.
 Now, for $ \S_{i+1} = \widetilde \S_i^{(1,0)}$, evolving on $X^{i+1}$ {whose} $\dim X^{i+1} =n_{i+1}=  n_i + 1$ (see the structure of the prolongation $\widetilde \S^{(1, 0)}$),  the distribution $\mathcal{D}_i^{k_{i}}$ is involutive and of rank  $2(k_i+1)$,
so the dimension of the nonlinear subsystem of  $ \S_{i+1}$ satisfies
$\dim \bar x^{i+1} \leq    n_{i} + 1-  2(k_i+1)= n_i  - 2 k_i-1 $,
being thus smaller than  $\dim \bar x^i$.
So the sequence $\S_0, \S_1, \ldots,  \S_i$ terminates after a finite number of prolongations giving
$\S_{i^*}$ that satisfies either (i) or~(iii).
\end{proof}

%
%

 \begin{algorithm}
 \begin{algorithmic}[1]
 \REQUIRE $\S_0, \ldots, \S_i$
 \STATE \textbf{switch}~$\S_i$  \textbf{do}

\STATE \quad \textbf{case} {$k_i$ does not exist}: 
\STATE \quad \quad \textbf{if} $\rho$ exists \textbf{then}
\STATE \quad \quad \quad  {$\S_0$ is LSOPI.}
\STATE \quad \quad \textbf{else}
\STATE \quad \quad \quad {$\S_0$ is not flat (in particular, not LSOPI).}
\STATE \quad \quad \textbf{endif}

\STATE \quad \textbf{case} $k_i$   exists but a LSOPI-distribution~$\mathcal{H}_i^{k_i}$ does not:
\STATE \quad \quad  {$\S_0$ is not LSOPI.}
\STATE \quad \textbf{case} $k_i$   exists and there exists a unique LSOPI-distribution~$\mathcal{H}_i^{k_i}$:
\STATE \quad \quad  {Define $\S_{i + 1} = \widetilde \S_i^{(1,0)}$.}
\STATE \quad \textbf{case} $k_i$   exists, but the LSOPI-distribution~$\mathcal{H}_i^{k_i}$ is not unique:
\STATE \quad \quad  {The algorithm does not give an answer.}
\STATE \textbf{end switch}
\end{algorithmic}
\caption{Possible cases at iteration~$i$ for checking LSOPI}\label{alg iteration i}
 \end{algorithm}


  {In order to express the conditions of Theorem~\ref{thm: algo} in a verifiable  form, that is, to transform  Algorithm~\ref{alg iteration i} into constructive Algorithm~\ref{algo detailed} below}, we need to describe explicitly the construction of $\S_{i+1} = \widetilde \S_i^{(1, 0)}$, see {line 11} of Algorithm~\ref{alg iteration i} 
  summarizing all possible cases for LSOPI at iteration~$i$.
 For $i\geq 0$,
consider the control system
$$
\S_i : \dot x^i = f^i(x^i)+u_1^ig_1^i(x^i)+u_2^ig_2^i(x^i), \quad x^i\in X^i.
$$
To simplify the notation,
 we will drop the index~$i$ when dealing with $f, g_1, g_2, x$, and~$u$,  but keep it for~$\S_i$, $\mathcal{D}_i^j$ and~$\mathcal{H}_i^j$.  Assume that~$\S_i$ satisfies
(ii) of Theorem~\ref{thm: algo}.
 Then, by Proposition~\ref{prop: equiv H feedback}(i),
we compute
the distribution  $\red{\widetilde {\mathcal{G}}_2^p} = \mspan\{\red{\tilde  g_{2}^p}\}$ of {corank} one in $\mathcal{D}_i^0$, where $\red{\tilde  g_{2}^p} = \b_2^1 g_1 + \b_2^2 g_2$ is such that  $ \mathcal{H}_i^{k_i}  =  \mathcal{D}_i^{k_i-1} + \mspan\{ad_f^{k_i} \red{\tilde  g_{2}^p}\}$. 
Define the {to-be-prolonged control}
$\red{\t u_1^p} =  \b_2^2  u_1 - \b_2^1 u_2$, and choose any  feedback\footnote{
We can always take the first column of $\b$  as either  $\left(
\begin{array}{c}
1\\
0
\end{array}
\right)$ or $\left(
\begin{array}{c}
0\\
1
\end{array}
\right)$.} $u = \b \t u$ such that $\b = \left(
\begin{array}{cc}
* & \b_2^1\\
* & \b_2^2
\end{array}
\right)$, \red{where~$*$ stands for any functions that render $\b$ invertible}.
Apply it to~$\S_i$
to get $\red{\widetilde \S_i^p} : \dot x = f(x)+\red{\t u_1^p}\t g_1(x)+\t u_2\red{\tilde  g_{2}^p}(x)$, and prolong the first control~${\red{\t u_1^p}}$ (see Definition~\ref{def: prolongation}) to obtain
$$
\widetilde{\S}_i^{(1, 0)} :
\left \{
\begin{array}{lcl}
\dot x &= & f(x)+ \red{\t u_1^p} \tilde g_1(x)+v_2 \red{\tilde  g_{2}^p}(x) \\
\red{\dot {\t u}_1^p}&= &v_1,
\end{array}
\right.
$$
which is called \emph{the prolongation
of~$\S_i$ associated\footnote{From now on,  any construction of a prolongation~$\widetilde{\S}^{(1, 0)}$ will always be associated with a LSOPI-distribution $\mathcal{H}^{k}$.}
to~$\mathcal{H}_i^{k_i}$,}
and which we denote by $ \S_{i+1}$.
For $ \S_{i+1}$, {we have $\mathcal{D}_{i+1}^{j} = \mathcal{H}_i^{j} + \mspan\{\partial_{\red{\t u_1^p}}\}$, for $\dleq 0 j k_i$. In particular,} $\mathcal{D}_{i+1}^{k_i} = \mathcal{H}_i^{k_i} + \mspan\{\partial_{\red{\t u_1^p}}\}$ and is involutive and thus the first noninvolutive distribution  $\mathcal{D}_{i+1}^{k_{i+1}}$ of  $ \S_{i+1}$ satisfies $k_{i+1} > k_i$.

\begin{remark}\label{rk: to-be-prolonged control 1}
Notice that \red{given~$\mathcal{H}^k$, a to-be-prolonged control~$\red{\t u_1^p}$}  is not unique and
 given up to {an} affine transformation (i.e., we can take any $\red{\t u_1^p}  = \a_1 +  \gamma(\b_2^2  u_1 - \b_2^1 u_2)$, with~$\a_1$ arbitrary and {any function} $\gamma\neq 0$)
because
$\t g_2^{\red p}$ that spans $\widetilde {\mathcal{G}}_2^{\red p}$ is given up to a multiplicative function,
and does not depend on transformations $f \mapsto f + g\a $.
{Indeed, if instead of prolonging $\b_2^2  u_1 - \b_2^1 u_2$, we prolong $ \a_1 +  \gamma(\b_2^2  u_1 - \b_2^1 u_2)$ as above, the proposed algorithm will reach the same conclusions since, at each step,  the obtained systems are static feedback equivalent.}
{Actually, the next $(i+1)$st-step of the proposed algorithm takes care of the freedom in choosing $\a,\b_2^1,  \b_2^2$ 
at the previous~$i$th-step.}
Finding~$\red{\t u_{1}^{i,p}}$
at each step  {thus} requires knowing \cyan{the ration\footnote{\cyan{Where we use the usual notation from projective geometry.}} $[\b_2^1 : \b_2^2]$}  {only} (of {the step~$i$}), which in turn is reduced to calculate \cyan{a vector field $\red{\t g_{2}^p}$ that spans $\widetilde{\mathcal{G}}^p_2$} \red{of the step~$i$} and, finally, to construct~$\mathcal{H}_i^{k_i}$.
We treat that problem 
in the next section.
\end{remark}

\renewcommand{\algorithmiccomment}{\STATE //}
 \begin{algorithm}
 \begin{algorithmic}[1]
 \REQUIRE  $  \S$ (with associated vector fields $f$, $g_1, {g_2}$)

 \COMMENT {\textit{Initialization}}
  \STATE{Find the smallest integer $k$ such that~$\mathcal{D}^k$ is not involutive}
 \STATE{Set is{LSOPI} := True}
 \STATE{Set kExists := True}

 \WHILE{is{LSOPI} = True \AND kExists = True} 
 \IF{there exists a {LSOPI}-distribution~$\mathcal{H}^k$}

 \IF{the  {LSOPI}-distribution~$\mathcal{H}^k$ is unique}

  \STATE{Compute~$\mathcal{H}^k$}
  \COMMENT {\textit{For a given~$\mathcal{H}^k$, computation of the prolongation}}

 \STATE{Compute  $\t g_2^{\red p} = \b_2^1 g_1 + \b_2^2 g_2$ such that $ \mathcal{H}^k  =  \mathcal{D}^{k-1} + \mspan \{ad_f^k \t g_2^{\red p} \}$}
    \STATE{Compute  the {to-be-prolonged control} $\t u_1^{\red p}  =  \b_2^2  u_1 - \b_2^1 u_2$}
    \STATE{Compute the prolongation $\widetilde{\S}^{(1, 0)}$ associated to~$ \mathcal{H}^k$}

    \COMMENT {\textit{Verification whether  $\widetilde{\S}^{(1, 0)}$ is static feedback linearizable}}

    \IF{$\widetilde{\S}^{(1, 0)}$ verifies condition  (FL1) of Theorem~\ref{thm F-lin systems}}
     \IF{{$\rho$ exists}}
         \COMMENT {\textit{$\widetilde{\S}^{(1, 0)}$ is static feedback linearizable, all its linearizability distributions are involutive and thus {$k$ does not exist}}} 
             \STATE{Set kExists := False}
      \ELSE
         \COMMENT {{\textit{$\widetilde{\S}^{(1, 0)}$ is not controllable, thus not flat and not LSOPI}}}
           \STATE{{Set is{LSOPI} := False}}
      \ENDIF

   \ELSE

    \COMMENT {\textit{$\widetilde{\S}^{(1, 0)}$ is not static feedback linearizable}}
    \STATE{Set~$\S$ := $\widetilde{\S}^{(1, 0)}$}
    \STATE{Compute $k$}
    \ENDIF
\ELSE

  \COMMENT {\textit{{The LSOPI-distribution~$\mathcal{H}^k$ is not unique thus the algorithm does not give an answer}}}
 \STATE{Set~is{LSOPI} := Not conclusive}
 \ENDIF
 \ELSE

 \COMMENT {\textit{{The LSOPI-distribution~$\mathcal{H}^k$ does not exist}}}
  \STATE{Set is{LSOPI} := False}

   \ENDIF
  \ENDWHILE
 \ENSURE is{LSOPI}
\end{algorithmic}
\caption {Checking {LSOPI}}
\label{algo detailed}
 \end{algorithm}

\color{black}

\subsection{Construction of {the LSOPI-distribution~$\mathcal{H}^{k}$}}
\label{ssec: construction H}

{Next}
we study properties {of~$\mathcal{H}^k$ and present its possible
constructions.} 
Recall that $k$ is the smallest integer such that~$\mathcal{D}^k$ is not involutive. 
{For a distribution $\mathcal{D}$, we denote by $\overline{\mathcal{D}}$ its involutive closure.}

\color{black}
\begin{proposition}\label{prop Dk noninvolutive}
 Consider the system~$\S$, given by~\eqref{eq: Sigma}.
  \color{black}
 Then one of the following cases holds:

\noindent  either \emph{{Case~I:}} $\mathcal{D}^j = \overline{\mathcal{D}^j}$, for $\dleq 0 j n$, i.e., $k$ does not exist,

\noindent or $k$ exists and then
 \color{black}

\noindent either \emph{{Case~II:}}  \hspace{0.1cm}$[\mathcal{D}^{k-1},\mathcal{D}^k] \not\subset \mathcal{D}^k$,


 or \emph{{Case~III}:} $[\mathcal{D}^{k-1},\mathcal{D}^k] \subset \mathcal{D}^k$  but $[\mathcal{D}^{k},\mathcal{D}^k] \not\subset \mathcal{D}^k$,  with the {former} condition being
%
  empty if $k= 0$.

 Moreover, if {$k$ exists and} there exists a LSOPI-distribution~$\mathcal{H}^k$ in~$\mathcal{D}^k$,  then we necessarily have
\begin{equation}\label{eq: r}
r = \mcork(\mathcal{D}^k \subset \mathcal{D}^k + [\mathcal{D}^k,\mathcal{D}^k])\leq 2,
\end{equation}
 and if, additionally~$\mathcal{D}^k$ is in {Case~II,} then
\begin{equation}\label{eq: r I}
{r_{{}_\mathrm{II}}} =\mcork(\mathcal{D}^k \subset \mathcal{D}^k + [\mathcal{D}^{k-1},\mathcal{D}^k]) = 1.
\end{equation}
\end{proposition}

{In \red{Remark~\ref{rk: cases LSOPI-distribution}}, {see Section~\ref{ssec: algo},} we distinguished different cases depending on the existence and uniqueness of~$\mathcal{H}^k$. Obviously, Case~I {of Proposition~\ref{prop Dk noninvolutive}} corresponds to condition~(i) of \red{Remark~\ref{rk: cases LSOPI-distribution}}. The above proposition asserts that if $k$ exists, then we are either in Case~II or in Case~III. Now we prove that  Case~II corresponds to~(ii) and Case~III to {(iv)} of  \red{Remark~\ref{rk: cases LSOPI-distribution}}.}
We would like to emphasize that for any system~$\S$ {for which $k$ exists}, {due to the above proposition}, the first noninvolutive distribution~$\mathcal{D}^k$ is either in {Case~II} or in {Case~III}.
Observe that a necessary condition for a system to be LSOPI is that $r\leq 2$ ({i.e.,} at most two new directions stick out {from~$\mathcal{D}^k$} when calculating $\mathcal{D}^k+[\mathcal{D}^k, \mathcal{D}^k]$).
{In the proof of  Proposition~\ref{prop Dk noninvolutive} we show that if $k$ exists, then always $[\mathcal{D}^{k-2},\mathcal{D}^k] \subset \mathcal{D}^k$. It follows that}
in {Case~II,} {the only new direction} is given by a vector field of the form $[ad_f^{k-1} g_i, ad_f^{k} g_j]$, $\dleq 1 {i,j} 2$.
For {Case~II,} we explicitly construct in Proposition~\ref{prop: H case II} below the (unique) LSOPI-distribution~$\mathcal{H}^k$ {described by} Proposition~\ref{prop: equiv H feedback}. 
{Case~III} is treated in  Proposition~\ref{prop: H case III}, which is completed by
{the Appendices~A and~B}. 

For  a  distribution $\mathcal{D}$  defined on a manifold $X$,
we denote by $\mathcal{C}(\mathcal{D})$ its characteristic distribution, i.e., 
$\mathcal{C}(\mathcal{D}) = \{ \xi \in \mathcal{D} : [ \xi,\mathcal{D}]\subset \mathcal{D}\}$.
It follows directly from the Jacobi identity that the characteristic distribution is always involutive.

\begin{proposition}[{Case~II, explicit construction of~$\mathcal{H}^k$;  Case~III}]\label{prop: H case II}
{\it Consider 
$\S: \dot x = f(x)+u_1g_1(x)+u_2g_2(x)$ and \red{let~$k$ be its non-involutivity index.}
{There exists a {LSOPI}-distribution~$\mathcal{H}^k$
if and only if}
 \begin{enumerate} \item[-]

 either {\normalfont Case~II} holds {with $r=2$} and, {additionally},
 \begin{enumerate}[\normalfont {(A}1)]
  \item ${r_{{}_\mathrm{II}}} = 1;$
 \item $\mrk \mathcal{C}(\mathcal{D}^k) = 2k-1$;
 \item $\mrk  (\mathcal{D}^{k-1} +  [f, \mathcal{C}(\mathcal{D}^k)]) = 2k+1$;
 \item $\mathcal{D}^{k-1} +  [f, \mathcal{C}(\mathcal{D}^k)]$ is involutive;
 \end{enumerate}
{called \nf{subCase {II}\,A}};
\item [-]
{or {\normalfont Case~II}} holds {with $r=1$} and, {additionally},
 \begin{enumerate}[\normalfont {(B}1)]
    \item $\mrk \mathcal C(\mathcal{D}^k) = 2k$;
    \item $\mrk  (\mathcal{D}^{k-1} + \mathcal C(\mathcal{D}^k))= 2k+1$;
    \item  $\mathcal{D}^{k-1} + \mathcal C(\mathcal{D}^k)$ is involutive;
\end{enumerate}
{called \nf{subCase {II}\,B}};
\item [-] 
or {\normalfont {Case~III}} holds.
\end{enumerate}

Moreover, in \nf{{Case~II}}, if a LSOPI-distribution~$\mathcal{H}^k$ exists, then it is unique and given by
$$
\mathcal{H}^k=\mathcal{D}^{k-1} +  [f, \mathcal{C}(\mathcal{D}^k)] \mbox{ in  \nf{subCase~II\,A}},
$$
and by
$$
\mathcal{H}^k = \mathcal{D}^{k-1} + \mathcal C(\mathcal{D}^k)\mbox{ in  {\normalfont subCase~II\,B}.}
$$
 In  {\normalfont {Case~III}},~$\mathcal{H}^k$ is never unique. Furthermore, \cyan{for any}
given~$\mathcal{H}^k$, unique in  \nf{{Case~II}}  or any in  \nf{{Case~III}}, we can explicitly identify the distribution  $\red{\widetilde {\mathcal{G}}_2^p} = \mspan\{\red{\t g_2^p}\}$ of {corank} one in~$\mathcal{D}^0$, where~$\red{\t g_2^p}$
is such that $ \mathcal{H}^k  =  \mathcal{D}^{k-1} + \mspan\{ad_f^k \red{\t g_2^p}\}$.}
\end{proposition}

\begin{center}
\begin{figure}
 \includegraphics[scale=0.8]{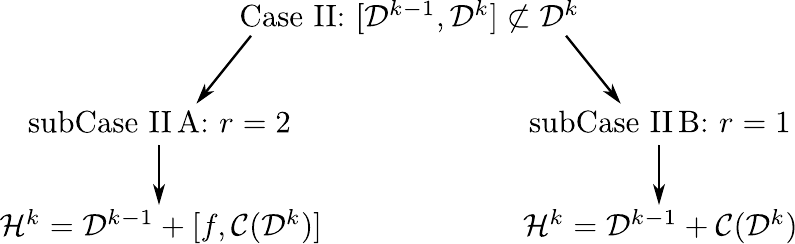}
 \caption{Construction of~$\mathcal{H}^k$ by Proposition~\ref{prop: H case II}.}
 \label{fig:caseII}
\end{figure}

\end{center}

Conditions (A1)-(A4) and (B1)-(B3) describe, resp., the subCases II\,A and II\,B  and characterize in an invariant way\footnote{{According to the proof of Case II\,A, we have $ \mathcal{C}(\mathcal{D}^k) \subset \mathcal{D}^{k-1}$, implying that $\mathcal{D}^{k-1} +  [f, \mathcal{C}(\mathcal{D}^k)]$ is indeed {feedback invariant and well defined}.}} all situations when~$\mathcal{D}^k$ of Case~II admits a LSOPI-subdistribution $\mathcal{H}^k \subset \mathcal{D}^k$, \cyan{see also Figure~\ref{fig:caseII}}. The characterization, being geometric, implies that~$\mathcal{H}^k$ (if it exists) is unique. \red{Thus, Theorem 2 applies to systems $\S$ that satisfies Case II and either the set of conditions (A), corresponding to subCase II\,A, or the set of conditions (B), corresponding to subCase II\,B, and moreover, for which all prolongations in the sequence $\S = \S_0, \S_1, \ldots, \S_{\ell-1}$ also verify Case II and either (A) or~(B).}

If {Case~III} holds, then  {according to} Proposition~\ref{prop: H case II}, {many involutive  subdistributions of the form} $\mathcal{H}^k = \mathcal{D}^{k-1} + \mspan\{ad_f^k \red{\t g_2^p}\}$ exist, but
not all   choices of~$\mathcal H^k$ lead to a conclusive answer whether the system is {LSOPI} or not.
 \red{The main difficulty is that for each system $\S_i$ for which $\mathcal{H}_i^{k_i}$ is not unique, we would need to verify whether the sequence constructed from each LSOPI-distribution $\mathcal{H}_i^{k_i}$ satisfies the LSOPI conditions. However, this cannot be done algorithmically, as the family of all involutive subdistributions~$\mathcal{H}^k$ of corank one in~$\mathcal{D}^k$ is parameterized by a functional parameter.
 It is therefore natural to ask whether, among the infinitely many LSOPI-distributions in Case III, only some lead to sequences of prolongations satisfying the LSOPI conditions, and how to distinguish them.
 }
We explain in Proposition~\ref{prop: H case III} {below}
how to
{identify
{a right} LSOPI-distribution~$\mathcal{H}^k$}
for {all subcases of {Case~III}   for which the system may be LSOPI (except for one very particular subcase for which we propose a  natural choice for~$\mathcal{H}^k$ but without guaranteeing that the algorithm is conclusive with that choice),
and show that for all remaining subcases the system is never LSOPI.}
%

In Proposition~\ref{prop: H case III} below, we use the notion of {growth} vector that we define as follows.
To a distribution $\mathcal{E}$, we associate the following sequence of distributions:
\begin{equation}\label{eq: def sequence  G}
 \mathcal{E}^0 = \mathcal{E} \mbox{ and } \mathcal{E}^{j+1} = \mathcal{E}^{j} + [\mathcal{E}^{j},\mathcal{E}^j], \, j\geq 0,
\end{equation}
and define {$r^j  = \mrk \mathcal{E}^j$, $j\geq 0$. The vector $(r^0, r^1, r^2, \ldots)$} 
{is}
called the {{growth} vector of $\mathcal{E}$.}


\begin{proposition}[{Case~III, explicit construction of~$\mathcal{H}^k$}]\label{prop: H case III}
 Consider $\S: \dot x = f(x)+u_1g_1(x)+u_2g_2(x)$, given by~\eqref{eq: Sigma}. Assume that its first noninvolutive distribution~$\mathcal{D}^k$ is in \nf{{Case~III}}. {Then} one of the following mutually exclusive conditions holds:
 \begin{enumerate}[\normalfont {(C}1)]

\item  {$\mrk \overline{\mathcal{D}}^k= 2k+3$ and $\overline{\mathcal{D}}^k = TX$}; 
\item the  {{growth}} vector of~$\mathcal{D}^k$ starts with $(2k+2, 2k+3,2k+4, \ldots)$;

\item the  {{growth}} vector of~$\mathcal{D}^k$ starts with $(2k+2, 2k+3,2k+ 5, \ldots)$;

\item  $\mrk \overline {\mathcal{D}}^k= 2k+3$, {$\overline {\mathcal{D}}^k\neq TX$},  and
$\mrk (\overline {\mathcal{D}}^k+ [f,\mathcal{D}^k]) = 2k+3$.

  \item  $\mrk \overline {\mathcal{D}}^k= 2k+3$, and
$\mrk (\overline {\mathcal{D}}^k+ [f,\mathcal{D}^k]) = 2k+4$, implying the existence of a non-zero vector field\footnote{\cyan{We denote it by $\t g_{2}$ instead of $\t g_{2}^p$ because only in some cases it defines a to-be-prolonged control.}} $\t g_{2} \in \mathcal{D}^0$ such that 
$ad_f^{k+1} \t g_{2} \in  \overline {\mathcal{D}}^k$ and
 \begin{list}{}{}
 \item  either \normalfont ({C5})$'$  $ad_f^{k+1} \t g_{2} \in  \mathcal{D}^k$, i.e., $\mrk \mathcal{D}^{k+1} = 2k+3$;
  \item  \hspace{0.37cm} or \normalfont ({C5})$''$  $ad_f^{k+1} \t g_{2} \not\in  \mathcal{D}^k$, i.e., $\mrk \mathcal{D}^{k+1} = 2k+4$;
 \end{list}

\item  $\mrk \overline {\mathcal{D}}^k= 2k+3$ and
$\mrk (\overline {\mathcal{D}}^k+ [f,\mathcal{D}^k]) = 2k+5$;

 \end{enumerate}

Moreover,
the following statements hold:
\begin{enumerate}[a)]
 \item If~$\S$ 
 satisfies \nf{({C1})}, then \cyan{there exists infinitely many LSOPI-distributions, the associated} prolongation $\widetilde {\S}^{(1,0)}$ can be  defined by  any LSOPI-distribution~$\mathcal{H}^k$ and
  $\widetilde {\S}^{(1,0)}$ is static feedback linearizable, implying that~$\S$ is LSOPI.

 \item If~$\S$ is LSOPI  and satisfies \nf{({C2})}, then the prolongation $\widetilde {\S}^{(1,0)}$ has to be defined by   $\mathcal{H}^k = \mathcal{C}(\mathcal{D}^k + [\mathcal{D}^{k}, \mathcal{D}^k])$.

 \item If~$\S$ satisfies  {either \nf{({C3})}, or \nf{({C4})}, or   \nf{({C5})} with  \nf{({C5})$'$},} or~\nf{({C6})}, then~$\S$ is not LSOPI.
 \end{enumerate}

\end{proposition}

The second part of the above proposition \cyan{describes} all cases except for ({C5}) with ({C5})$''$, i.e., $\mrk \mathcal{D}^{k+1} = 2k+4$ (see ({C5})$''$ above); {the latter}  contains subcases and is studied  in {Appendices~{A}} {and~B}.

\begin{center}
\begin{figure}[h]
 \includegraphics[scale=0.72]{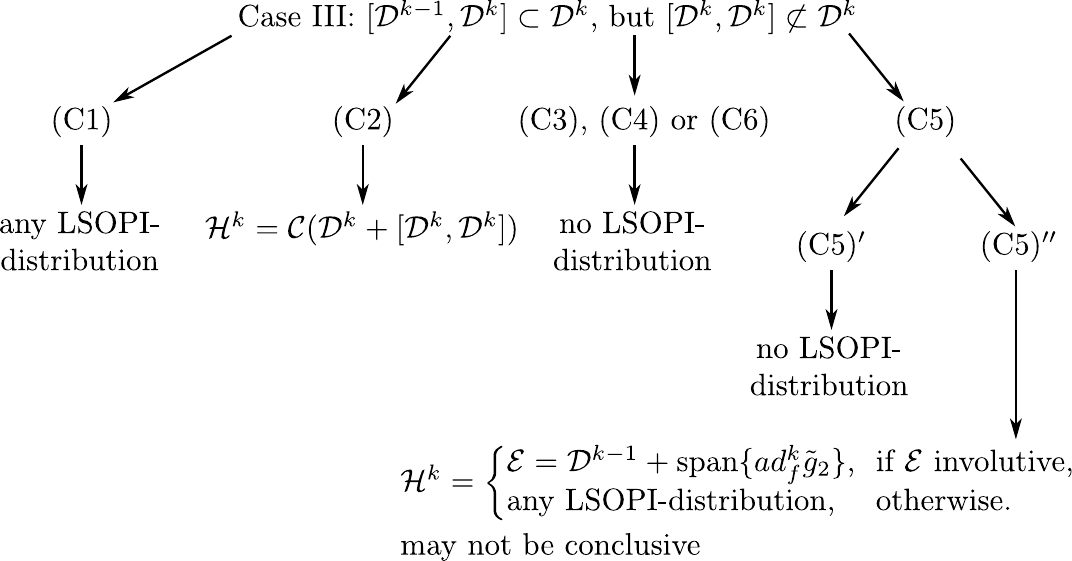}
 \caption{Construction of~$\mathcal{H}^k$ by Proposition~\ref{prop: H case III} and relation~\eqref{eq: H case C5}.}
 \label{fig:caseIII}
\end{figure}

\end{center}

\color{black}
A particular subcase of {Case~III} is when the {{growth}} vector of~$\mathcal{D}^k$ starts with $(2k+2, 2k+3,2k+4, \ldots)$ as considered in {({C2})}. In this scenario,
there exists a choice of~$\mathcal{H}^k$ that is more natural than others. Specifically, the distribution~$\mathcal{D}^k$ contains (with an inclusion of corank one) the characteristic distribution of $\mathcal{D}^k + [\mathcal{D}^{k}, \mathcal{D}^k]$, which is involutive by construction (and also contains~$\mathcal{D}^{k-1}$).
Therefore, it is natural to define~$\mathcal{H}^k$ as $\mathcal{H}^k = \mathcal{C}(\mathcal{D}^k + [\mathcal{D}^{k}, \mathcal{D}^k])$. Moreover, according to statement~{b)} of Proposition~\ref{prop: H case III}, any other choice of~$\mathcal{H}^k $ will {actually} lead to a prolongation for which  no LSOPI-distribution exists. 
%
Similarly, if ({C5}) holds, then a particular vector field
${\t g_{2}}\in \mathcal{D}^0$
can be distinguished.
If
${\t g_{2}}$  {additionally} satisfies ({C5})$''$ and the distribution
%
%
$\mathcal{E} = \mathcal{D}^{k-1} + \mspan\{ad_f^k {\t g_{2}}\}$
is involutive, we show in Appendix~{A}
that  its associated prolongation ${\widetilde \S^{(1,0)}}$ is the only one for which~$\mathcal{D}_p^{k+1}$ is also involutive (that is, we gain at least {two new}  involutive distributions for the prolongation).  Thus in the case (C5) with  (C5)$''$, we propose to construct a LSOPI-distribution~$\mathcal{H}^k$ by
\begin{equation}\label{eq: H case C5}
\mathcal{H}^k =
\left\{\begin{array}{ll}
         \mathcal{E}  = \mathcal{D}^{k-1} + \mspan\{ad_f^k {\t g_{2}}\},& \mbox{ if $  \mathcal{E}$ involutive,}  \\
          \mbox{any LSOPI-distribution},& \mbox{ otherwise.}
         \end{array}\right.
\end{equation}
\cyan{Subcases of Case III described by Proposition~\ref{prop: H case III} as well as the above comment are summarized in Figure~\ref{fig:caseIII}.}

\subsection{\red{Constructive LSOPI-algorithm}} \label{ssec: algo constructive}

Combining Theorem~\ref{thm: algo} and Propositions~\ref{prop: H case II} and \ref{prop: H case III} leads to the {following main result,
where the sets of conditions (A1)-(A4) and (B1)-(B3) are denoted by~(A) and (B), resp..} 

\begin{theorem}\label{thm: algo constructive}
 Consider the system $\S: \dot x = f(x)+u_1g_1(x)+u_2g_2(x)$ and assume that a sequence of prolonged systems $\S= \S_0, \S_1, \ldots,\S_{i-1},  \S_i$ has been constructed.
 Then {for~$\S_i$:}
 \begin{enumerate}[\normalfont (i)]
  \item either {$k_i$ does not exist},
  that is,  all linearizability distributions $\mathcal{D}_i^j$  are involutive;
  \item [\nf{(ii)}$'$] or  $k_i$ exists and~$\S_i$ satisfies   either \nf{{(A)}} or  \nf{{(B)}} or  {\nf{({C1})} or \nf{({C2})}},
  {which} uniquely identify the LSOPI-distribution $\mathcal{H}_i^{k_i}$ via Proposition~\ref{prop: H case II} for \nf{{(A)}} and  \nf{{(B)}}, {or}    via Proposition~\ref{prop: H case III} for {\nf{({C1})} and \nf{({C2})}}, and then set $\S_{i + 1} =\widetilde \S_i^{(1,0)}$;
  \item [\nf{(iii)}$'$]  or  $k_i$   exists  but either a LSOPI-distribution~$\mathcal{H}_i^{k_i}$ does not exist or~$\S_i$ satisfies {either} \nf{(C3)}, or \nf{({C4})},  {or \nf{({C5})} with \nf{({C5})$'$}}, or \nf{({C6})};

   \item [\nf{(iv)}$'$] {or $k_i$   exists and~$\S_i$ satisfies \nf{(C5)} with  \nf{(C5)$''$}   and then define {$\mathcal{H}_i^{k_i}$ by~\eqref{eq: H case C5}}
   and set  $\S_{i + 1} = \widetilde \S_i^{(1,0)}$;}
 \end{enumerate}
 and the sequence $\S_0, \S_1, \ldots,  \S_i$ always terminates at, say,~$\S_{i^*}$ which satisfies either {\normalfont (i)} or {\normalfont (iii)$'$}.

 If~$\S_{i^*}$ satisfies \nf{(i)} {and there exists $\rho$ such that $\mathcal{D}_{i^*}^{\rho} = TX^{i^*}$, where $X^{i^*}$ is the state-space of~$\S_{i^*}$}, then~$\S_{i^*}$ is static feedback linearizable, implying that~$\S$ is LSOPI;
  {if~$\S_{i^*}$ satisfies \nf{(i)} and  $\rho$ does not exist, then~$\S$ is not flat (in particular, not  LSOPI).}

  If~$\S_{i^*}$ satisfies  {\normalfont (iii)$'$}  and {if at each step $\dleq 0 { i} i^*-1$, the prolongation~$\S_{ i+1}$ was constructed as in \nf{(ii)$'$},}   then~$\S$ is not LSOPI. If~$\S_{i^*}$ satisfies  {\normalfont (iii)$'$} and there exists   $\dleq 0 { i} i^*-1$ such that the prolongation $\S_{ i+1}$ was constructed as in \nf{(iv)$'$}, then the algorithm does not give a conclusive answer.

\end{theorem}

 On one hand, Theorem~\ref{thm: algo constructive} allows to weaken assumption~(ii) of Theorem~\ref{thm: algo} {due to replacing it by (ii)$\red '$ and completing it with (iv)$'$, that is,}  by adding many cases, where~${\mathcal{H}_i^{k_i}}$ is not unique. On the other hand, it provides verifiable conditions to check all assumptions. \red{Observe that none \cyan{among the subcases (C1),\ldots, (C6) has been} covered by Theorem~\ref{thm: algo} and  moreover, in the formulation of condition (ii)$'$ of Theorem~\ref{thm: algo constructive}, we also specify how the LSOPI-distribution must be constructed (while only its existence was required in the statement of Theorem~\ref{thm: algo}).}
 {The algorithm of constructing~$\S_{i^*}$, {based on Theorem~\ref{thm: algo constructive},} is almost always constructive,  {namely, if} at each step~$i$, the corresponding system~$\S_i$ satisfies either (i) or (ii)$\red '$ or (iii)$\red '$.
 If, however, at a certain step~$i$, the corresponding system~$\S_i$ fulfils  (iv)$'$, that is, (C5) with (C5)$''$, and the algorithm terminates with~$\S_{i^*}$ satisfying  (iii)$\red '$, then we cannot conclusively determine whether the system~$\S$ is not LSOPI or a wrong choice of a LSOPI-distribution~$\mathcal{H}_i^{k_i}$ was made at the step~$i$ (and another choice, {say~$\widetilde{\mathcal{H}}_i^{k_i}$}, could have allowed to continue the algorithm).  We will distinguish and characterize special distributions~$\mathcal{H}^k$ in Appendix~A and describe subcases {of the case (C5) with (C5)$''$}, {for which} the algorithm is conclusive, in Appendix~B.}

Algorithm~\ref{algo detailed constructive} {summarizes}  the details of the LSOPI process {based on Theorem~\ref{thm: algo constructive};  in particular, it distinguishes  different casses of constructing the right candidate for} the LSOPI-distribution~$\mathcal{H}^k$. 

\renewcommand{\algorithmiccomment}{\STATE //}
 \begin{algorithm}
 \begin{algorithmic}[1]
 \REQUIRE  $ \S$ (with associated vector fields $f$, $g_1,  {g_2}$)
  \COMMENT {\textit{Initialization}}
  \STATE{Find the smallest integer $k$ such that~$\mathcal{D}^k$ is not involutive}
 \STATE{Set isLSOPI := True  }
 \STATE{Set kExists := True}
\STATE{Set isConclusive := True}

 \WHILE{isLSOPI = True \AND kExists = True}
 \IF{there exists a  {LSOPI}-distribution~$\mathcal{H}^k$}

 \IF{either {(A)} or {(B)} or {({C1})} or {({C2})} 
 is satisfied}

  \STATE{Compute~$\mathcal{H}^k$ via  Proposition~\ref{prop: H case II} for {(A)} and  {(B)}, resp., via Proposition~\ref{prop: H case III} for~{({C1})} {and} {({C2})} 
  }

  \STATE{Perform instructions of lines 10--{26} of Algorithm~\ref{algo detailed}}

\ELSIF{{(C5) with (C5)$''$} is satisfied}

 \STATE{Set~isConclusive := False }

    \STATE{{Compute~$\mathcal{H}^k$  by~\eqref{eq: H case C5}}}


   \STATE{Perform instructions of lines 10--{26} of Algorithm~\ref{algo detailed}}

\ELSE

  \COMMENT {\textit{Either {\normalfont {({C3})}} or   {\normalfont {({C4})}} {or \nf{({C5})} with \nf{({C5})$'$}} {or   {\normalfont {({C6})}}}   is satisfied}}
   \STATE{Set isLSOPI := False}

 \ENDIF
 \ELSE

 \COMMENT {\textit{{The LSOPI-distribution~$\mathcal{H}^k$ does not exist}}}
  \STATE{Set isLSOPI := False}

   \ENDIF
  \ENDWHILE
\IF{isConclusive = False  \AND isLSOPI = False }
    \STATE{Set answerLSOPI := Not conclusive }
\ELSE
  \STATE{Set answerLSOPI := isLSOPI}
  \ENDIF

 \ENSURE answerLSOPI
\end{algorithmic}
\caption {Checking {LSOPI}}
\label{algo detailed constructive}
 \end{algorithm}


\section{\red{Interplay between flatness, LSOPI, LSOP {and L$\ell$P} properties for two-input control systems}}\label{sec:interplay}

In this section, we will discuss relations between flatness, LSOPI-systems and LSOP-systems.
Fix the state-dimension \( n \) and set \( m = 2 \).
Denote by $\mathcal{FLAT}$ the class of all  flat, by
$\mathcal{LSOPI}$ the class of all LSOPI-systems, see Definition~\ref{def: LSOPI}, by
$\mathcal{LSOP}$ the class of all LSOP-systems, see  Remark~\ref{rk:LSOP}, and by
$\mathcal{L}\ell\mathcal{P}$ the class of all \green{L$\ell$P-systems, that is, systems} linearizable by an \( \ell \)-fold prolongation.
{For the control system~$\S$,
the latter {property} means the existence}
around~$x_0$ of an invertible static feedback transformation $u = \alpha(x) + \beta(x)\tilde u$ bringing~$\Sigma$  into 
$\widetilde{\Sigma}: \dot x =\tilde  f(x)+\tilde  u_1 \tilde  g_1(x) +\tilde   u_2\tilde  g_2(x)$, and for which the prolonged system
\begin{equation}
\label{eq: p-fold prolongation}
\widetilde {\Sigma}^{(\ell, 0)} :
\left \{
\begin{array}{l@{\,}cl}
\dot x &= &\tilde f(x)+ \t u_1^0 \tilde g_1(x)+v_2\tilde  g_2(x), \\
\dot {\t u}_1^j &= &{\t u}_1^{j+1}, \quad {\dleq 0 j \ell-2,}\\
\dot {\t u}_1^{\ell-1}&= &v_1,
\end{array}
\right.
\end{equation}
for a certain $ \ell \geq  1$,
is locally static feedback linearizable.
Observe that for linearizability via an invertible $\ell$-fold prolongation we need only one feedback transformation ({applied to} the original system~$\S$), while for LSOPI-systems or, more generally, for LSOP-systems,
we apply an invertible feedback transformation for each  system~$\Sigma_i$ of the sequence $\Sigma = \S_0, \S_1, \ldots, \S_\ell$.
Therefore, linearizability via an invertible $\ell$-fold prolongation can be considered as a subcase of LSOP\footnote{
If~$\S$ is  dynamic feedback linearizable via an invertible $\ell$-fold prolongation,  then~$\S$
 is also LSOP.
Indeed, for $\dleq 0 i \ell-1$, we clearly have $\Sigma_{i+1} = \widetilde {\Sigma}^{(i+1, 0)}$ which is simply the prolongation obtained by prolonging~$v_1$ of
$$ \S_i = \widetilde {\Sigma}^{(i, 0)} :
\left\{\begin{array}{l@{\,}c@{\,}l}
        \dot x&=& \tilde f(x)+\t u_1^0  \tilde g_1(x)+v_2\tilde  g_2(x)\\
        \dot {\t u}_1^j &=&  \t u_1^{j+1}, \quad \dleq 0 j i-2, \\
        \dot {\t u}_1^{i-1} &=& v_1.
       \end{array}
\right.
$$
}.

For a flat control system~$\S$ 
and \green{its} flat {output}~$\varphi$, whenever necessary to specify the number of derivatives of~$u$ on which the components of $\varphi$ depend, we {say} that the system~$\S$   is $(x, u, \ldots, u^{(q)})$-flat if the $q$th-derivative \green{of $u$} is the highest involved {in $\vp$ and $\S$ does not admit any flat output depending only on lower derivatives of $u$}. In the particular case $\varphi_i =\varphi_i(x)$, for $\dleq 1 i m$, we {say} that the system is $x$-flat.

It can be easily shown that for
$\S$, \green{with $m=2$,} $x$-flatness
 actually {reduces} to dynamic linearizability  via an~$\ell$-fold prolongation of a suitably chosen control, for a certain~$\ell$.
Moreover, it has been proved in \cite{gstottner2020linearization, gstottner2022flat} that the above property is also valid for $(x,u)$-flat systems with two inputs.
\green{We have the following relations between various notions of feedback linearization.}


\begin{theorem}\label{thm: equiv flat LSOP} For systems with \( m = 2 \), the following relations hold:
\begin{enumerate}[\normalfont(1)]

\item \( \mathcal{LSOPI} \subset \mathcal{LSOP} \) and the inclusion is strict;
    \item \( \mathcal{FLAT} = \mathcal{LSOP} \);
 \item $\mathcal{L}\ell\mathcal{P} \subset \mathcal{FLAT}$.
\end{enumerate}
\end{theorem}

\begin{remark}
 Relations (1) and (2) hold for any $m$ but we state and discuss them for $m=2$ since throughout the paper all systems satisfy $m=2$. As proved in~\cite{gstottner2020linearization,gstottner2022flat}, the class of $x$-flat and $(x,u)$-flat systems coincide with $\mathcal{L}\ell\mathcal{P}$ and it is an open question whether it holds for $(x, u, \ldots, u^{(q)})$-flat systems, that is, whether (3) is actually the equality.
\end{remark}

We present the proof of Theorem~\ref{thm: equiv flat LSOP} below (rather than in Section~\ref{sec_proofs}) since it  provides a clearer understanding of the interplay between the aforementioned properties.
%

\color{black}

\begin{proof}
\red{{\it Proof of} (1). It is obvious that \cyan{$\mathcal{LSOPI}$ forms a subclass of $\mathcal{LSOP}$} (recall {that} for LSOP-systems
we may have  $k_i \geq k_{i+1}$, see Definition~\ref{def: LSOPI} and Remark~\ref{rk:LSOP}). Example~\ref{ex: not necessary conditions} below  shows that there are systems that are LSOP without being LSOPI, thus the inclusion is \cyan{indeed} strict.}

\red{{\it Proof of} (2). It is clear that \cyan{being LSOP} implies flatness.
Let us now consider a flat control system of the form
$\S: \dot x = f(x)+u_1g_1(x)+u_2g_2(x)$.}

Suppose  that $\S$ is $(x, u, \ldots, u^{(q)})$-flat\footnote{
\red{See \cite{gstottner2022flat} for an example of a  two-input system that is flat \cyan{(more precisely, $(x,u,\dot u)$-flat)} without being $x$- or $(x,u)$-flat.}
}, with $q\geq 1$ (\cyan{the cases corresponding to $x$- and $(x,u)$-flatness follow from \cite{gstottner2020linearization,gstottner2022flat}}),
and let $(\varphi_1,\varphi_2)$, where   $\varphi_i = \varphi_i(x, u, \ldots, u^{(q)})$ for $i=1, 2$, be  a flat output of $\S$.
Prolong each control $(q+1)$ times
to get the following control system (where there is no tilde on $\S^{(q+1,q+1)}$ since no preliminary feedback transformation is applied):
$$
\Xi = \S^{(q+1,q+1)}: \left\{
\begin{array}{lcl}
\dot x &= &f(x)+u_1^0g_1(x)+u_2^0g_2(x)\\
\dot u_i^j&= &u_i^{j+1}, \quad \dleq 0 j q-1, \\
\dot u_i^q&= &v_i,  \quad \dleq 0 i 2.
\end{array}
\right. 
$$
 that we denote by~$\Xi$. The states and controls  of  $\Xi$ are, respectively,  $\xi = (x^\top, u^0,$ $ \ldots, u^q)^\top\in X \times \mathbb{R}^{2(q+1)}$ and  $v =(v_1, v_2)^\top\in \mathbb{R}^2$.
For $\Xi$, the pair  $(\varphi_1,\varphi_2)$ is  a flat output   depending on the state $\xi = (x^\top, u^0, \ldots, u^{q})^\top$ only (that is, an $x$-flat output).
Let~$\rho_1$ and $\rho_2$ denote the relative degrees of  $\varphi_1$ and~$\varphi_2$, that is, $\vp_1^{(\rho_1)}$ and $\vp_2^{(\rho_2)}$ are the lowest time-derivatives depending explicitly on~$v$. Consider the decoupling matrix
given by $D = (L_{G_j}L_F^{\rho_i-1}\varphi_i)_{\dleq 1 {i,j} 2}$, where $F$, $G_1$ and $G_2$ denote the drift and control vector fields of $\Xi$. 
We have $1\leq \mrk  D(\cyan \xi) \leq 2$ and according to Assumption~\ref{assum: constant rk and acces}, $ \mrk  D(\cyan \xi) $ is constant in a neighborhood of~$\cyan \xi_0$. 
If $\mrk  D(\cyan \xi) = 2$, then flatness implies that the system $\Xi$ is necessarily static feedback linearizable (i.e., $\rho_1 + \rho_2 = n + 2(q+1)$). Moreover, since each control variable has been prolonged the same number of times, it can be easily shown that the original system $\S$ is also static feedback linearizable, and thus actually $x$-flat, contradicting the $(x, u, \ldots, u^{(q)})$-flatness assumption.
Therefore, $\mrk  D(\cyan \xi) = 1$. Permute $\vp_1$ and $\vp_2$, if necessary, and put
$$
\vp_1^{(\rho_1)} = L_F^{\rho_1} \varphi_1  + (L_{G_1}L_F^{\rho_1-1}\varphi_1) v_1 + (L_{G_2}L_F^{\rho_1-1}\varphi_1)v_2 = \t v_1.
$$
Notice that $\t v_1 = \vp_1^{(\rho_1)} = \psi(x, u^0, \ldots, u^q, v)$ and that $\vp_2^{(\rho_2)}$ depends explicitly on $\t v_1$. Following, for instance~\cite{gstottner2020linearization}, we deduce that it is the control $\t v_1$ that has to be prolonged $\ell$-times   in order to dynamically linearize $\Xi$ (with the value of  $\ell$  that can be explicitly computed in \cyan{terms} of $\rho_1$, $\rho_2$, and $\dim \xi$).

In conclusion, in  order to construct a static feedback \cyan{linearizable} system $\widetilde \Xi^{(\ell, 0)}$,  we have applied the  transformations summarized by the following diagram:
$$\begin{scriptsize}
\S \xrightarrow{\hspace{-0.1cm}\tiny\begin{array}{c}
                 \text{$(q+1)$-fold}\\
                  \text{prol. of $u_1$}
                \end{array}\hspace{-0.1cm}
} \S^{(q+1, 0)}\xrightarrow{\hspace{-0.1cm}\tiny\begin{array}{c}
                 \text{$(q+1)$-fold}\\
                  \text{prol. of $u_2$}
                \end{array}\hspace{-0.1cm}} \Xi = \S^{(q+1, q+1)}\xrightarrow{\tiny \t v_1 = \psi(x, u^0, \ldots, u^q, v)} \widetilde \Xi  \xrightarrow{
                \hspace{-0.1cm}\tiny\begin{array}{c}
                 \text{$\ell$-fold}\\
                  \text{prol. of $\t v_1$}
                \end{array}\hspace{-0.1cm}
              } \widetilde\Xi^{(\ell, 0)},
\end{scriptsize}              
$$
which can be further  decomposed into a sequence of successive one-fold prolongations. It follows that flatness implies the LSOP property, and finally, that \cyan{the classes $\mathcal{FLAT}$ and $\mathcal{LSOP}$ coincide}.
%
%
 \end{proof}

\cyan{Relations, given by Theorem~\ref{thm: equiv flat LSOP}, between the \green{four classes of (dynamically) linearizable systems}
can thus be represented by:
\begin{center}
\begin{small}
\begin{tabular}{*{3}{r@{\;}}}
 $\mathcal{LSOPI}$& & \\
  & \rotatebox[origin=c]{-45}{\large$\subset$} & \\
  & & $\mathcal{LSOP}$  $\green = $ $\mathcal{FLAT}$.  \\
    & \rotatebox[origin=c]{45}{\large$\subset$} & \\
 $\mathcal{L}\ell\mathcal{P}$ & &
\end{tabular}
\end{small}
\end{center}
As we have seen in the above proof, in general,
linearization via an~$\ell$-fold prolongation requires an a priori knowledge of the system's flat output in order to determine the control that needs to be prolonged $\ell$ times. Therefore,} when no flat output is known, this method \red{demands} transforming the system into a specific normal form or solving certain PDEs. As a result, the approach is typically non-constructive, motivating the notion of LSOPI introduced in this paper which, as we have seen, is constructive in almost all cases.

\smallskip

 \red{To illustrate the  strict inclusion between  \cyan{$\mathcal{LSOPI}$ and  $\mathcal{LSOP}$},
%
consider the  class of  four-dimensional {control-affine} systems with two inputs, for which the problem of flatness was solved in \cite{pomet1997dynamic}.
The results of \cite{pomet1997dynamic}}
%
%
can be interpreted in terms of dynamic linearizability  via an~$\ell$-fold prolongation of a suitably chosen control.

\begin{example}[LSOP-system, but not {LSOPI}] \label{ex: not necessary conditions}
From the normal forms prensented in~\cite{pomet1997dynamic}, it can be immediately seen that an $x$-flat
{control-affine} system with four states and two inputs
becomes static feedback linearizable after at most two prolongations of a suitable input, while an $(x,u)$-flat system (i.e., all possible flat outputs depend explicitly on~$u$)  is {dynamically} linearizable via a 3-fold prolongation.
\cyan{Therefore all  $x$-flat and $(x,u)$-flat two-input control  systems~$ \S$ in $\mathbb{R}^4$ are} clearly LSOP.  However, unlike the cases where $\ell \leq 2$ (corresponding to $x$-flatness), where~$\S$ can also be shown to be  {LSOPI},  in the case $\ell=3$ (which describes $(x,u)$-flatness)
$\S$ is not {LSOPI}.
In order to see this,
any $(x,u)$-flat two-input control-affine system in dimension four
can be locally transformed via a change of coordinates and an invertible static feedback transformation into the following {form~\cite{pomet1997dynamic}:}
$$
 \S : \left\{\begin{array}{lcl}
\dot x_1 &=&u_1\\
\dot x_2 &=& a(\overline x_3) x_4 + b(\overline x_3) + x_3 u_1\\
\dot x_3 &=& c(\overline x_3) x_4 + d(\overline x_3) + x_4 u_1\\
\dot x_4 &=&u_2,
\end{array}\right.
$$
where  $\overline x_3= (x_1, x_2,x_3)$, $a, b, c, d$ are smooth functions not depending on $x_4$, and~$a$ is not vanishing {at} the nominal point around which we work.
 For the above form, it can be shown that the 3-fold prolongation ${\S}^{(3,0)}$,  obtained by prolonging three times the control $u_1$ as $v_1 = \dot u_1^2 = \frac{\md^3 u_1}{\md t^3}$, see~\eqref{eq: p-fold prolongation}, is static feedback linearizable. It follows that $ {\S}$ is LSOP.
We show next that $ {\S}$ is not {LSOPI}. \red{To distinguish different systems~$\S_i$ constructed by the  algorithm, we use Notation~\ref{notation Sigma i} (for instance,  for the original system $\S_0= \S$, we use the index~0, i.e., we write $\S_0 : \dot x^0 = u_1^0 g_1^0 + u_2^0 g_2^0 $, $x_0\in X^0$, where $x^0 = x$ and $X^0= X$).}

For~$\S_0 = {\S}$, the distribution $\mathcal{D}_0^0$ is not involutive (implying $k_0 = 0$), its {growth} vector is~(2, 3,~4) and $g_2^0 = \partial_{x_4}$ spans the characteristic distribution of $\mathcal{D}_0^0 + [\mathcal{D}_0^0, \mathcal{D}_0^0]$. Thus we are in  Case {III for which ({C2}) holds} and according to our algorithm, the to-be-prolonged control is simply $u_1^{0, \red p} =u_1^{\red 0}$. We obtain the prolongation $\S_1 ={\S_0^{(1,0)}}$  (as in the \red{above proof}, there is no tilde above~$\S_0^{(1,0)}$ because no preliminary feedback is needed) \cyan{controlled by $u_1^1$ and $u_2^1=  u_2$}:
$$
\S_1 : \left\{\begin{array}{lcl}
\dot x_1 &=&u_1^0\\
\dot x_2 &=& a(\overline x_3) x_4 + b(\overline x_3) + x_3 u_1^0\\
\dot x_3 &=& c(\overline x_3) x_4 + d(\overline x_3) + x_4 u_1^0\\
\dot x_4 &=&u_2^1\\
\dot u_1^0 &=&u_1^1.
\end{array}\right.
$$
For~$\S_1$, we  have $\mathcal{D}_1^0 = \mspan\{\partial_{u_1^0}, \partial_{x_4}\}$ and $\mathcal{D}_1^1 =\mathcal{D}_1^0 +  \mspan\{\partial_{x_1} + x_3 \partial_{x_2}+ x_4 \partial_{x_3}, a(\overline x_3) \partial_{x_2} + (c(\overline x_3) +  u_1^0) \partial_{x_3}\}$ which is clearly noninvolutive, so~$k_1 = 1$, and we are in the {subCase~II\,B} with ${r = 1}$ and $[\mathcal{D}_1^0, \mathcal{D}_1^1 ] \not\subset \mathcal{D}_1^1$.
By straightforward computations, it can be shown  that  $\mathcal{D}_1^1$ does not contain any involutive subdistribution of the form $\mathcal{H}_1^1 = \mspan\{\partial_{u_1^0}, \partial_{x_4}, \cyan\b_1(x, u_1^0) (\partial_{x_1} + x_3 \partial_{x_2}+ x_4 \partial_{x_3}) + \cyan\b_2(x, u_1^0)  (a(\overline x_3) \partial_{x_2} + (c(\overline x_3) +  u_1^0) \partial_{x_3})\}$, where $\cyan\b_1, \cyan\b_2$ are smooth functions non vanishing simultaneously. Therefore, according to Proposition~\ref{prop: equiv H feedback}, we cannot construct a prolongation~$\S_2$ for which the first noninvolutive distribution {would be} of order $k_2 > k_1$ implying that~$\S_0$ is not {LSOPI}.

Notice that
when prolonging again~$\S_1$ by adding $\dot u_1^1 =u_1^2$, where $u_1^2$ is the new control, and keeping $u_2^2 = u_2^1$ unchanged, for the resulting prolongation~$\S_2$, we have  $\mathcal{D}_2^1 =\mspan\{\partial_{u_1^1}, \partial_{u_1^0}, \partial_{x_4}, a(\overline x_3) \partial_{x_2} + (c(\overline x_3) +  u_1^0) \partial_{x_3}\}$ which is again noninvolutive (now we are in  {Case~III}:  $r = 1$ and $[\mathcal{D}_2^0, \mathcal{D}_2^1 ]\subset \mathcal{D}_2^1$), so $k_2 =k_1 = 1$. Prolonging again $\dot u_1^2= u_1^3$  and keeping $u_2^3 = u_2^2$ unchanged, we  obtain $\S_3={\S}^{(3,0)}$ which is static feedback linearizable.\hfill $\rhd$
\end{example}

The above example shows that there are systems that are LSOP without being {LSOPI} and
\red{we understand that the obstruction preventing
the LSOPI-algorithm from working on all flat systems, is the requirement of ``involutivity gain``,  translating  by $k_{i+1} \geq k_i$. According to Proposition~\ref{prop: equiv H feedback}, a necessary condition for the ''involutivity gain'' is the existence of an LSOPI-involutive distribution~$\mathcal{H}^k$.  This is a strong structural condition that has to be satisfied for each system $\S_i$, $\dleq 1  i \ell-1$, of the sequence  $\S_0, \S_1, \ldots, \S_\ell$. For LSOP (equivalently, for flatness) this condition is no longer necessary, as shown by Example~\ref{ex: not necessary conditions}.}
\medskip

The following example shows that
\cyan{there exists prolongations for which}
$k_{p}< k$, where $k$ and $k_p$ are, resp.,  the non-involutivity indices of the system~$\cyan \S$ and of its prolongation \cyan{$\widetilde \S^{(1,0)}$}.

\begin{example}[The prolongation may destroy involutivity]\label{ex:involutivity destroyed}
Consider the following control system 
$$
\S: \left\{\begin{array}{lcllcl}
\dot w &= &   x_1^1x_2^1\\
\dot x_1^1 &=&\dot x_1^2 & \dot x_2^1 &=&\dot x_2^2  \\
\dot x_1^2 &=&\dot x_1^3& \dot x_2^2 &=&\dot x_2^3\\
\dot x_1^3 &=& u_1& \dot x_2^3 &=& u_2
           \end{array}\right. 
$$
for which we have $\mathcal{D}^0 = \mspan\{\partial_{x_1^3}, \partial_{x_2^3}\}$,  $\mathcal{D}^1 = \mathcal{D}^0 + \mspan\{\partial_{x_1^2}, \partial_{x_2^2}\}$,  
$\mathcal{D}^2 = \mathcal{D}^1 + \mspan\{\partial_{x_1^1}, \partial_{x_2^1} + x_1^1\partial_w\}$ and $\mathcal{D}^2$ is the first noninvolutive distribution, hence the non-involutivity index of $\S$ is $k= 2$.

Let us now apply the following  invertible feedback transformation
 $u_1 = \t u_1^p + x_1^3 \t u_2$, $u_2 = \t u_2$ and consider the prolongation
 \cyan{$\dot{\t u}_1^p = v_1$,  $\t u_2 = v_2$ yielding}
$$
\widetilde \S^{(1,0)}: \left\{\begin{array}{lcllcl}
\dot w &= &   x_1^1x_2^1\\
\dot x_1^1 &=&\dot x_1^2 & \dot x_2^1 &=&\dot x_2^2  \\
\dot x_1^2 &=&\dot x_1^3& \dot x_2^2 &=&\dot x_2^3\\
\dot x_1^3 &=& \t u_1^p + x_1^3 v_2& \dot x_2^3 &=& v_2\\
\dot{\t u}_1^p  &=&  v_1.
           \end{array}\right. 
$$
For $\widetilde \S^{(1,0)}$, we have $\mathcal{D}_p^0 = \mspan\{\partial_{\t u_1^p}, x_1^3\partial_{x_1^3}+ \partial_{x_2^3}\}$ and $\mathcal{D}_p^1 = \mspan\{\partial_{\t u_1^p},\partial_{x_1^3}, \partial_{x_2^3},$ $x_1^3\partial_{x_1^2}+  \partial_{x_2^2}\}$, which is clearly noninvolutive. It follows that the   \cyan{non-involutivity index} of  $\widetilde \S^{(1,0)}$ is  $k_p= 1 < k=2$. \hfill $\rhd$
\end{example}

Notice however that the above example falls \cyan{into} Case III with condition (C1) satisfied. Therefore
\cyan{$\mathcal{D}^2$ contains infinitely many LSOPI-distributions~$\mathcal{H}^2$ and the prolongation associated to
any~$\mathcal{H}^2$ is static feedback linearizable}. Example~\ref{ex:involutivity destroyed} demonstrates that a ``wrong'' prolongation, i.e., one not associated with a LSOPI-distribution \cyan{$\mathcal{H}^k$}, can actually destroy involutivity.

\color{black}

\section{Example: \red{the chained form}}\label{sec_examples}

%
%
%
%
 \cyan{In this section we show how our results can be applied to}
 systems  static feedback equivalent to the chained form. 
 To simplify the calculations, we will consider the simplest case, that is, systems of the following {chained} form in dimension~4:
 $$
\S: \left\{\begin{array}{lcl}
\dot x_1 &=& u_1\\
\dot x_2 &=& x_3  u_1\\
\dot x_3 &=& x_4 u_1 \\
\dot x_4 &=&u_2,
\end{array}\right.
$$
with \cyan{$x \in X\subset \mathbb{R}^4$},
and suppose that we work around $u_0 = (u_{10}, u_{20})^\top$ such that $u_{10}\neq 0$.
\red{Similarly to the examples of Section~\ref{sec:interplay}, to distinguish different systems~$\S_i$ constructed by the LSOPI algorithm, we use Notation~\ref{notation Sigma i}.}
\cyan{It is obvious that $\S_0=\S$ is LSOPI with $\ell =2$  and that the sequence  of systems  $\S_1, \S_2$ of Definition~\ref{def: LSOPI}  can be constructed  by first, prolonging $u_1^0$ (with no preliminary feedback) yielding $\S_1 = \S_0^{(1,0)}$, with $\dot u_1^0 = u_1^1$ and $u_2^1 = u_2^0$, and second, by prolonging $u_1^1$ (again with no preliminary feedback) to get $\S_2 = \S_2^{(1,0)}$, with $\dot u_1^2 = u_1^1$ and $u_2^2 = u_2^1$, which is static feedback linearizable. Checking the calculations leading to the previous assertion are left to the reader.  The main objective of this section is to illustrate Remark~\ref{rk: to-be-prolonged control 1} (non-{uniqueness} of the to-be-prolonged control) as well as statement a) of Proposition~\ref{prop: H case III} (in case (C1) any LSOPI-distribution~$\mathcal{H}^k$ leads to a static feedback linearizable prolongation). Below we will use the hat symbol for the sequence constructed using general prolongations (involving arbitrary functions).}

The distribution $\mathcal{D}_0^0 = \mspan\{g_1^0, g_2^0 \} =\mspan\{\partial_{x_1} + x_3 \partial_{x_2} + x_4 \partial_{x_3},  \partial_{x_4}\}$ is clearly noninvolutive and its {growth} vector is (2,3,4).
We are in the {Case~III} and {condition ({C2}) of Proposition~\ref{prop: H case III} holds}, {so} we compute the characteristic distribution  $\mathcal{C}(\mathcal{D}_0^0 + [\mathcal{D}_0^0 , \mathcal{D}_0^0 ]) = \mspan\{\partial_{x_4}\}$. Hence {a} to-be-prolonged  control is \cyan{$u_1^0$}.
Recall that  {a}  to-be-prolonged  control is not unique (since the vector field spanning $\mathcal{C}(\mathcal{D}_0^0 + [\mathcal{D}_0^0 , \mathcal{D}_0^0 ])$ is not unique and given up to a multiplicative function), see also  Remark~\ref{rk: to-be-prolonged control 1}.
Let us show that if, instead of prolonging
\cyan{$u_1^0$},
we prolong~$\hat u_1^{\red{0}}$  where $u_1^{0} = \gamma_1(x) \hat u_1^{\red{0}}$, with  $\gamma_1(x)$ any {arbitrary non-zero} function, the proposed algorithm still {detects $\S$ as a LSOPI-system}. The system\footnote{We denote it using the hat symbol since it is~$\hat u_1^{\red{0}}$ that is prolonged.} \cyan{$\widehat \S_1 = \widehat\S_0^{(1,0)}:  \dot x^1 = f^1+ u_1^1g_1^1+ u_2^1g_2^1$} obtained after the first iteration
and for which we drop the hat symbol 
of~$\hat u_1^{\red{0}}$, is given by:
 $$
\widehat \S_1 = \cyan{\widehat\S_0^{(1,0)}}: \left\{\begin{array}{lcl}
\dot x_1 &=& \gamma_1(x)  u_1^0\\
\dot x_2 &=& \gamma_1(x) x_3   u_1^0\\
\dot x_3 &=&  \gamma_1(x) x_4 u_1^0\\
\dot x_4 &=&u_2^1\\
\dot {u}^0_1 &=&u_1^1,
\end{array}\right.
$$
and is considered on $X^1 = X^0 \times \mathbb{R}$, with {the} state $x^{1} = {(x^\top, u_1^0)^\top}$ and {the} control $u^1 = (u_1^1,u_2^1)^{\top}$.
%
For~$\widehat\S_1$, we have
\cyan{$g_1^1 = \partial_{u_1^0}$, $g_2^1 = \partial_{x_4}$, $ad_{f^1} g_1^1 = - \gamma_1 g_1^0$, $ad_{f^1} g_2^1 =  -(\gamma_1  u_1^0 \partial_{x_3} + {\pfrac{\gamma_1}{x_4}} u_1^0 g_1^0) $. Thus
$\mathcal{D}_1^0 = \mspan\{g_1^1, g_2^1 \}$ is, as expected, involutive. On the other hand,
$
\mathcal{D}_1^1 = \mspan\{g_1^1, g_2^1, ad_{f^1} g_1^1,  ad_{f^1} g_2^1\} =
\mspan\{\partial_{u_1^0},  \partial_{x_4},  \partial_{x_3}, {\partial_{x_1} + x_3  \partial_{x_2}}  \}
$}
 is noninvolutive, and such that $r = \mcork(\mathcal{D}_1^1 \subset \mathcal{D}_1^1 + [\mathcal{D}_1^1,\mathcal{D}_1^1]) = 1$ and $ \mathcal{D}_1^1 + [\mathcal{D}_1^1,\mathcal{D}_1^1] = TX^1$ (recall that $g_1^0$ denotes the control vector field of the original system~$\S_0$). {Therefore we are in case ({C1}) of Proposition~\ref{prop: H case III}, statement a)} applies and any  {LSOPI-distribution~$\mathcal{H}_1^1$}
leads to a prolongation that is locally static feedback linearizable.
Indeed, let us take
 \color{black}
$\mathcal{H}_1^1 = \mspan\{\partial_{u_1^0},  \partial_{x_4},  \a_1\partial_{x_3} +  \a_2(\partial_{x_1} + x_3  \partial_{x_2})  \} $,
where~$\a_1$ and~$\a_2$ are smooth functions depending on $x_1, x_2$ and $x_3$ only, not vanishing simultaneously.
\cyan{In order to define  a to-be-prolonged control $\t u_1^{1, p}$ associated to this choice of~$\mathcal{H}_1^1$, recall that we need to compute a vector field $\t g_2^{1, p} = \b_2^1 g_1^1 + \b_2^2 g_2^1$ such that $ad_{f^1} g_2^{1, p} = \a_1\partial_{x_3} +  \a_2(\partial_{x_1} + x_3  \partial_{x_2}) \mmod \mathcal{D}_1^0$, i.e., $\mathcal{H}_1^1 = \mathcal{D}_1^0 +  \mspan\{ad_{f^1} g_2^{1, p}\}$ (see also the comment after the proof of Theorem~\ref{thm: algo}). To that end,} denote $\bar g_1^0= \partial_{x_1} + x_3  \partial_{x_2}$, and notice that
$  \partial_{x_3} = \frac{1}{\gamma_1  u_1^0}[g_2^1 -{\pfrac{\gamma_1}{x_4}}\frac{1}{\gamma_1 } u_1^0 g_1^1, f^1 ] \mmod \mathcal{D}_1^0 $ and
$ \bar g_1^0 = \frac{1}{\gamma_1} [g_1^1, f^1 ] - x_4 \partial_{x_3} $.
Hence
we have $ \a_1\partial_{x_3} +  \a_2 \bar g_1^0 = \frac{1}{\gamma_1 }[(\a_1 - \frac{\a_2 - \a_1 x_4}{\gamma_1}{\pfrac{\gamma_1}{x_4}}) g_1^1     +   \frac{\a_2 - \a_1 x_4}{u_1^0}  g_2^1 , f^1 ] \mmod \mathcal{D}_1^0 $.
We can thus define:
$
\t g_2^{1, \red p} = \b_2^1 g_1^1 + \b_2^2 g_2^1, \mbox{ where } \b_2^1 =  \a_1 - \frac{\a_2 - \a_1 x_4}{\gamma_1}{\pfrac{\gamma_1}{x_4}}, \; \b_2^2 = \frac{\a_2 - \a_1 x_4}{u_1^0},
$ 
and a to-be-prolonged control
$
\t u_1^{1, \red p} = \b_2^2 u_1^1  - \b_2^1 u_2^1.
$
Again, instead of $\t u_1^{1, \red p}$  we prolong~$\hat u_1^{1, \red p}$, where $\t u_1^{1, \red p} = \gamma_2(x, u_1^0) \hat u_1^{1, \red p}$, with  $\gamma_2(x, u_1^0)$ any {non-zero} arbitrary  function. 
Let us suppose that  $\b_2^2(x_0)\neq 0 $, i.e., $(\a_2 - \a_1 x_4) (x_0)\neq 0$ (notice that if $\b_2^2(x_0)= 0 $, then we necessarily have $\b_2^1(x_0) = \a_1 (x_0) \neq 0$, and the rest of the proof is similar).
{This yields $ u_1^1 = \frac{\gamma_2 }{\b_2^2} \hat u_1^{1, \red p}  -  \frac{\b_2^1}{\b_2^2} u_2^1$. Thus,} we get (where again we drop the hat  \red{and the superscript $p$} of $\hat u_1^{1, \red p}$, {and \red{where} $u_2^2 = u_2^1$}):
 $$
\cyan{\widehat \S_2  = \widehat\S_1^{(1,0)}}: \left\{\begin{array}{lcl}
\dot x_1 &=& \gamma_1(x)  u_1^0\\
\dot x_2 &=& \gamma_1(x) x_3   u_1^0\\
\dot x_3 &=&  \gamma_1(x) x_4 u_1^0\\
\dot x_4 &=&u_2^2\\
\dot {u}^0_1 &=&  \frac{\gamma_2(x, u_1^0) }{\b_2^2(x)} u_1^1  -  \frac{\b_2^1(x)}{\b_2^2(x)} u_2^2  \\
\dot {u}^1_1 &=& u_1^2,
\end{array}\right.
$$
which is considered on $X^2 = X^1 \times \mathbb{R} = { X \times \mathbb{R}^2}$. For $\cyan{\widehat \S_2}$, we have:
$\mathcal{D}_2^0 = \mspan\{\partial_{u_1^1}, \partial_{x_4} -   \frac{\b_2^1}{\b_2^2}  \partial_{u_1^0}\}$ and, by a straightforward computation, it can be shown that
\begin{equation*}
 \begin{array}{l@{\,}c@{\,}l}
  \mathcal{D}_2^1 &=&  \mspan\{\partial_{u_1^1},  \partial_{u_1^0}, \partial_{x_4}, \partial_{x_3} + \frac{\a_1}{\a_2 - \a_1 x_4} g_1^0\} 
  \\
  &=& \mspan\{\partial_{u_1^1},  \partial_{u_1^0}, \partial_{x_4}, (\a_2 - \a_1 x_4)\partial_{x_3} + \a_1 g_1^0\}.
 \end{array}
\end{equation*}
Recall that~$\a_1$ and $\a_2$  depend on $x_1, x_2,$ and $x_3$ only; hence $[\partial_{x_4}, (\a_2 - \a_1 x_4)\partial_{x_3} +  \a_1 g_1^0] = - \a_1 \partial_{x_3} + \a_1 \partial_{x_3} = 0$. It follows that $\mathcal{D}_2^1$ is involutive and of constant rank 4.
 Finally,
\begin{equation*}
 \begin{array}{l@{\,}c@{\,}l}
   \mathcal{D}_2^2 &=& \mathcal{D}_2^1 +\mspan\{ \gamma_1  g_1^0,
\gamma_1  u_1^0 \partial_{x_2} +\zeta g_1^0\}= \mspan\{\partial_{u_1^1},  \partial_{u_1^0}, \partial_{x_4}, \partial_{x_3}, \partial_{x_2}, \partial_{x_1}\},
 \end{array}
\end{equation*}
where $\zeta$ is a {function whose explicit form is irrelevant}. All 
$\mathcal{D}_2^j$, \cyan{$\dleq 0 j 2$}, are involutive, of constant rank and $\mathcal{D}_2^2 = TX^2$. Therefore 
$\S_2$ is static feedback linearizable and~$\S$ is {LSOPI}. \hfill $\rhd$

\color{black}

\begin{remark}
The above example illustrates Remark~\ref{rk: to-be-prolonged control 1} and shows that indeed the to-be-prolonged control is not unique and that we can prolong any of them. Notice also that although we {considered}  systems in the {chained}  form in dimension~4, the above procedure generalizes to arbitrary dimension.
\end{remark}


\section{Proofs}\label{sec_proofs}

\subsection{Proof of Proposition~\ref{prop: nec cond flatness}}
\label{ssec: proof prop nec cond flatness}

Consider the flat control system $\S : \dot x = f(x)+u_1g_1(x)+u_2g_2(x)$  at $(x_0,  u_0^l) \in X \times U^l$, for  a certain $l\geq -1$, whose first noninvolutive distribution is~$\mathcal{D}^k$. It is well known that the involutive distributions $\mathcal{D}^{\red j}$, for $ \red j\leq k-1$, are feedback invariant (see \cite{jakubczyk1980on}). The feedback invariance of~$\mathcal{D}^k$ follows from Proposition~7.1 in \cite{nicolau2016two}, \red{proving statement (1)}.

\red{To prove statement (2), observe that, by definition of the sequence of distribution $\mathcal{D}^j$, for two-input control systems
we necessarily have $\mcork (\mathcal D^{k-1} \subset \mathcal D^k)\leq 2$. Thus
only the following three cases are possible:\vspace{-0.1cm}
\begin{list}{-}{}
 \item either $\mcork (\mathcal D^{k-1} \subset \mathcal D^k) = 0$, and in that case we actually have  $ \mathcal D^k = \mathcal D^{k-1}$ and $ \mathcal D^k$ would be involutive, contradicting the hypothesis of Proposition~\ref{prop: nec cond flatness};\vspace{-0.1cm}
 \item
or $\mcork (\mathcal D^{k-1} \subset \mathcal D^k) = 1$;\vspace{-0.1cm}
\item or  $\mcork (\mathcal D^{k-1} \subset \mathcal D^k) = 2$, and in that case statement (2) of Proposition~\ref{prop: nec cond flatness} holds.\vspace{-0.1cm}
\end{list}
}

\red{Therefore let us} assume $\mcork (\mathcal D^{k-1} \subset \mathcal D^k) = 1$ and let $q$ be the smallest integer such that $\mcork (\mathcal D^{q-1} \subset \mathcal D^q) = 1$. It is clear that $\dleq 1 q k$.
For~$\S$, we have the following sequence of inclusions :
$$
\mathcal{D}^0 \underset{2}{\subset} \mathcal{D}^1 \underset{2}{\subset} \cdots  \underset{2}{\subset} \mathcal{D}^{q-1}\underset{1}{\subset} \mathcal{D}^q \underset{1}{\subset} \cdots \underset{1}{\subset} \mathcal{D}^{k-1} \underset{1}{\subset} \mathcal{D}^k \varsubsetneq   TX,
$$
where all distributions except~$\mathcal{D}^k$ are involutive, and the integers appearing under the ``$\subset$''  symbol  denote the corank of the corresponding inclusions.

 \color{black}
For a distribution $\mathcal{D}$ and a function $\psi$, we will write $\md \psi \perp \mathcal{D}$  meaning that $\langle \md \psi , g \rangle = 0 $ for any $g\in \mathcal{D}$.
The distribution~$\mathcal{D}^{k-1}$ is involutive and of constant rank $k+q$, 
so choose independent functions $w = (w_1, \ldots, w_{n-(k+q)})$ such that
$\md w_i \perp \mathcal{D}^{k-1}$, $\dleq 1 i n-(k+q)$.
The distribution $\mathcal{D}^{k-2}$ is also involutive and annihilated by {all differentials} $\md w_i$ and $\md L_f w_i$, so choose~$\psi_1$ of the form {$\psi_1= \md L_f w_1$, 
after a suitable re-ordering of $w_i$'s,} 
such that $\md \psi_1 \perp \mathcal{D}^{k-2}$ and $L_{ad_f^{k-1} g_{\red s}}\psi_1 (x_0)\neq 0$, for either ${\red s}=1$ or ${\red s}=2$ (such $\psi_1$ exists because there exists $ad_f^k g_{\red s} \not\in \mathcal{D}^{k-1}$). Set $z_1^j = L_f^{j-1}\psi_1$, for $\dleq 1 j k$. It follows that $\md z_1^1, \ldots, \md z_1^{k-q+1} \perp \mathcal{D}^{q-2}$ and, since $\mrk  \mathcal{D}^{q-2} = 2 q - 2$, there exists a function~$\psi_2$, independent of \cyan{all} $w_i$ and $z_1^j$, such that $\md \psi_2\perp \mathcal{D}^{q-2}$. Set $z_2^j =  L_f^{j-1}\psi_2$, for $\dleq 1 j q$, and $\t u = L_f^{\rho}\psi + (L_gL_f^{\rho-1}\psi) u$, where $\rho = (k, q)$, {$\psi = (\psi_1, \psi_2)^\top$,  $u = (u_1, u_2)^\top$, $\t u = (\t u_1, \t u_2)^\top$, $g = (g_1, g_2)$}. In \cyan{the} $(w,z)$-coordinates the system reads:
$$
\begin{array}{l@{}c@{}l l@{}c@{}l}
 \dot w& =& h(w, z_1^1)\\
 \dot z_1^j& =&  z_1^{j+1}, \, \dleq 1 j k-1, & \dot z_2^j& =&  z_2^{j+1}, \, \dleq 1 j q-1, \\
  \dot z_1^{k}& =&\t u_1,& \dot z_2^{q}& =&\t u_2.
\end{array}
$$
If $q < k$, then the \cyan{above} form is obvious, and if $q = k$, then the form follows from the equation ${\dot w_1} = z_1^1$ and the condition $\mcork(\mathcal{D}^{k-1} \subset \mathcal{D}^k) = 1$. The above system
consists of two independent subsystems: one with state variables $(w, z_1^j)$, which is nonlinear, and one with state variables $z_2^j$, which is linear. Therefore its flatness (and thus that of~$\S$) is equivalent to the flatness of the single-input $(w,z_1)$-subsystem.
Now recall that, by hypothesis, the distribution
\cyan{$\mathcal{D}^k = \mspan\{\partial_{z_1^j},\partial_{z_2^p}, \sum_{i}{\pfrac{h_i}{z_1^1}\partial_{w_i}}, \dleq 1 j k, \dleq 1 p q\} $}
%
is not involutive implying that the distribution of the same order $k$ associated to the single-input $(w,z_1)$-subsystem and given by \cyan{$\mspan\{\partial_{z_1^j}, \sum_{i}{\pfrac{h_i}{z_1^1}\partial_{w_i}}, \dleq 1 j k\} $} is not involutive either. Since for
single-input control systems, static feedback linearizability and flatness are equivalent, we deduce that the   $(w,z_1)$-subsystem is not flat. Hence, the original two-input system~$\S$ is not flat either, contradicting our assumptions. \red{Given the \cyan{three cases listed at the beginning of the proof of statement (2), we thus deduce that for any flat $\S$ with non-involutivity index $k$}, we necessarily have  $\mcork(\mathcal{D}^{k-1} \subset \mathcal{D}^k) = 2$.}

\color{black}
\subsection{Proof of Proposition~\ref{prop: equiv H feedback}}
\label{ssec: proof prop equiv H feedback}

(i) $\Rightarrow$~(ii). Consider 
$\S : \dot x = f(x)+u_1g_1(x)+u_2g_2(x)$ and suppose that its first noninvolutive distribution~$\mathcal{D}^k$ contains an involutive subdistribution~$\mathcal{H}^k$ verifying   $\mathcal{D}^{k-1} \subset \mathcal{H}^k \subset \mathcal{D}^k$, with both inclusions of corank one. It follows that there exist smooth functions $\b_2^1, \b_2^2$ non vanishing simultaneously such that $\mathcal{H}^k =\mathcal{D}^{k-1} +  \mspan\{\b_2^1 ad_f^k g_1 +  \b_2^2 ad_f^k g_2\}$. We put $\red{\t g_2^p} = \b_2^1 g_1 +  \b_2^2  g_2$ which is clearly non zero and gives $\mathcal{H}^k =\mathcal{D}^{k-1} +  \mspan\{ad_f^k \red{\t g_2^p}\}$.
{Then following the proof of Proposition 7.2 of~\cite{nicolau2016two}, it can be shown that  the involutivity of $\mathcal {H}^k = \mathcal{D}^{k-1} +\mspan \{ad_f^k \red{\t g_2^p}\}$  implies that all 
$\mathcal{H}^j = \mathcal{D}^{j-1}+\mspan \{ad_f^j \red{\t g_2^p}\}$, for $\red 0 \leq j \leq k-1$, are also involutive.}
Now choose any invertible feedback $u =\a+ \b \t u$, \cyan{with $ \t u = (\t u_1^p, \t u_2)$,} such that $\b = \left(
\begin{array}{cc}
* & \b_2^1\\
* & \b_2^2
\end{array}
\right)$ \cyan{and $\a$ arbitrary},
apply it to~$\S$  to get $\widetilde \S : \dot x = {\t f(x)}+\red{\t u_1^p}\t g_1(x)+\t u_2\red{\t g_2^p}(x)$, and prolong the first control to get
$$
\widetilde{\S}^{(1, 0)} :
\left \{
\begin{array}{lcl}
\dot x &= & {\t f(x)}+ \red{\t u_1^p} \tilde g_1(x)+v_2 \red{\t g_2^p}(x) \\
\dot {\t u}_1^{\cyan p}&= &v_1.
\end{array}
\right.
$$
It is immediate {that} the linearizability distributions $\mathcal D_p^{j}$, associated to $\widetilde{\S}^{(1, 0)}$, are of the form
$$\left.
\begin{array}{lcl}
\mathcal D_p^{j} &= & \mspan\{\partial_{\red{\t u_1^p}}\} + \mathcal {H}^j, \; \dleq 0 j k.\end{array}
\right.$$
The involutivity of $\mathcal{H}^j$ implies that of $\mathcal D_p^{j}$ and, in particular, we have obtained the sequence of involutive distributions  $\mathcal{D}_p^0 \subset \ldots \subset \mathcal{D}_p^k$ and $\mrk  \mathcal{D}_p^k = 2k+2$.

(ii) $\Rightarrow$ (i). Consider the control system $\S : \dot x = f(x)+u_1g_1(x)+u_2g_2(x)$ and  let~$\mathcal{D}^k$ be  its first noninvolutive distribution. Suppose that there exists an invertible static feedback transformation $u = {\a(x)} + \beta(x)\tilde u$ such that {for} the prolongation
$$
\widetilde{\S}^{(1, 0)} :
\left \{
\begin{array}{lcl}
\dot x &= & {\t f(x)} + \red{\t u_1^p} \tilde g_1(x)+v_2 \tilde  g_{2}^{\red p}(x) \\
\dot {\t u}_1^{\cyan p}&= &v_1
\end{array}
\right.
$$
{the distributions} $\mathcal{D}_p^0 \subset \ldots \subset \mathcal{D}_p^k$ are involutive and $\mrk  \mathcal{D}_p^k = 2k+2$.
  For simplicity of notation, we will drop the {tildes} \red{as well as the superscript $p$ of $\t u_1^p$ and $\tilde  g_{2}^{\red p}$}.
 If $k=0$ (that is, the first noninvolutive distribution of~$\S$ is  $\mathcal{D}^0 = \mspan \{ g_1, g_2\}$), then the proof is immediate. Let us suppose $k\geq 1$.
  For $\widetilde{\S}^{(1, 0)}$, we have
  $
\mathcal D_p^0 =  \mspan \{\partial_{u_1}, g_2 \}$ and
$
\mathcal D_p^1 =  \mspan \{\partial_{u_1}, g_1, g_2, ad_fg_2+u_1[g_1, g_2]\}.$
Since $k\geq 1$, the distribution $\mathcal{D}^0 = \mspan \{ g_1, g_2\}$ is involutive, thus  $[g_1, g_2] \in  \mathcal{D}^0$ and $\mathcal D_p^1 = \mspan \{\partial_{u_1}, g_1, g_2, ad_fg_2\}.$ It is easy to prove (by an induction argument) that, for $ 1\leq j\leq k$,
$$
\mathcal D_p^j = \mspan \{\partial_{u_1}, g_1,\cdots, ad_f^{j-1}g_1, g_2, \ldots, ad_f^{j}g_2\}.
$$
Since the intersection of involutive distributions is an involutive distribution, it follows that $\mathcal D_p^j \cap~TX = \mspan  \{g_1, \ldots,ad_f^{j-1}g_1,g_2, \ldots, ad_f^{j}g_2\}$ is involutive, for $1\leq { j}\leq k$. We deduce that
$$\mathcal H^k=  \mspan  \{g_1, \ldots,ad_f^{k-1}g_1,g_2, \ldots, ad_f^{k}g_2\}$$
is  involutive. It is immediate that $\mathcal{D}^{k-1} \subset \mathcal {H}^k\subset \mathcal{D}^{k}$, where both inclusions are of corank one, otherwise either $\mathcal {H}^k=\mathcal{D}^{k}$ and~$\mathcal{D}^{k}$ would be involutive  (contradicting the noninvolutivity of~$\mathcal{D}^{k}$) or  $\mathcal {H}^k=\mathcal{D}^{k-1}$ and then $\mrk \mathcal{H}^k$ would be at most $2k$, and since $ \mathcal D_p^k = \mspan \{\partial_{u_1}\} +\mathcal{H}^k$, $\mrk  \mathcal{D}_p^k$ would be at most $2k+1$,   contradicting the assumption 
$\mrk  \mathcal{D}_p^k = 2k+2$.

\subsection{Proof of Proposition~\ref{prop Dk noninvolutive}}
\label{ssec: proof prop Dk noninvolutive}

 {If all distributions $\mathcal{D}^j$, for $\dleq 0 j n$,  are involutive, then $k$ does not exist and the system is in Case~I. If $k$ exists,}
 by applying the Jacobi identity and from the involutivity of the distributions $\mathcal{D}^j$, for $\dleq 0 j k-1$,
 it can be easily shown that  we actually have $[\mathcal{D}^{k-2},\mathcal{D}^k] \subset \mathcal{D}^k$, thus  when calculating $[\mathcal{D}^k, \mathcal{D}^k] $, the new directions that stick out of~$\mathcal{D}^k$  are necessarily obtained with brackets either of the form $[ad_f^{k-1} g_i, ad_f^k g_j]$, for some $\dleq 1 {i,j} 2$, or of the form $[ad_f^{k} g_1, ad_f^k g_2]$. The noninvolutivity of~$\mathcal{D}^k$ implies immediately that we are either in {Case~II }or in {Case~III} (which, {clearly,} are mutually exclusive). {Case~II} may occur only if $k\geq 1$. If $k=0$, then we necessarily are in {Case~III}.
Notice, in particular,  that $r= \mcork(\mathcal{D}^k \subset \mathcal{D}^k + [\mathcal{D}^{k}, \mathcal{D}^k ] ) \leq~4$. Indeed, by applying the Jacobi identity, we {have} $[ad_f^{k-1} g_1,ad_f^{k} g_2] = [ad_f^{k-1} g_2,ad_f^{k} g_1] \mmod \mathcal{D}^k$, so at {most} 3 new directions can be added by vector fields coming from  $ [\mathcal{D}^{k-1}, \mathcal{D}^k ]$.
Let us now suppose that~$\mathcal{D}^k$ admits  a LSOPI-distribution~$\mathcal{H}^k$. According to Proposition~\ref{prop: equiv H feedback},
  there exists $\red{\widetilde {\mathcal{G}}_2^p} = \mspan\{\red{\t g_2^p}\}$ of {corank} one in~$\mathcal{D}^0$, 
such that $ \mathcal{H}^k  =  \mathcal{D}^{k-1} + \mspan\{ad_f^k \red{\t g_2^p}\}$.
Then
$$
\mathcal{D}^k = \mathcal{H}^k  + \mspan\{ad_f^k g_1\},
$$
where $g_1$ is any vector field such that $\mathcal{D}^0=  \mspan\{  g_1,\red{\t g_2^p} \}$.
Repeating the above arguments and from the involutivity of~$ \mathcal{H}^k$, it follows that new  directions that may stick out of~$\mathcal{D}^{k}$ are necessarily given  either by $[ad_f^{k-1} g_1,ad_f^{k} g_1]$ or by $[ad_f^{k} g_1,ad_f^{k} \red{\t g_2^p}]$ or by both, and then we  have $r\leq 2$ and if~$\mathcal{D}^k$ is in {Case~II,} then necessarily ${r_{{}_\mathrm{II}}} = \mcork(\mathcal{D}^k \subset  \mathcal{D}^k +  [\mathcal{D}^{k-1}, \mathcal{D}^k ] ) =1$.

\subsection{Proof of Proposition~\ref{prop: H case II}}
 \label{ssec: proof prop H}

 \color{black}

\textit{Proof of} {subCase~II\,A}. {\it Necessity.}
 Consider the system $\Sigma : \dot x =  f(x)+u_1  g_1(x)+ u_2  g_2(x)$ and let~$\mathcal{D}^k$ be  its first noninvolutive distribution. Assume  that \cyan{Case~II  holds with $r=2$, and that there exists a} LSOPI-distribution~$\mathcal{H}^k$.
Proposition~\ref{prop Dk noninvolutive} asserts that  ${r_{{}_\mathrm{II}}}=1$, proving (A1). Following the same reasoning as in Proposition~\ref{prop: equiv H feedback}, we can define a {non-zero} vector field $\red{\t g_2^p} \in \mathcal{D}^0$ such that $\mathcal{H}^k =\mathcal{D}^{k-1}  + \mspan\{ad_f^k \red{\t g_2^p}\}$.

Denote $\mathcal{E} = \mathcal{D}^{k-1} + [f,\mathcal{C}(\mathcal{D}^k)]$.
We will show next that we actually have $ \mathcal{H}^k = \mathcal{E} = \mathcal{D}^{k-1} + [f,\mathcal{C}(\mathcal{D}^k)]$.
Applying {the} Jacobi identity to $[f, [ad_f^{k-1} g_1,$ $ ad_f^{k-1} {\red{\t g_2^p}}]]$, it can be proved that  $[ad_f^{k-1}\red{\t g_2^p}, ad_f^{k}g_1]\in \mathcal{D}^{k}$, which together with the involutivity of~$\mathcal{H}^k$, gives immediately  $[ad_f^{k-1}\red{\t g_2^p}, \mathcal{D}^{k}] \subset \mathcal{D}^{k}$ and, in particular, that $ad_f^{k-1}\red{\t g_2^p} \in \mathcal{C}(\mathcal{D}^{k})$.
{Recall that $[\mathcal{D}^{k-2},\mathcal{D}^{k}]\subset \mathcal{D}^{k}$. Hence  $\mathcal{D}^{k-2}+ \mspan \{ad_f^{k-1}\red{\t g_2^p}\}\subset  \mathcal{C}(\mathcal{D}^{k})$. We show next that we necessarily have
$
\mathcal{C}(\mathcal{D}^{k}) = \mathcal{D}^{k-2}+ \mspan \{ad_f^{k-1}\red{\t g_2^p}\}.
$
From the involutivity of~$\mathcal{H}^k$ and since
$r=2$ and ${r_{{}_\mathrm{II}}}=1$, it follows that the two new directions completing $\mathcal{D}^{k} $ to $\mathcal{D}^{k} + [\mathcal{D}^{k}, \mathcal{D}^{k}]$ are necessarily given by
$[ad_f^{k-1}g_1, ad_f^{k}g_1]$ and $[ad_f^{k}g_1, ad_f^{k}\red{\t g_2^p}]$.
{Any  $\xi\in \mathcal{C}(\mathcal{D}^{k})$ satisfies $\xi\in \mathcal{D}^{k}$ and is} \cyan{thus}
of the form $\xi = \a_1 ad_f^{k-1} g_1 + \a_2 ad_f^{k} g_1 +\a_3 ad_f^{k} \red{\t g_2^p} + h$, where $h\in \mathcal{D}^{k-2}+ \mspan \{ad_f^{k-1} \red{\t g_2^p}\}$, and $\a_i$ are smooth functions. 
We have $[ad_f^{k-1} g_1, \xi] = \a_2[ad_f^{k-1} g_1,  ad_f^{k} g_1] \mmod \mathcal{D}^k$
and $[ad_f^{k} g_1, \xi] = \a_1[ad_f^{k} g_1,  ad_f^{k-1} g_1]  + \a_3 [ad_f^{k}g_1, ad_f^{k}\red{\t g_2^p}] \mmod \mathcal{D}^k$. Since $\xi$ is a  characteristic vector field, we deduce  $\a_i = 0$, for $i=1,2, 3$. It follows that}
$$
\mathcal{C}(\mathcal{D}^{k}) = \mathcal{D}^{k-2}+ \mspan \{ad_f^{k-1}\red{\t g_2^p}\},
$$
proving (A2),
and that
$$
\mathcal{E} = \mathcal{D}^{k-1} + [f, \mathcal{C}(\mathcal{D}^{k})] = \mathcal{D}^{k-1} +\mspan \{ad_f^{k}\red{\t g_2^p}\}.
$$
Hence $\mrk ( \mathcal{D}^{k-1} +  [f, \mathcal{C}(\mathcal{D}^{k})] )= 2k+1$, showing (A3).
Now recall that the involutive subdistribution~$\mathcal H^k$ is given by
{$\mathcal H^k=  \mathcal{D}^{k-1} +\mspan \{ad_f^{k}\red{\t g_2^p}\}$, so} we actually have $\mathcal H^k=\mathcal E=  \mathcal{D}^{k-1} + [f, \mathcal{C}(\mathcal{D}^{k})]$, implying the involutivity of $ \mathcal{D}^{k-1} + [f, \mathcal{C}(\mathcal{D}^{k})]$ and proving (A4).

\textit{Sufficiency.}
Consider the control system $\Sigma : \dot x = f(x) +u_1  g_1(x)+ u_2  g_2(x)$  whose first noninvolutive distribution~$\mathcal{D}^k$ satisfies {(A1)-(A4)}.
From $r=2$ and ${r_{{}_\mathrm{II}}}=1$, $\mrk \mathcal{C}(\mathcal{D}^k) = 2k-1$, and~$\mathcal{D}^k$ is not involutive, of rank at most $2k+2$, we deduce that necessarily $\mrk \mathcal{D}^k = 2k +2$. %
%
We will prove that conditions (A2)-(A3) enable us to define a {non-zero} vector field $\red{\t g_2^p} \in \mathcal{D}^0$ such that the  distribution $\mathcal {H}^k  =   \mathcal{D}^{k-1} + [f, \mathcal{C}(\mathcal{D}^{k})]$, which is involutive according to (A4), can be written as
$$
\mathcal {H}^k = \mathcal{D}^{k-1}+\mspan \{ad_f^k\red{\t g_2^p}\}.
$$
To define~$\red{\t g_2^p}$, notice first that clearly  $\mathcal{D}^{k-2}\subset \mathcal{C}(\mathcal {D}^k)$. From $[ad_f^k g_1,ad_f^k g_2 ]\not \in \mathcal{D}^k + [ \mathcal{D}^{k-1},  \mathcal{D}^{k}]$, it can be deduced that $\mathcal{C}(\mathcal {D}^k) \subset \mathcal{D}^{k-1}$, {implying that $\mathcal{D}^{k-1} + [f, \mathcal{C}(\mathcal {D}^k)]$ is well defined,} and finally that we have
$$
\mathcal{C}(\mathcal {D}^k) = \mathcal{D}^{k-2} + \mspan \{\zeta\},
$$
with $\zeta$ of the form $\zeta = \a_1 ad_f^{k-1}g_1+ \a_2 ad_f^{k-1}g_2$, where~$\a_1$ and $\a_2$ are smooth functions not vanishing simultaneously. It follows that $\zeta =ad_f^{k-1} (\a_1g_1 + \a_2 g_2) \mmod \mathcal{D}^{k-2}$ and we put $\red{\t g_2^p} = \a_1g_1 + \a_2 g_2$,
which gives
$$
\mathcal{C}(\mathcal {D}^k)  = \mathcal{D}^{k-2} + \mspan \{ad_f^{k-1}\red{\t g_2^p}\},
$$
and, as claimed,
$$
\mathcal {H}^k =\mathcal{D}^{k-1} + [f, \mathcal{C}(\mathcal{D}^{k})]= \mathcal{D}^{k-1}+\mspan \{ad_f^k\red{\t g_2^p}\}.
$$

 \color{black}

\textit{Proof of} {subCase~II\,B}. {\it Necessity.} Let us consider
$\Sigma : \dot x =  f(x)+u_1  g_1(x)+ u_2  g_2(x)$ and let~$\mathcal{D}^k$ be  its first noninvolutive distribution.
\color{black}
Assume  that {Case~II }holds with $r = {r_{{}_\mathrm{II}}} =  1$ (i.e., $ [ \mathcal{D}^{k-1},  \mathcal{D}^{k}]\not\in \mathcal{D}^{k}$ and  $[ad_f^k g_1,ad_f^k g_2 ]\in \mathcal{D}^k + [ \mathcal{D}^{k-1},  \mathcal{D}^{k}]$), and that there exists a LSOPI-distribution~$\mathcal{H}^k$. 
Denote $\mathcal{E} =  \mathcal{C}(\mathcal{D}^k) +\mathcal{D}^{k-1}$. We show next that $ \mathcal{H}^k = \mathcal{E} = \mathcal{C}(\mathcal{D}^k) +\mathcal{D}^{k-1}$.
Following the same reasoning as in the necessity part of the proof of {subCase~II\,A},
define a {non-zero} vector field $\red{\t g_2^p} \in \mathcal{D}^0$ such that $\mathcal{H}^k =\mathcal{D}^{k-1}  + \mspan\{ad_f^k \red{\t g_2^p}\}$,  prove that $ad_f^{k-1}\red{\t g_2^p} \in \mathcal{C}(\mathcal{D}^{k})$, and finally deduce that $[ad_f^{k-1}g_1, ad_f^{k}g_1]\not\in \mathcal{D}^{k}$.
From this together with $r = \mcork(\mathcal{D}^{k} \subset\mathcal{D}^{k}+  [\mathcal{D}^{k}, \mathcal{D}^{k}]) = 1$, it follows that
$[ad_f^{k-1}g_1, ad_f^{k}g_1]$ completes $\mathcal{D}^{k} $ to $\mathcal{D}^{k} + [\mathcal{D}^{k}, \mathcal{D}^{k}]$  implying that there exists a smooth function~$\a$ such that $[ad_f^{k}\red{\t g_2^p}, ad_f^{k}g_1] = \a [ad_f^{k-1}g_1, ad_f^{k}g_1] \mmod {\mathcal D}^k$.
{Hence}
$[ad_f^{k}\red{\t g_2^p}-\a ad_f^{k-1}g_1, ad_f^{k}g_1] =0\mmod {\mathcal D}^k$
{and we deduce:}
\begin{equation}
\label{eq:inclusion}
\mathcal{C}(\mathcal{D}^{k}) \supset \mathcal{D}^{k-2}+ \mspan \{ad_f^{k-1}\red{\t g_2^p} ,ad_f^{k}\red{\t g_2^p}-\a ad_f^{k-1}g_1 \}.
\end{equation}
{It follows that  $\mcork (\mathcal{C}(\mathcal{D}^{k}) \subset \mathcal{D}^k) = 2$, since inclusion of corank one would imply~$\mathcal{D}^k$
involutive, and thus~\eqref{eq:inclusion} is, actually, equality of both distributions, implying $\mrk \mathcal{C}(\mathcal{D}^{k}) = 2k$.}
This proves (B1) and $\mrk (\mathcal{C}(\mathcal{D}^{k})+ \mathcal{D}^{k-1} )= 2k+1$, \cyan{giving} (B2).
Now recall that
$\mathcal H^k= {\mathcal{D}^{k-1}+ \mspan \{ad_f^{k}\red{\t g_2^p}\}}.$
It follows immediately that we
have $\mathcal H^k=\mathcal E= \mathcal{C}(\mathcal{D}^{k})+ \mathcal{D}^{k-1}$, implying the involutivity of $\mathcal{C}(\mathcal{D}^{k})+ \mathcal{D}^{k-1}$ and {thus} proving (B3).

\textit{Sufficiency.}
Consider 
$\Sigma : \dot x = f(x) +u_1  g_1(x)+ u_2  g_2(x)$  whose first noninvolutive distribution~$\mathcal{D}^k$ satisfies {(B1)-(B3)}.
Since $\mrk \mathcal{C}(\mathcal{D}^k) = 2k$ and~$\mathcal{D}^k$ is not involutive of rank at most $2k+2$, we deduce that necessarily $\mrk \mathcal{D}^k = 2k +2$. 
As in the proof of  {subCase~II\,A},  we show that   conditions (B1)-(B2) enable us to define a {non-zero} vector field $\red{\t g_2^p} \in \mathcal{D}^0$ such that the involutive subdistribution $\mathcal {H}^k  = \mathcal{D}^{k-1} + \mathcal{C}(\mathcal{D}^k)$ 
can be written as
$
\mathcal {H}^k = \mathcal{D}^{k-1}+\mspan \{ad_f^k\red{\t g_2^p}\}.
$
We have already noticed that
$\mathcal{D}^{k-2}\subset \mathcal{C}(\mathcal{D}^k)$, see
{the statement at the beginning of the proof of Proposition~\ref{prop Dk noninvolutive},}
and since $\mrk (\mathcal{C}(\mathcal{D}^k) \cap \mathcal{D}^{k-1}) = 2k-1$, we have
$$
\mathcal{C}(\mathcal{D}^k) \cap \mathcal{D}^{k-1} = \mathcal{D}^{k-2} + \mspan \{\zeta\},
$$
with $\zeta$ of the form $\zeta = \a_1 ad_f^{k-1}g_1+ \a_2 ad_f^{k-1}g_2$, where~$\a_1$ and $\a_2$ are smooth functions not vanishing simultaneously. It follows $\zeta =ad_f^{k-1} (\a_1g_1 + \a_2 g_2) $ $\mmod \mathcal{D}^{k-2}$ and we put $\red{\t g_2^p} = \a_1g_1 + \a_2 g_2$. We can always suppose $\a_2(x_0)\neq0$ (otherwise permute $g_1$ and~$g_2$).

Since $ ad_f^{k-1}\red{\t g_2^p} \in \mathcal{C}(\mathcal{D}^k)$, we have $ [ad_f^{k-1}\red{\t g_2^p},\mathcal{D}^{k}] \subset \mathcal{D}^{k}$ and it can be shown, by the involutivity of~$\mathcal{D}^{k-1}$ and applying  the Jacobi identity, that $ [ad_f^{k-1}g_1, $ $ad_f^{k}\red{\t g_2^p}] \in \mathcal{D}^{k}$. Therefore, the new direction completing~$\mathcal{D}^{k}$ to $\mathcal{D}^{k} + [\mathcal{D}^{k}, \mathcal{D}^{k}]$ is given by  $ [ad_f^{k-1}g_1, ad_f^{k}g_1]$ and there exists a smooth function~$\a$ such that $ [ad_f^{k}\red{\t g_2^p},ad_f^{k}g_1] = \a [ad_f^{k-1}g_1, ad_f^{k}g_1]$ $ \mmod \mathcal{D}^{k}$. This gives $ [ad_f^{k}\red{\t g_2^p} - \a ad_f^{k-1}g_1,$ $ad_f^{k}g_1] = 0 \mmod \mathcal{D}^{k}$ and it can be easily verified that
$$
\mathcal{C}(\mathcal{D}^k)  = \mathcal{D}^{k-2} + \mspan \{ad_f^{k-1}\red{\t g_2^p}, ad_f^{k}\red{\t g_2^p} - \a ad_f^{k-1}g_1\},
$$
which gives, as claimed,
$$
\mathcal {H}^k = \mathcal{D}^{k-1}+\mspan \{ad_f^k\red{\t g_2^p}\}.
$$

{
{\it Uniqueness of the  LSOPI-distribution} (subCases II\,A and II\,B).
In the necessity part, we proved that if a LSOPI-distribution  $\mathcal {H}^k$ exists, then it is given by  $\mathcal {H}^k = \mathcal{D}^{k-1} + [f, \mathcal{C}(\mathcal {D}^k)]$ in subCase~II\,A, and by    $\mathcal {H}^k = \mathcal{D}^{k-1} + \mathcal{C}(\mathcal {D}^k)$ in subCase~II\,B, so it is unique (and the second term of each {expression} is obviously feedback invariant).}


 \textit{Proof of} {Case~III}. Consider the control system $\Sigma : \dot x =  f(x)+u_1  g_1(x)+ u_2  g_2(x)$ and suppose that its first noninvolutive distribution~$\mathcal{D}^k$ is such that {Case~III} is satisfied. By Frobenius theorem, {see, e.g.~\cite{jakubczyk1980on}}, all involutive distributions $\mathcal{D}^0\subset \cdots \subset \mathcal{D}^{k-1}$ can be simultaneously rectified and since we necessarily have\footnote{Immediate consequence of~$\mathcal{D}^k$  noninvolutive and $[ \mathcal{D}^{k-1}, \mathcal{D}^k]\subset  \mathcal{D}^k$.}
 $\mrk \mathcal{D}^k = 2k+2$, the system~$\Sigma$ can be brought via a change of coordinates and  a suitable invertible static feedback transformation into
  $$
 \begin{array}{lcl l }
  \dot{\bar w} &=& \overline f (w,  z^1)\\
  {\dot w_{i}} &=& z_i^1\\
  \dot z_i^j &= & z_i^{j+1}, \quad& \dleq 1 j k-1, \\
  \dot z_i^{k}  &= & v_i, \quad &\dleq 1 i 2,
 \end{array}
 $$
where  $\dim w = n- 2k \geq 3$,  {$w = (\bar w^\top, {w_{1}},  {w_{2}})^\top$  and  $z^1 = (z_1^1, z_2^1)^\top$}.
 \cyan{For the above system,} we have
 $$\mathcal{D}^j = \mspan\{\partial_{z_i^{k-j}}, \ldots,\partial_{z_i^{k}}, \dleq 1 i 2\}, \mbox{ for } \dleq 0 j k-1,
 $$
 and
 $$
 \mathcal{D}^k =  \mspan\{\partial_{z_i^{1}}, \ldots,\partial_{z_i^{k}},  {\partial_{w_i}} +{   \pfrac{\overline f } {z_i^{1}} {\partial_{\bar w}}}, \dleq 1 i 2\}.
 $$
{The} condition  $[ \mathcal{D}^{k-1}, \mathcal{D}^k]\subset  \mathcal{D}^k$ implies that  $\overline f(w,  z^1)$ is affine with respect {to~$z_1^{1}$ and~$z_2^{1}$}, that is, we have $\overline f(w,  z^1) = \overline f_0(w ) + z_1^1 \overline f_1(w ) + z_2^1 \overline f_2(w)$ and~$\mathcal{D}^k$ can be written  as
 $$
 \mathcal{D}^k =  \mspan\{\partial_{z_i^{1}}, \ldots,\partial_{z_i^{k}},
 {\partial_{w_i}} +\overline f_i(w) {{\partial_{\bar w}}}, \dleq 1 i 2\}.
 $$
 It is easy to see that~$\mathcal{D}^k$ contains LSOPI-distributions~$\mathcal{H}^k$ {and that they are} never unique.
 Indeed,
 for~$\mathcal{H}^k$, we can take any distribution of the form
  $$
  \begin{array}{l@{\,}c@{\,}l}
 \mathcal{H}^k &=&   \mathcal{D}^{k-1} +  {\mspan\{\xi\}}
 = \mspan\{ \xi, \partial_{z_i^{1}}, \ldots,\partial_{z_i^{k}},
\dleq 1 i 2
 \},
  \end{array}
 $$
 {where $ \xi = \b_1(w)
 \left(\partial_{{w_1}} +\overline f_1(w) {{\partial_{\bar w}}}\right) + \b_2(w)
 \left(\partial_{{w_2}} +\overline f_2(w) {{\partial_{\bar w}}}\right)$, with  $\b_1$ and~$\b_2$  arbitrary smooth functions, depending on  $w$ only, and not vanishing simultaneously (actually we can take either $\b_1 = 1$ or $\b_2 = 1$).}

\subsection{Proof of Proposition~\ref{prop: H case III}}
\label{ssec: proof prop choice H r 1}

Consider the control system $\S : \dot x = f(x)+u_1g_1(x)+u_2g_2(x)$ and suppose that its first noninvolutive distribution~$\mathcal{D}^k$ verifies {Case~III} \cyan{(therefore, according to Proposition~\ref{prop: H case II}, a LSOPI-distribution $\mathcal{H}^k$ always exists, but it is never unique). It is clear that under the Case~III} assumption, the
{system~$\S$ is} in one of the~6 cases {(C1),$\ldots$, (C6) of Proposition~\ref{prop: H case III} and that they are mutually exclusive.}
{We show next the} {existence of~$\t g_2$ claimed in} condition~({C5}).
%
%

 {\it Proof of } {({C5})}.
We show that  the conditions  $\mrk \overline {\mathcal{D}}^k= 2k+3$ and
$\mrk (\overline {\mathcal{D}}^k+ [f,\mathcal{D}^k]) = 2k+4$ imply the existence of a {non-zero} vector field $ {\t g_2} \in \mathcal{D}^0$ such that $ad_f^{k+1} {\t g_2} \in \overline {\mathcal D}^k$.
We can always assume (permute~$g_1$ and~$g_2$, if necessary) that $ad_f^{k+1}g_1 \not\in \overline {\mathcal{D}}^k$. Hence there exists a smooth function~$\alpha$ such that  $ad_f^{k+1}g_2 = \alpha ad_f^{k+1}g_1 \mmod \overline {\mathcal{D}}^k$. It follows  that  
$ad_f^{k+1}(g_2  - \alpha g_1) = 0\mmod \overline {\mathcal{D}}^k$. The vector field ${\t g_2} = g_2  - \alpha g_1$ is clearly {non-zero} (since~$g_1$ and~$g_2$ are independent everywhere on $X$) and satisfies $ad_f^{k+1}{\t g_2} \in \overline {\mathcal D}^k$.

 {\it Proof of} {a)}. {The proof of {statement {a)}}
 is contained in the proof of Theorem 3.4 of \cite{nicolau2016two} stating that if {Case~III} holds and, moreover, 
 $\mathcal{D}^k +[\mathcal{D}^{k}, \mathcal{D}^k] = TX$, then 
 $\S$ is actually  dynamically linearizable via a one-fold prolongation
 without any structural condition (only a regularity condition is needed).}

 {\it Proof of}~b). Consider 
 $\S : \dot x = f(x)+u_1g_1(x)+u_2g_2(x)$
 and suppose that its first noninvolutive distribution~$\mathcal{D}^k$ verifies {Case~III}  and {(C2)}.
 {We show first that 
 $\mathcal{C}(\mathcal{D}^k + [\mathcal{D}^{k}, \mathcal{D}^k])$ is indeed a LSOPI-distribution.}
  We introduce the notation $\mathcal{E}^0 =  \mathcal{D}^k $ and $\mathcal{E}^{j+1} = \mathcal{E}^j + [\mathcal{E}^j, \mathcal{E}^j]$ for $j = 0, 1$.
 Recall that {Case~III}  corresponds to $[\mathcal{D}^{k-1}, \mathcal{D}^k]\subset \mathcal{D}^k$ {and by} successive applications of the Jacobi identity, it follows immediately that we also have  $[\mathcal{D}^{k-1}, \mathcal{E}^1]\subset \mathcal{E}^1$, thus  $\mathcal{D}^{k-1} \subset \mathcal{C}( \mathcal{E}^1)$, where $\mathcal{C}( \mathcal{E}^1)$ is the characteristic distribution of~$\mathcal{E}^1$.
 Since $\mrk \mathcal{E}^2 = 2k+4$, we can always assume without loos of generality (permute~$g_1$ and~$g_2$, if necessary) that  $[ad_f^{k} g_1, [ad_f^{k} g_1, ad_f^{k} g_2]] \not\in \mathcal{E}^1$.
 Hence  there exists a smooth function $\alpha$ such that  $[ad_f^{k}g_2, [ad_f^{k} g_1, ad_f^{k} g_2]]  = \alpha [ad_f^{k} g_1, [ad_f^{k} g_1, ad_f^{k} g_2]] \mmod \mathcal{E}^1$.
 It follows  that
 $[ad_f^{k}(g_2  - \alpha g_1), [ad_f^{k} g_1, ad_f^{k} g_2]] = 0\mmod \mathcal{E}^1$.
 The vector field $\t g_2 = g_2  - \alpha g_1$ is clearly {non-zero} (since~$g_1$ and~$g_2$ are independent everywhere on $X$) and satisfies $[ad_f^{k}\t g_2, \mathcal{E}^1] \subset \mathcal{E}^1$.
 Hence $ ad_f^{k}\t g_2 \in  \mathcal{C}( \mathcal{E}^1)$ and we deduce that
  $$
  \mathcal{C}( \mathcal{E}^1)   \supset \mathcal{D}^{k-1} + \mspan\{ad_f^k \t g_{2}\}
  $$
  {and the inclusion above is, actually, the equality of both distributions, since~$ \mathcal{E}^1$ is non involutive and $ \mrk \mathcal{E}^1 = 2k+3$.} 
Recall that $ \mathcal{E}^1 = \mathcal{D}^k + [\mathcal{D}^{k}, \mathcal{D}^k]$,  so  we have just shown that $\mathcal{D}^{k-1} \subset \mathcal{C}(\mathcal{D}^k + [\mathcal{D}^{k}, \mathcal{D}^k]) \subset \mathcal{D}^k$, with both inclusions of corank one, {and since the characteristic distribution is {always} involutive, it follows that  $ \mathcal{C}(\mathcal{D}^k + [\mathcal{D}^{k}, \mathcal{D}^k])$ is, indeed, a LSOPI-distribution.} 

 We  prove {next} that any\footnote{\cyan{There are many LSOPI-distributions $\mathcal{H}^k$, see Proposition~\ref{prop: H case II}.}} LSOPI-distribution $\mathcal{H}^k \neq \mathcal{C}(\mathcal{D}^k + [\mathcal{D}^{k}, \mathcal{D}^k])$  leads to a prolongation for which  $\mathcal{D}_p^{k+1}$ is noninvolutive and  does not admit a LSOPI-distribution $\mathcal{H}_p^{k+1}$.
    Let us suppose $k\geq 1$. If $k=0$ (i.e., the first noninvolutive distribution of~$\S$ is  $\mathcal{D}^0 = \mspan \{ g_1, g_2\}$), then the proof is similar and left to the reader.
  Transform~$\Sigma$ via a  static feedback {transformation} into the form $\widetilde{\Sigma} : \dot x = f(x) +\tilde u_1\tilde   g_1(x)+\tilde  u_2 \tilde g_2(x)$, where $\tilde g_2$ is such that $ad_f^k \t g_{2} \in \mathcal{C}(\mathcal{D}^k + [\mathcal{D}^{k}, \mathcal{D}^k])$ and $\tilde g_1$ is any {non-zero} vector field {satisfying} $\mathcal{D}^
  0= \mspan\{\t g_1, \t g_2\}$, and {such that}
  $\mathcal{H}^k = \mathcal{D}^{k-1} + \mspan\{ad_f^k \t g_{1}\}$ is a LSOPI-distribution. Construct the prolongation associated to~$\mathcal{H}^k$:
   \begin{equation}\label{eq: S p proof C-i}
\widetilde{\S}^{(0, 1)} :
\left \{
\begin{array}{lcl}
\dot x &= &f(x)+ v_1 \tilde g_1(x)+ \t u_2 \tilde  g_{2}(x) \\
\dot {\t u}_2&= &v_2,
\end{array}
\right.
 \end{equation}
  where $v_1 = \t u_1$ (notice the notation $\widetilde{\S}^{(0, 1)}$ indicating that~$\t u_2$ is prolonged, while  $\t u_1$ is kept unchanged).
For simplicity of notation, we will drop the {tildes}.
Following the proof of Proposition~\ref{prop: equiv H feedback}, 
for $ 0\leq j\leq k$, we have
   \begin{equation}\label{eq: Dp j proof C-i}
   \begin{array}{l@{\,}c@{\,}l}
    \mathcal D_p^j &=& \mspan \{\partial_{u_2}, g_1,\ldots, ad_f^{j}g_1, g_2,\ldots, ad_f^{j-1}g_2\}
    = \mspan \{\partial_{u_2}\} +\mathcal{H}^j,
   \end{array}
 \end{equation}
which are involutive (recall that according to Proposition~\ref{prop: equiv H feedback}, all distributions {$\mathcal{H}^j = \mathcal{D}^{j-1} + \mspan \{ ad_f^j g_1\}$, for $ 0\leq j\leq k-1$, associated to the LSOPI-distribution $\mathcal{H}^k = \mathcal{D}^{k-1} + \mspan \{ ad_f^k g_1\}$}, are also involutive).
Hence
   \begin{equation}\label{eq: Dp k+1 proof C-i}
\mathcal D_p^{k+1} = \mspan \{\partial_{u_2}, g_1,\ldots, ad_f^{k+1}g_1, g_2,\ldots, ad_f^{k}g_2\}.
 \end{equation}
If $ad_f^{k+1}g_1 \not\in \mathcal{E}^1$, where 
{$ \mathcal{E}^1 = \mathcal{D}^k + [\mathcal{D}^{k}, \mathcal{D}^k]$,} 
then it is clear that $[ad_f^kg_1, ad_f^kg_2]\not\in \mathcal D_p^{k+1}$ implying that~$\mathcal D_p^{k+1}$ is not involutive.
If $ad_f^{k+1}g_1  \in \mathcal{E}^1$, then at least the bracket $[ad_f^kg_1,[ad_f^kg_1, ad_f^kg_2]]\not\in \mathcal D_p^{k+1}$ (recall that $ ad_f^kg_2 \in \mathcal{C}(\mathcal{E}^1)$),  yielding again the noninvolutivity of~$\mathcal D_p^{k+1}$. It follows that in both cases $k_p = k+1$ and, moreover, $[ \mathcal D_p^{k}, \mathcal D_p^{k+1}]\not\in \mathcal D_p^{k+1}$, so if~$\mathcal D_p^{k+1}$ admits a LSOPI-distribution $ \mathcal H_p^{k+1}$, then $ \mathcal H_p^{k+1}$ is unique, {see Proposition~\ref{prop: H case II}}.
Suppose that $ \mathcal H_p^{k+1}$ exists. Then $ \mathcal H_p^{k+1}$ is of the form
\begin{equation}\label{eq: Hp k+1 general}
\begin{array}{l@{\,}l}
\mathcal H_p^{k+1} =  \mspan \{&\partial_{u_2}, g_1,\ldots, ad_f^{k}g_1, g_2,\ldots, ad_f^{k-1}g_2,
\a_1 ad_f^{k+1}g_1 + \a_2ad_f^{k}g_2 \},
\end{array}
\end{equation}
where~$\a_1$ and $\a_2$ are smooth functions (of $x$ and $u_2$) not vanishing simultaneously.
{Suppose that there are points $(x_0, u_0)$ such that}
  $\a_1(x_0, u_{20}) \neq 0$, then
 we actually have
\begin{equation}\label{eq: case C-i Hp k+1}
\begin{array}{l@{\,}l}
\mathcal H_p^{k+1} =  \mspan \{&\partial_{u_2}, g_1,\ldots, ad_f^{k}g_1, g_2,\ldots, ad_f^{k-1}g_2,
ad_f^{k+1}g_1 + \a  ad_f^{k}g_2 \},
\end{array}
\end{equation}
where $\a $ is a smooth function.
Now recall that for the system~$\S$, we have  $[ \mathcal{D}^{k-1}, \mathcal{D}^k]\subset \mathcal{D}^k$,  thus  $
[ad_f^{k-1} g_2, ad_f^kg_2] \in \mathcal{D}^k
$
and there exists a smooth function $\eta_1$ such that
$$
[ad_f^{k-1} g_2, ad_f^kg_2] = \eta_1  ad_f^kg_2 \mmod \mathcal{H}^k,
$$
where 
$\mathcal{H}^k = \mathcal{D}^{k-1} + \mspan\{ad_f^k  g_{1}\}$. 
By applying the Jacobi identity and the fact that $[ad_f^{k-1} g_2, ad_f^kg_1] \in \mathcal{H}^k$, we deduce
$$
\begin{array}{lll}
[ad_f^{k-1} g_2, ad_f^{k+1}g_1] &=& [ad_f^{k} g_1, ad_f^kg_2] + [f, [ad_f^{k-1} g_2, ad_f^{k}g_1]]
\\
&= &[ad_f^{k} g_1, ad_f^kg_2] + \eta_2 ad_f^kg_2 + \eta_3  ad_f^{k+1}g_1 \mmod \mathcal{H}^k,
 \end{array}
$$ 
for    smooth functions $\eta_2$ and $\eta_3$.
Using the above relations, compute
\begin{small}
$$
\begin{array}{l}
[ad_f^{k-1} g_2,  ad_f^{k+1}g_1 + \a ad_f^{k}g_2] = \\
\hspace{0.9cm}=  [ad_f^{k-1} g_2,  ad_f^{k+1}g_1 ] + \a [ad_f^{k-1} g_2, ad_f^{k}g_2] + (L_{ad_f^{k-1} g_2} \a)ad_f^{k}g_2  \\
\hspace{0.9cm}= [ad_f^{k} g_1, ad_f^kg_2] + \gamma_1 ad_f^{k+1}g_1 +   \gamma_2 ad_f^kg_2  \mmod \mathcal{H}^k,
\end{array}
$$
\end{small}
for   suitable smooth functions $\gamma_1$ and $\gamma_2$ (that can be expressed in {terms} of~$\eta_i$ and~$\a$).
Suppose first that  $ad_f^{k+1}g_1 \not\in \mathcal{E}^1$.
Then $[ad_f^{k} g_1, ad_f^kg_2] \wedge ad_f^{k+1}g_1 \wedge ad_f^kg_2 \neq 0\mmod \mathcal{H}^k$ and from $\mathcal H_p^{k+1}  = \mathcal{H}^k +  \mspan \{\partial_{u_2},ad_f^{k+1}g_1 + \a  ad_f^{k}g_2 \}$, we deduce that $[ad_f^{k-2} g_2,  ad_f^{k+1}g_1 + \a ad_f^{k}g_2]\not\in \mathcal H_p^{k+1}$. Therefore, $\mathcal H_p^{k+1}$ is not involutive,
{and thus not a} LSOPI-distribution.

Let us now assume that $ad_f^{k+1}g_1 \in \mathcal{E}^1$. It follows that
$ad_f^{k+1}g_1$ can be written as
\begin{equation}
 \label{eq:adf k+1 g1}
 ad_f^{k+1}g_1 = \eta_4  [ad_f^{k} g_1, ad_f^kg_2] + \eta_5  ad_f^kg_2 \mmod \mathcal{H}^k,
\end{equation}
hence
$$
\begin{array}{l@{\,}l}
\mathcal H_p^{k+1} =  \mspan \{ & \partial_{u_2}, g_1,\ldots, ad_f^{k}g_1, g_2,\ldots, ad_f^{k-1}g_2, 
\eta_4  [ad_f^{k} g_1, ad_f^kg_2] + { \t  \eta_5}  ad_f^{k}g_2 \}, 
\end{array}
$$
{where  $\t  \eta_5 = \eta_5 +\a$.} 
Compute
$ \zeta =  [ad_f^{k} g_1,  \eta_4  [ad_f^{k} g_1, ad_f^kg_2] + { \t  \eta_5}  ad_f^{k}g_2] = \eta_4  [ad_f^{k} g_1,$ $ [ad_f^{k} g_1, ad_f^kg_2]] + (L_{ad_f^{k} g_1}  \eta_4 +  { \t \eta_5}) [ad_f^{k} g_1, ad_f^kg_2] + (L_{ad_f^{k}g_1} { \t  \eta_5})  ad_f^{k}g_2$. Then $\zeta \in~\mathcal H_p^{k+1}$ implies that we necessarily have $\eta_4 =0$ (recall that
$ad_f^{k}g_2$ is characteristic for $\mathcal{E}^1 =\mathcal{D}^k +[\mathcal{D}^{k}, \mathcal{D}^k] $ and $[ad_f^{k} g_1, [ad_f^{k} g_1, ad_f^kg_2]]\not\in \mathcal{E}^1$), and finally that
$$
\begin{array}{l@{\,}c@{\,}l}
\mathcal H_p^{k+1} &= & \mspan \{\partial_{u_2}, g_1,\ldots, ad_f^{k}g_1, g_2, \ldots, ad_f^{k-1}g_2,   ad_f^{k}g_2 \}
= \mspan \{\partial_{u_2} \} + \mathcal{D}^k,
\end{array}
$$
which is not involutive, contradicting the fact that $\mathcal H_p^{k+1}$ is a LSOPI-dis\-tri\-bu\-tion.
{It follows that there are no points  $(x_0, u_0)$ such that
  $\a_1(x_0, u_{20}) \neq 0$ and {thus~$\a_1$} is actually identically zero. In that case  we  have   $\mathcal H_p^{k+1} =  \mspan \{\partial_{u_2}\} + \mathcal{D}^k$, and $\mathcal H_p^{k+1}$ is not involutive, {and thus not an LSOPI-distribution}.}

We have just shown that in all possible cases, any LSOPI-distribution $\mathcal{H}^k \neq \mathcal{C}(\mathcal{D}^k + [\mathcal{D}^{k}, \mathcal{D}^k])$  leads to a prolongation for which  $\mathcal{D}_p^{k+1}$ is noninvolutive and  does not admit a LSOPI-distribution $\mathcal{H}_p^{k+1}$. {Therefore, the only  LSOPI-distribution~$\mathcal{H}^k$
defining the prolongation $\widetilde \S^{(1,0)}$ that may possess a LSOPI-distribution $\mathcal{H}_p^{k_p}$ is  $\mathcal{H}^k = \mathcal{C}(\mathcal{D}^k + [\mathcal{D}^{k}, \mathcal{D}^k])$.}

 {\it {Proof of}}  {c) for  (C5) with (C5)$'$}.
 {Assume first that~$\S$ verifies (C5),
 and transform~$\Sigma$} via a  static feedback {transformation} into $\widetilde{\Sigma} : \dot x = f(x) +\tilde u_1\tilde   g_1(x)+\tilde  u_2 \tilde g_2(x)$, where $\tilde g_2$ is defined by {({C5})} and $\tilde g_1$ is any {non-zero} vector field {satisfying $\mathcal{D}^
  0= \mspan\{\t g_1, \t g_2\}$ and} {a LSOPI-distribution is either 
  $\mathcal{H}^k = \mathcal{D}^{k-1} + \mspan\{ad_f^k \t g_{1}\}$ or $\mathcal{H}^k = \mathcal{D}^{k-1} + \mspan\{ad_f^k \t g_{2}\}$.}
  
  {Let us first suppose that $\mathcal{H}^k = \mathcal{D}^{k-1} + \mspan\{ad_f^k \t g_{1}\}$  
  is  LSOPI.} To simplify the notation, drop the tildes.
  As in the proof of 
  {statement {b)}} for the case {({C2})}, the prolongation $\widetilde{\S}^{(0, 1)}$ associated to~$\mathcal{H}^k$
 is of the {form~\eqref{eq: S p proof C-i}}, the distributions $\mathcal D_p^j$, for $\dleq 0 j k$, are given by~\eqref{eq: Dp j proof C-i} and are all involutive, while~$\mathcal D_p^{k+1}$ is given by~\eqref{eq: Dp k+1 proof C-i}.

From {({C5})}, we necessarily have $ad_f^{k+1}g_1 \not\in \overline{\mathcal{D}}^k$, then it is clear that $[ad_f^kg_1, ad_f^kg_2]\not\in \mathcal D_p^{k+1}$ implying that~$\mathcal D_p^{k+1}$ is not involutive.
Suppose that a LSOPI-distribution  $ \mathcal H_p^{k+1}$ exists. Then $ \mathcal H_p^{k+1}$ is of the form
\begin{equation}\label{eq: Hp k+1 general (C5)}
\begin{array}{l@{\,}l}
\mathcal H_p^{k+1} =  \mspan \{&\partial_{u_2}, g_1,\ldots, ad_f^{k}g_1, g_2,\ldots, ad_f^{k-1}g_2,
\a_1 ad_f^{k+1}g_1 + \a_2ad_f^{k}g_2 \},
\end{array}
\end{equation}
where~$\a_1$ and $\a_2$ are smooth functions (of $x$ and $u_2$) not vanishing simultaneously, with~$\a_1$ non identically zero.
Following the proof of 
{statement {b)}}, it can be shown that
\begin{equation}\label{eq: zeta C-iii}
\begin{array}{l@{\;}c@{\;}l}
\zeta_1&= &[ad_f^{k} g_1,  \a_1 ad_f^{k+1}g_1 + \a_2 ad_f^{k}g_2]
=
\a_1 [ad_f^{k} g_1,   ad_f^{k+1}g_1 ] + \a_2[ad_f^{k} g_1, ad_f^kg_2] 
\\
& &\qquad + (L_{ad_f^{k} g_1}\a_1) ad_f^{k+1}g_1 +  (L_{ad_f^{k} g_1}\a_2)  ad_f^kg_2,
\\
\zeta_2 &=& [ad_f^{k-1} g_2,  \a_1 ad_f^{k+1}g_1 + \a_2 ad_f^{k}g_2]
=\a_1[ad_f^{k} g_1, ad_f^kg_2] + \gamma_1 ad_f^{k+1}g_1 \\
&& \qquad+   \gamma_2 ad_f^kg_2+ \gamma_3 ad_f^{k+1} g_2  \mmod \mathcal{H}^k.
\end{array}
\end{equation}
{Suppose {now} that} {\nf{({C5})}$'$ holds,} {that is,} $ ad_f^{k+1}g_2 \in \mathcal{D}^k$ (recall that~$g_2$ is the special vector field~$\t g_2$ defined by {({C5})}, for which we dropped the tilde), 
then
from\footnote{
\red{Recall that~$\a_1$ is non identically zero, see the comment just after relation~\eqref{eq: Hp k+1 general (C5)}.}}
$ \a_1\neq 0$ and $ ad_f^{k+1}g_1\wedge [ad_f^{k} g_1, ad_f^kg_2]\neq 0 \mmod \mathcal{D}^k$, it follows that $\zeta_2 \not\in \mathcal H_p^{k+1} $, implying that $\mathcal H_p^{k+1}$ is not involutive, contradicting the  LSOPI-distribution assumption. 

{Let us now assume {that}~$\S$ satisfies (C5) with (C5)$'$ and  that $ {\mathcal{E}}= \mathcal{D}^{k-1} + \mspan\{ad_f^k \t g_{2}\}$, with  $\t g_{2}$ defined by (C5), is involutive, implying that ${\mathcal{E}}$ is a LSOPI-distribution. {Define $\mathcal{H}^k =\mathcal{E}$ and} construct its associated prolongation $\widetilde{\S}^{( {1, 0})}$. Using similar arguments as before, it can be shown that for   $\widetilde{\S}^{( {1, 0})}$, we have $\mathcal{D}_p^k =  \mspan \{\partial_{u_1}\} + \mathcal{H}^k$ and $\mathcal{D}_p^{k+1} = \mathcal{D}_p^k + \mspan \{ad_f^k g_{1}, ad_f^{k+1} \t g_{2}\} =\mspan \{\partial_{u_1}\} + \mathcal{D}^k$ which is noninvolutive and such that $\mcork (\mathcal{D}_p^{k}\subset \mathcal{D}_p^{k+1}) =1$. By Proposition~\ref{prop: nec cond flatness}, it follows that~$\widetilde{\S}^{( {1, 0})}$ is not flat
{and therefore}~$\S$ is not LSOPI.
}

{{\it Proof of}  c) for  (C6)}. Suppose that~$\S$ satisfies ({C6}), so we have $\mrk \overline {\mathcal{D}}^k= 2k+3$ and
$\mrk (\overline {\mathcal{D}}^k+ [f,\mathcal{D}^k]) = 2k+5$, and suppose that $k\geq 1$ (like for the case {({C2})}, if $k=0$, then similar arguments apply and   the proof is left to the reader). Consider any LSOPI-distribution~$\mathcal{H}^k$ and vector fields $g_1, g_2$ such that
  $\mathcal{D}^
  0= \mspan\{g_1, g_2\}$ and
  $\mathcal{H}^k = \mathcal{D}^{k-1} + \mspan\{ad_f^k g_{2}\}$.
  Construct the prolongation associated to~$\mathcal{H}^k$ and given by
   \begin{equation}
{\S}^{(1,0)} :
\left \{
\begin{array}{lcl}
\dot x &= &f(x)+ u_1  g_1(x)+ v_2   g_{2}(x) \\
\dot {u}_1&= &v_1
\end{array}
\right.
 \end{equation}
for which the distributions
$
\mathcal D_p^j = \mspan \{\partial_{u_1}, g_1,\ldots,  ad_f^{j-1}g_1,$ $ g_2,\ldots, ad_f^{j}g_2\} =  \mspan \{\partial_{u_2}\} +\mathcal{H}^j,
$
 for $ 0\leq j\leq k$, are involutive, and
\begin{equation}
 \label{eq: Dp noninv C6}
\mathcal D_p^{k+1} = \mspan \{\partial_{u_1}, g_1,\ldots, ad_f^{k}g_1, g_2,\ldots, ad_f^{k+1}g_2\}.
\end{equation}
According to the conditions describing case {({C6})}, we have  $[ad_f^{k} g_1, ad_f^kg_2] \wedge ad_f^{k+1}g_1 \wedge ad_f^{k+1}g_2 \neq 0\mmod \mathcal{D}^k$, implying that~$\mathcal D_p^{k+1}$ is not involutive.
Suppose that  a LSOPI-distribution $ \mathcal H_p^{k+1}$ exists. Then $ \mathcal H_p^{k+1}$ is of the form
\begin{equation}\label{eq: Hp k+1 general C-iv}
\begin{array}{l@{\,}l}
\mathcal H_p^{k+1} =  \mspan \{& \partial_{u_1}, g_1,\ldots, ad_f^{k-1}g_1, g_2,\ldots, ad_f^{k}g_2,
\a_1 ad_f^{k}g_1 + \a_2ad_f^{k+1}g_2 \},
\end{array}
\end{equation}
where~$\a_1$ and $\a_2$ are smooth functions (of $x$ and $u_1$) not vanishing simultaneously.
%
%
{Suppose that there are points  $(x_0, u_0)$ such that  $\a_2(x_0, u_{10}) \neq 0$. Then} we actually have
$$
\begin{array}{l@{\,}l}
\mathcal H_p^{k+1} =  \mspan \{&\partial_{u_1}, g_1,\ldots, ad_f^{k- 1}g_1, g_2,\ldots, ad_f^{ k}g_2, 
\a ad_f^{ k}g_1 +  ad_f^{k + 1}g_2 \},
\end{array}
$$
where $\a $ is a smooth function.

Similarly to {the proof of~b)}, we have
$
[ad_f^{k-1} g_1, ad_f^kg_1] = \eta_1  ad_f^kg_1 \mmod \mathcal{H}^k
$
({indeed,} recall that $[\mathcal{D}^{k-1}, \mathcal{D}^k]\subset \mathcal{D}^k$)
and
$ 
[ad_f^{k-1} g_1, ad_f^{k+1}g_2] = - [ad_f^{k} g_1, ad_f^kg_2] + \eta_2 ad_f^kg_1 + \eta_3  ad_f^{k+1}g_2 \mmod \mathcal{H}^k,
$
for    smooth functions $\eta_1, \eta_2$ and $\eta_3$. From this, we get
\begin{equation}\label{eq: bracket noninv C-iv}
\begin{array}{l}
[ad_f^{k-1} g_1,  \a ad_f^{k}g_1 + ad_f^{k+1}g_2] 
=
-[ad_f^{k} g_1, ad_f^kg_2]
+ \gamma_1 ad_f^{k}g_1 +   \gamma_2 ad_f^{k+1}g_2  \mmod \mathcal{H}^k,
\end{array}
\end{equation}
for  \cyan{some}  smooth functions $\gamma_1, \gamma_2$.  Recall that {by (C6),} $[ad_f^{k} g_1, ad_f^kg_2]\wedge ad_f^{k}g_1 \wedge  ad_f^{k+1}g_1 \wedge ad_f^{k+1}g_2 \neq 0\mmod \mathcal{H}^k$. Hence, $[ad_f^{k-1} g_1,  \a ad_f^{k}g_1 +  ad_f^{k+1}g_2]\not\in \mathcal H_p^{k+1}$. Therefore, $\mathcal H_p^{k+1}$ is not involutive,  {and thus} $\mathcal H_p^{k+1}$ is not a LSOPI-distribution.
  {It follows that $\a_2$ is identically zero, but then  $\mathcal H_p^{k+1} =  \mspan \{\partial_{u_2}\} + \mathcal{D}^k$, and $\mathcal H_p^{k+1}$ is not involutive {either, giving a} contradiction.}

{{\it Proof of}  c) for  (C3)}.
{Let us now assume that~$\S$ satisfies ({C3})}
and repeat the proof of \red{the} case {({C6})} until {obtaining~\eqref{eq: Dp noninv C6}. Then according to the condition describing (C3), we have  $[ad_f^{k} g_1, ad_f^kg_2] \wedge [ ad_f^{k}g_1, [ad_f^{k} g_1, ad_f^kg_2]] \wedge  [ ad_f^{k}2, [ad_f^{k} g_1, ad_f^kg_2]] \neq 0\mmod \mathcal{D}^k$, implying that~$\mathcal D_p^{k+1}$ is not involutive.
Suppose that  a LSOPI-distribution $ \mathcal H_p^{k+1}$ exists and repeat the proof of  case  ({C6}) from~\eqref{eq: Hp k+1 general C-iv} until~\eqref{eq: bracket noninv C-iv} (notice, in particular,  that for obtaining~\eqref{eq: bracket noninv C-iv} only the involutivity of~$\mathcal{H}^k$ and the fact that~$\mathcal{D}^k$ satisfies Case~III were used).}
We know that  $[ad_f^{k} g_1, ad_f^kg_2]  \wedge ad_f^{k}g_1  \neq 0\mmod \mathcal{H}^k$, and that $ad_f^{k+1} g_2\not\in  \mathcal{H}^k$. It follows that the only possibility for $\mathcal{H}_p^{k+1}$ to be involutive is that $ad_f^{k+1} g_2 \in \mathcal{D}^k+[\mathcal{D}^k,\mathcal{D}^k]$, implying
$$
ad_f^{k+1}g_2 = \eta_4  [ad_f^{k} g_1, ad_f^kg_2] + \eta_5  ad_f^kg_1 \mmod \mathcal{H}^k,
$$
{which is~\eqref{eq:adf k+1 g1}, with the role of $g_1$ and $g_2$ permuted.}
Now recall that the {{growth}} vector of~$\mathcal{D}^k$ starts with $(2k+2, 2k+3,2k+ 5, \ldots)$, thus 
{$\zeta$ defined by the same formula as is the proof of b),  with   $g_1$ and $g_2$ permuted, satisfies
$\zeta \in \mathcal H_p^{k+1}$, which implies that we necessarily have $\eta_4 =0$.  Finally, we get}
$$
\begin{array}{l@{\,}c@{\,}l}
 \mathcal H_p^{k+1} &=&  \mspan \{\partial_{u_1}, g_1,\ldots, ad_f^{k-1}g_1, g_2,\ldots, ad_f^{k}g_2,   ad_f^{k}g_1 \} 
= \mspan \{\partial_{u_1} \} + \mathcal{D}^k,
\end{array}
$$
which is not involutive, contradicting  the fact that $\mathcal H_p^{k+1}$ is a LSOPI-dis\-tri\-bu\-tion.   {Like for case (C6), we conclude that $\a_2$ is identically zero, but then  $\mathcal H_p^{k+1} =  \mspan \{\partial_{u_2}\} + \mathcal{D}^k$, and $\mathcal H_p^{k+1}$ is not involutive, {yielding a} contradiction.}
\smallskip

To sum up, we have shown that if~$\S$ is such that its first noninvolutive  distribution~$\mathcal{D}^k$ satisfies  {({C3})} or {({C6})}, then  any LSOPI-distribution~$\mathcal{H}^k$  leads to a prolongation for which  $\mathcal{D}_p^{k+1}$ is noninvolutive and  does not admit a LSOPI-distribution $\mathcal{H}_p^{k+1}$. Therefore,~$\S$ is not LSOPI.

{{\it Proof of}  c) for  (C4)}.
{Suppose now {that}~$\S$ satisfies ({C4}), that is, we have} $\mathcal{D}^k +[\mathcal{D}^{k}, \mathcal{D}^k] \neq TX$ and $\mrk \overline {\mathcal{D}}^k= \mrk (\overline {\mathcal{D}}^k+ [f,\mathcal{D}^k]) = 2k+3$. The system is {thus} not strongly accessible (see~\cite{sussmann1972controllability} for a formal definition and an algebraic {criterion}), thus not flat either (recall that flat systems are always strongly accesible, see \cite{fliess61vine},  and also, e.g., \cite{aranda1994infinitesimal,pomet1997dynamic}) {and, in particular, not LSOPI,} showing statement  c) for case ({C4}).

\appendix

\section{Construction of {a special~$\mathcal{H}^k$} in case (C5)}

\color{black}

\begin{proposition}
\label{prop: case C5 non unique}
Consider 
$\S: \dot x = f(x)+u_1g_1(x)+u_2g_2(x)$ given by~\eqref{eq: Sigma}, and assume that its first non involutive distribution~$\mathcal{D}^k$  is such that $\mathcal{D}^k +[\mathcal{D}^{k}, \mathcal{D}^k] \neq TX$. 
There exists  a prolongation ${\widetilde \S^{(1,0)}} $ satisfying the conditions of Proposition~\ref{prop: equiv H feedback}{\normalfont(ii)} with $\mathcal{D}_p^{k+1}$ also involutive and such that {$\mcork (\mathcal{D}_p^{k+1} \subset\mathcal{D}_p^{k+2} )= 2$} 
if and only if~$\mathcal{D}^k$
satisfies \normalfont{{({C5})}} and
\begin{enumerate}[\normalfont {({C5}}a{)}]
 \item   {$\mrk (\overline {\mathcal D}^k + [f, \overline{\mathcal D}^k])= 2k+5$,}
\end{enumerate}
and 
if $k\geq 1$ (and only in that case),~$\mathcal{D}^k$ additionally verifies {condition} {(C5)$''$ and}:
\begin{enumerate}[\normalfont {({C5}}a{)}]
\setcounter{enumi}{1}
   \item 
   The distribution
 $\mathcal{D}^{k-1} + \mspan\{ad_f^k \t g_{2}\},$
with $\t g_{2}$ defined {by~{\normalfont (C5)}}, is involutive and in particular, it is a LSOPI-distribution.
\end{enumerate}
\end{proposition}

\begin{remark}
If $\mathcal{D}^k +[\mathcal{D}^{k}, \mathcal{D}^k]  \neq TX$, {it follows from}
the above result
that for all {cases (C1)-(C6)} of Proposition~\ref{prop: H case III}  (except {for} that described  by Proposition~\ref{prop: case C5 non unique}) no matter the choice of~$\mathcal{H}^k$ defining 
$\widetilde{\S}^{(1,0)}$, the distribution~$\mathcal{D}_p^{k+1}$ would (already!) be  noninvolutive and only one involutive distribution is gained for the prolongation\footnote{Recall that 
${\widetilde \S^{(1,0)}}$ {is} such that its associated distribution $\mathcal{D}_p^{k} $ is given by $\mathcal{D}_p^{k} = \mspan\{\partial_{\t u_1}\} +  \mathcal{H}^k $ and {thus} involutive due to the involutivity of~$ \mathcal{H}^k $.}.
{Condition ({C5}) with ({C5})$''$,} {(C5a) and (C5b),} gives necessary and sufficient conditions for the existence of a {unique} involutive subdistribution~$\mathcal{H}^k$, {which we call special,}  leading to a prolonged system ${\widetilde \S^{(1,0)}}$ {
for which~$\mathcal{D}_p^{k+1}$ is also} involutive  and explicitly constructs that  involutive subdistribution.
\end{remark}

\begin{proof}

 {\it Necessity.} Consider the control system $\S : \dot x = f(x)+u_1g_1(x)+u_2g_2(x)$,  
 {and assume that its first non involutive distribution~$\mathcal{D}^k$  is such that $\mathcal{D}^k +[\mathcal{D}^{k}, \mathcal{D}^k] \neq TX$.} Suppose that there exists an invertible static feedback transformation\footnote{{We use the hat symbol, instead of tilde,  to distinguish its associated vector field $\hat g_2$ from~$\t g_2$ defined by condition (C5).}}
 $u = \a(x) +  \beta(x)\hat u$  such that the {associated} prolongation
 \begin{equation}\label{eq: S p proof C-i hat}
\widehat{\S}^{(1, 0)} :
\left \{
\begin{array}{lcl}
\dot x &= &\hat f(x)+ \hat u_1 \hat g_1(x)+v_2 \hat  g_{2}(x) \\
\dot {\hat u}_1&= &v_1
\end{array}
\right.
 \end{equation}
is such that  all distributions $\mathcal{D}_p^0 \subset \cdots \subset \mathcal{D}_p^k \subset \mathcal{D}_p^{k+1}$ are involutive and 
{$\mcork (\mathcal{D}_p^{k+1} \subset\mathcal{D}_p^{k+2} )= 2$}. 
For simplicity of notation, we will drop the {hat symbol}.
We first suppose that $k\geq 1$. Repeating  the proof of the implication~(ii) $\Rightarrow$ (i) of Proposition~\ref{prop: equiv H feedback}, we deduce that the distribution 
$$
\mathcal H^k= \mathcal D_p^k \cap~TX   = \mspan  \{g_1, \ldots,ad_f^{k-1}g_1,g_2,\ldots, ad_f^{k}g_2\}
$$
is  involutive and such that both inclusions  $\mathcal{D}^{k-1} \subset \mathcal {H}^k\subset \mathcal{D}^{k}$ are of corank one, i.e., $\mathcal H^k$ is a LSOPI-distribution.
The involutivity of
$$
\begin{array}{l@{\,}l}
\mathcal D_p^{k+1} = \mspan \{\partial_{u_1}, g_1, \ldots,ad_f^{k}g_1, g_2, \ldots, ad_f^{k}g_2,
ad_f^{k+1}g_2\}
\end{array}
$$
implies that of $\mathcal{D}^k + \mspan \{ad_f^{k+1}g_2\}$. This yields that we necessarily have $\overline {\mathcal{D}}^k = \mathcal{D}^k + \mspan \{ad_f^{k+1}g_2\}$ and $\mrk \overline {\mathcal{D}}^k = 2k+3$ (in particular, that $ad_f^{k+1}g_2\not\in\mathcal{D}^k$ {and that $[ad_f^kg_1, ad_f^kg_2] = \a ad_f^{k+1}g_2 \mmod \mathcal{D}^k$, with~$\a$ a non vanishing smooth function).
Compute $\mathcal{D}_p^{k+2} =  \mathcal D_p^{k+1} + \mspan \{ ad_f^{k+1}g_1,  ad_f^{k+2}g_2\}= \mspan \{\partial_{u_1}\} +\overline {\mathcal{D}}^{k}+[f,\overline{\mathcal{D}}^{k}]$. Then from $\mcork (\mathcal{D}_p^{k+1} \subset\mathcal{D}_p^{k+2} )= 2$, it follows immediately that $\mrk {\mathcal D}^{k+1}= 2k+4$,  $\mrk (\overline {\mathcal D}^k + [f, {\mathcal D}^k])= 2k+4$  and 
 $\mrk (\overline {\mathcal D}^k + [f,\overline {\mathcal D}^k])= 2k+5$.} 
{To sum up, we have} $ad_f^{k+1}g_1 \not\in \overline {\mathcal D}^k$ and   $\overline {\mathcal{D}}^k = \mathcal{D}^k + \mspan \{ad_f^{k+1}g_2\}$, where $g_1$ and~$g_2$ are the vector fields $\hat g_1$ and $\hat g_2$ used to define the prolongation $\widehat{\S}^{(1, 0)}$ in~\eqref{eq: S p proof C-i hat}, for which we dropped the {hat symbol}. Therefore $\t g_2 = \gamma  g_2$ (where  $ \t g_2 \in \mathcal{D}^0$ is the {non-zero} vector field such that $ad_f^{k+1} \t g_2 \in \overline {\mathcal D}^k$, {defined by ({C5})}, with $\gamma$ a non vanishing function, and it follows that
$$
\mathcal{D}^{k-1} + \mspan\{ad_f^k \t g_2\} =  \mathcal{D}^{k-1} + \mspan\{ad_f^k  g_2\} = \mathcal{H}^k,
$$
is involutive. 

Let us now suppose that $k=0$, i.e., the distribution $\mathcal{D}^0 = \mspan \{ g_1, g_2\}$ is noninvolutive implying  $[g_1, g_2] \not\in~\mathcal{D}^0$. Since~$\mathcal D_p^1 = \mspan \{\partial_{u_1}, g_1, g_2, ad_fg_2+u_1[g_1, g_2]\}$ is assumed involutive, it follows that we actually have
$\mathcal D_p^1 = \mspan \{\partial_{u_1}, g_1, g_2, [g_1, g_2], ad_fg_2\} = \mspan \{\partial_{u_1}, g_1, g_2, [g_1, g_2]\} $
yielding the involutivity of $\mathcal{D}^0+[\mathcal{D}^0,\mathcal{D}^0] = \mspan \{g_1, g_2, [g_1, g_2]\}$. So $\overline{\mathcal{D}}^0 = \mathcal{D}^0+[\mathcal{D}^0,\mathcal{D}^0]$, and is of rank 3. 
Moreover, we have $ad_fg_2 \in \overline{\mathcal{D}}^0$, and thus a non zero vector field ${\t g_2} \in \mathcal{D}^0$ such that $ad_f{\t g_2} \in  \overline{\mathcal{D}}^0$, whose existence is claimed in {condition ({C5})}, can be taken as ${\t g_2} = g_2$.
Notice that condition {({C5})$''$} $ad_fg_2 \not\in {\mathcal{D}}^0$ is not necessary if $k=0$, and (C5b) is {trivially satisfied} (indeed, we have $\mspan\{ {\t g_2}\} = \mspan\{g_2\}$ involutive).

 {\it Sufficiency.}
Consider the control system $\Sigma : \dot x = f(x) +u_1  g_1(x)+ u_2  g_2(x)$ with~$\mathcal{D}^k$ satisfying 
{(C5) with ({C5})$''$} and {(C5a)-(C5b)}.  Transform~$\Sigma$ via a  static feedback {transformation} into the form $\widetilde{\Sigma} : \dot x = \t f(x) +\tilde u_1\tilde   g_1(x)+\tilde  u_2 \tilde g_2(x)$, where $\tilde g_2$ is defined in the statement of {({C5})}. By Proposition~\ref{prop: equiv H feedback}, the involutivity of  $\mathcal{H}^j= \mathcal{D}^{j-1}+\mspan \{ad_f^j \t g_2\}$,  for $1\leq j \leq k-1$, follows from that of $\mathcal{H}^k= \mathcal{D}^{k-1}+\mspan \{ad_f^k \t g_2\}$ and it is immediate to see that {for {the} prolongation
associated to~$\mathcal{H}^k$, we have $ \mathcal D_p^j = \mspan \{\partial_{u_1}\} + \mathcal{H}^j$, for $\dleq 0 j k$, which are involutive. {Now}
 $\mrk (\overline {\mathcal{D}}^k+ [f,\mathcal{D}^k]) = 2k+4$, together with  ({C5})$''$, implies that $\overline {\mathcal D}^k = {\mathcal D}^k + \mspan \{ad_f^{k+1} \t g_2\}$. It follows that
$
\mathcal D_p^{k+1} = \mspan \{\partial_{u_1}\} + \overline {\mathcal D}^k
$
 is also involutive.} {Finally, (C5a) yields $\mcork (\mathcal{D}_p^{k+1} \subset\mathcal{D}_p^{k+2} )= 2$.}
\end{proof}

\color{black}

\section{Subcases of {(C5)$''$ {for which} the algorithm is conclusive}}

\begin{proposition}\label{prop: C-iii}
Consider 
$\S: \dot x = f(x)+u_1g_1(x)+u_2g_2(x)$ {and} let~$\mathcal{D}^k$ be its first noninvolutive distribution. Assume that $k\geq 1$ {and that~$\mathcal{D}^k$ satisfies  {\normalfont {({C5})}} with (C5)$''$} 
and let~$\t g_2$ be the vector field defined by {\normalfont {({C5})}}.  Then the distribution $\mathcal{D}^{k+1}$ is feedback invariant and {$   \mathfrak{r} =  \mcork (\mathcal{D}^{k+1} \subset \mathcal{D}^{k+1}+ [\mathcal{D}^k, \mathcal{D}^{k+1}]) \leq 2.$
Assume that either  $\mathfrak{r} = 2$ or   $\mathfrak{r} = 1$  
and, in {the latter} case, assume additionally $[ad_f^{k} \t g_{2}, \mathcal{D}^{k+1}] \subset  \mathcal{D}^{k+1}$.
Then if~$\S$ is LSOPI}, {the distribution $\mathcal{E} =  \mathcal{D}^{k-1} + \mspan\{ad_f^k \t g_{2}\}$, where $\t g_{2}$ satisfies $ad_f^{k+1} \t g_{2} \in  \overline {\mathcal{D}}^k$, has to be involutive and
 the prolongation has to be defined by   $\mathcal{H}^k =\mathcal{E}$ (which implies,  by Proposition~\ref{prop: case C5 non unique},  that $\mathcal{D}_p^{k+1}$ is also involutive).} %

\end{proposition}

\begin{center}
\begin{figure}
 \includegraphics[scale=0.65]{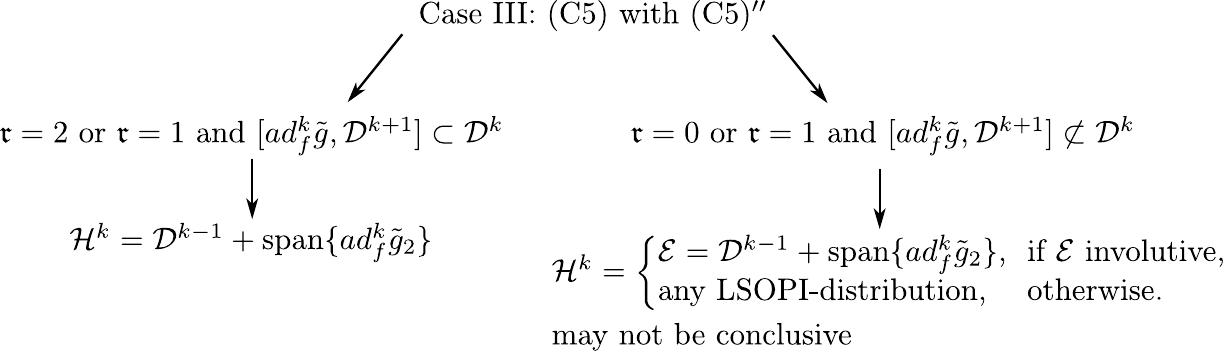}
 \caption{Construction of~$\mathcal{H}^k$ by Proposition~\ref{prop: C-iii}.}
 \label{fig:caseIII annex}
\end{figure}
\end{center}

\begin{remark}
 Suppose that 
 $\mathcal{D}^{k-1} + \mspan\{ad_f^{k} \t g_{2}\}$, where~$\t g_2$ is defined by  {({C5})}, is involutive (i.e., it is a LSOPI-distribution).
 Then whenever one of the cases of Proposition~\ref{prop: C-iii} is satisfied,  the only  approach to construct a sequence of systems $\S_0, \ldots, \S_i$ that may terminate at $\S_\ell$, which is static feedback linearizable, is to define the LSOPI-distribution as $\mathcal{H}^k = \mathcal{D}^{k-1} +  \mspan\{ad_f^{k} \t g_{2}\}$. If $\mathcal{D}^{k-1} + \mspan\{ad_f^{k} \t g_{2}\}$ is not involutive, and we are in one of the above cases, then the system is not LSOPI.
 {Finally, 
 the only cases where the algorithm is not conclusive are those when (C5) with (C5)$''$ holds, and either $\mathfrak{r} = 1$ and   {$[ad_f^{k} \t g_{2}, \mathcal{D}^{k+1}] \not\in  \mathcal{D}^{k+1}$} or $\mathfrak{r} = 0$, \cyan{see Figure~\ref{fig:caseIII annex}}.}
\end{remark}

\begin{proof}
{The invariance of $\mathcal{D}^{k+1}$ with respect to transformations $g \mapsto \b g $
and $f\mapsto  f + g\alpha$,  and \cyan{the bound $\mathfrak{r} \leq 2$ are straightforward and their proofs are} left to the reader.}
{Now repeat the proof of Proposition~\ref{prop: H case III}, statement c) for (C5) with (C5)$'$, until relations~\eqref{eq: zeta C-iii}, and notice that to obtain~\eqref{eq: zeta C-iii} only condition (C5) was used. Suppose now that {\normalfont (C5)$''$} holds and that} either  
{$\mathfrak{r} = 2$ or   $\mathfrak{r} = 1$}
and  $[ad_f^{k+1} g_{2}, \mathcal{D}^{k+1}]\in  \mathcal{D}^{k+1}$, then it can be shown that
  $
  ad_f^{k}g_2\wedge ad_f^{k+1}g_1\wedge  ad_f^{k+1}g_2\wedge [ad_f^{k} g_1, ad_f^{k+1}g_1] \neq 0 \mmod \mathcal{H}^k,
  $ 
 from which we deduce that
$\zeta_1 \not\in \mathcal H_p^{k+1} $, {see the first relation of~\eqref{eq: zeta C-iii}},
{contradicting the involutivity of $\mathcal H_p^{k+1}$.}
 We have shown that {when the above conditions hold}, any LSOPI-distribution~$\mathcal{H}^k$,  such that  $\mathcal{H}^k \neq \mathcal{D}^{k-1} + \mspan\{ad_f^{k} g_{2}\}$, leads to a prolongation for which  $\mathcal{D}_p^{k+1}$ is noninvolutive and  does not admit a LSOPI-distribution $\mathcal{H}_p^{k+1}$. {Thus if~$\S$ is LSOPI, then we necessarily have to define~$\mathcal{H}^k $ by $\mathcal{D}^{k-1} + \mspan\{ad_f^{k} g_{2}\}$}.
\end{proof}


\bibliographystyle{plain}
\bibliography{bibliographie}

\end{document}